\documentclass{article}
\usepackage[top=3cm,bottom=3cm,left=2.5cm,right=2.5cm]{geometry}
\usepackage[utf8]{inputenc}
\usepackage[english]{babel}
\usepackage{mathtools}
\usepackage[all,ps]{xy}
\usepackage{amsthm}
\usepackage{amsfonts}
\usepackage{mathrsfs}
\usepackage{amssymb}
\usepackage{enumitem}
\usepackage{pdfpages}
\usepackage[a-1b]{pdfx}
\usepackage[pdfa]{hyperref}
\usepackage[numbered]{bookmark}

\newtheorem{theorem}[subsubsection]{Theorem}
\newtheorem{corollary}[subsubsection]{Corollary}
\newtheorem{lemma}[subsubsection]{Lemma}
\newtheorem{proposition}[subsubsection]{Proposition}

\theoremstyle{definition}
\newtheorem{definition}{Definition}[subsection]

\theoremstyle{remark}

\def\lim{\varinjlim} %direct limit

\makeatletter
\newcommand*{\doublerightarrow}[2]{\mathrel{
  \settowidth{\@tempdima}{$\scriptstyle#1$}
  \settowidth{\@tempdimb}{$\scriptstyle#2$}
  \ifdim\@tempdimb>\@tempdima \@tempdima=\@tempdimb\fi
  \mathop{\vcenter{
    \offinterlineskip\ialign{\hbox to\dimexpr\@tempdima+1em{##}\cr
    \rightarrowfill\cr\noalign{\kern.5ex}
    \rightarrowfill\cr}}}\limits^{\!#1}_{\!#2}}}
\makeatother

\CompileMatrices

\title{A Proof of the Countable Telescope Conjecture for Module Categories}
\author{Pier Federico Pacchiarotti, advised by Prof.\ Jan Trlifaj}
\date{December 2021}

\begin{document}

\maketitle
\thispagestyle{empty}

%\clearpage
\pagenumbering{arabic}
%\newpage
\begin{abstract}
    The Countable Telescope Conjecture arose in the framework of stable homotopy theory, as a tool conceived to study the chromatic filtration. It turned out, however, to trigger extremely fertile research within the framework of Module Categories.
    The project aims at presenting an almost self-contained review of the recent work of Saroch on the Countable Telescope Conjecture for Module Categories. After recalling some preliminaries, we report various devices of independent interest that will lead to a proof of the aforementioned result. This will be the outcome of inductive refinements of families of particularly well-behaved dense systems of modules, our witnessing-notion for localness. The procedure will be reminiscent of Cantor diagonal argument in the implementation of a variant of Shelah's Compactness Principle. Then, we briefly review the main applications to Enochs Conjecture of the just developed theory, and we will also state a weaker version of it.
    
    The project closely follows the work of Saroch, \cite{Saroch-tools}; however, for the sake of completeness and conciseness, we slightly modified some well-known proofs applying the newly developed tools.
\end{abstract}

\newpage

\tableofcontents

\newpage
\addcontentsline{toc}{section}{Introduction}

\section*{Introduction}

In 1984 the Countable Telescope Conjecture was first officially announced by Ravenel as a statement about the so called chromatic filtration, in the field of stable homotopy theory. It was subsequently generalized to the framework of compactly generated triangulated categories by Neeman, under the name of Generalized Smashing Conjecture.

Unfortunately, it turned out to be false in such generality, as Keller exhibited a counterexample in 1994.

However, it was still unknown whether the conjecture could still be valid in some particular setting of interest, as for instance in the stable module category of a self-injective Artin algebra. So, in 2003 this led Krause and Solberg to formulate a new version of it, closer to the one about which we will discuss in this project.

In particular, they considered it in the context of cotorsion pairs of modules. The latter is a fundamental concept in representation theory of modules, which can be regarded as an attempt to introduce some homological duality within the category of $R$-modules, that is well-known to be not self-dual at the categorical level.
\smallskip

To do so, one considers those pairs of  $Ext^1_R$-orthogonal classes of $R$-modules, say $(\mathcal{A},\mathcal{B})$, such that $\mathcal{A}^{\perp}=\mathcal{B}$ as well as $^{\perp}\mathcal{B}= \mathcal{A}$ w.r.t. $Ext^{1}_R$, or equivalently $\mathcal{A}^{\perp}=\mathcal{B}$ and $^{\perp}{(\mathcal{A}^{\perp})}=\mathcal{A}$. Remark that the latter viewpoint amounts precisely to a 'Galois connection' supported by a proper class: $(\_)^{\perp}: (Mod-R,\subseteq) \rightarrow (Mod-R,\subseteq)$ with adjoint $^{\perp}(\_)$.

So far, we obtained just an analogue of the classic notion of torsion theories for the bifunctor $Ext^1$. However, the importance of the latter formalism is due to the following inherent homological duality relative to left-right approximations of modules. Classical references for the well-known theory of cotorsion pairs are for instance \cite{Enochs}, \cite{Xu}, \cite{Trlifaj&Goebel}.

In representation theory, we call a right (resp. left) approximation of an $R$-module $M$ with some class $\mathcal{C}$ - or $\mathcal{C}$-preenvelope (resp. $\mathcal{C}$-precover) - any $R$-linear morphism with domain $M$ and codomain in $\mathcal{C}$ which factorises any other morphism $M\rightarrow C$ for each $C\in \mathcal{C}$ (resp. require the arrow-dual property). We say that such an approximation is special whenever it is further a monomorphism with cokernel in $^{\perp}\mathcal{C}$ (resp. an epimorphism with kernel in $\mathcal{C}^{\perp}$), properties that need not hold true for a generic preenvelope (resp. precover).

For a fixed choice of $M, \mathcal{C}$ as before, one can consider the category of all special $\mathcal{C}$-preenvelopes of $M$ endowed with the natural subobject order relation induced by inclusion of the codomains. Whenever it exists, the minimum for such a 'lattice' (supported by a proper class) is called $\mathcal{C}$-envelope of $M$, and can be characterized by the fact that any morphism $M\rightarrow \mathcal{C}$ factors uniquely through it.

The following result by Salce clarifies the intrinsic duality property of such approximations. Namely, given any cotorsion pair $(\mathcal{A},\mathcal{B})$, each module has a special $\mathcal{A}$-precover iff each module has a special $\mathcal{B}$-preenvelope; in this case, we say that the cotorsion pair is complete.

This is a frequent occurrence in interesting cases, since a variation of Quillen's small object argument allows to prove that each cotorsion pair of the form $({^{\perp}(\mathcal{S}^{\perp})},\mathcal{S}^{\perp})$ for some set $\mathcal{S}$ is indeed complete, so that one has an abundance of complete cotorsion pairs.

Moreover, the pairs of this kind turn out to be of great interest from the viewpoint of representation theory. Indeed, Enochs proved that for a class of modules $\mathcal{C}$ that is closed under extensions and arbitrary direct limits, provided that a module has a special $\mathcal{C}^{\perp}$-preenvelope, then it also has a $\mathcal{C}^{\perp}$-envelope. An analogous result holds for (pre)covers, although this time we can weaken the assumptions on $\mathcal{C}$, since, due to the properties of precovers, we do not need the closure under extensions. 

A particular application of the aforementioned machinery is the Flat Cover Conjecture, which allows to compute $Ext$ via minimal flat resolutions.
\smallbreak

In our framework, then, the Telescope Conjecture as formulated by Krause and Solberg represents an attempt to determine the 'dimension' of the right class in some particularly interesting cotorsion pairs.

More explicitly, they aimed at investigating whether for a complete hereditary cotorsion pair $(\mathcal{A},\mathcal{B})$ over an Artin algebra $R$ it could hold $\mathcal{A} = \lim (\mathcal{A} \cap mod-R)$, provided that both classes are closed under $\lim$. (Here $mod-R$ is the class of strongly finitely presented $R$-modules.)

This was proved true for some particular cases, among which noteworthy is that of tilting cotorsion pairs (see \cite{Trlifaj&Goebel}).

Notwithstanding the great importance of such a class, it was remarkable that no further specialization of the considered underlying ring was required.

Subsequently, then, the Countable Telescope Conjecture for Module Categories acquired its classical formulation, with no assumptions on the ring $R$, as stated in 2008 in the paper \cite{classical ctc} by Angeleri H\"ugel, Saroch and Trlifaj.
\smallskip

\smallskip

More precisely, let $\mathcal{C}=(\mathcal{A},\mathcal{B})$ be a hereditary cotorsion pair. The classical Telescope Conjecture for Module Categories (TCM) states that any complete hereditary cotorsion pair whose classes are closed under direct limits is of finite type, meaning that it is generated by a set of strongly finitely presented modules.
\bigskip

In the aforementioned paper, the conjecture was proved for the relevant instance of cotorsion pairs over a right-Noetherian ring whose right class has bounded injective dimension, comprising in particular cotilting cotorsion pairs.

However, the advancements towards the general statement are again only limited to special cases satisfying some further finiteness assumptions, whereas in \cite{Saroch&Stovicek} Saroch and Stovicek proved a weaker and arguably more natural version of it, namely the Countable Telescope Conjecture for Module Categories (CTCM).

It states that any hereditary cotorsion pair is of countable type whenever the right class is closed under unions of well-ordered chains.

Differently from the rest of the literature, here a cotorsion pair of countable type is intended to be one that is generated by a class of countably generated modules with countably presented \emph{first} (instead of all) syzygies.

However, this weakened requirement still allows us to conclude the definability of the class $\mathcal{B}$.

Moreover, in the sense of \cite{Saroch&Stovicek}, the notion of hereditariness is weakened to the requirement that $\mathcal{A}={^{\perp_k}\mathcal{B}}$ and $\mathcal{A}^{\perp_k} = \mathcal{B}$ holds at least for $k = 1,2$, so that the left class is closed under kernels of epimorphisms iff the right one is closed under cokernels of monomorphisms; and the latter equivalence of definitions is due to the homological symmetry of cotorsion pairs.

This new notion suffices to prove the result, so that it seems a more 'natural' definition, in order to include also non-hereditary cotorsion pairs in the picture.
\smallbreak

With these premises, then, we advocate that the statement provided by Saroch and Stovicek is weaker only from the set-theoretical point of view.

However, from a meta-mathematical perspective, notice that this implies, at least in principle, a less rigid structure: in algebra there are several arguments related to some notion of finite 'dimension' which no longer hold in an infinite setting, as for instance Jacobson-Azumaya Theorem, duality isomorphisms to enrich braided categories of modules or the procedures based on the existence of upper bounds in directed colimits.

It is relevant to remark, in any case, that for many of these results there are partial generalizations highlighting the relevant algebraic 'underlying structure' in each proper set-theoretical setting.

With reference to the previous examples, recall an infinite version of Nakayama's Lemma, in its formulation regarding generating sets.

Moreover, worthy of notice is the valuable intuition (that becomes specifically formalized within the single categorical context of interest) that objects are dualizable whenever their 'size' is not greater than the 'additivity' of the category at stake (see 'Dualizable Object' in \emph{nLab}).

Finally, the last example leads us to the need to introduce the notion of $\lambda$-continuous directed limits. See the section \emph{$\lambda$-Continuous Limits} for an expanded motivation.

Hence, we could reasonably argue that the set-theoretical weakening eventually lets some inherent structure of cotorsion pairs emerge, thence what we mean by the expression 'more natural'.

Thus, we will follow the conventions and formulations introduced by the two authors of the aforementioned paper.

Finally, observe that, as achieved in \cite{Saroch&Stovicek}, our refined notion of hereditariness allows us to prove the countability type of the cotorsion pair $\mathcal{C}$ even requiring only the much weaker closure of the right class under well-ordered unions, i.e. well-ordered directed systems where all transition maps are not generic compatible morphisms, but actually compatible monomorphisms.
\smallbreak

The current version of the CTCM was proposed and proved by Saroch in \emph{'Approximations and Mittag-Leffler Conditions. The Tools'} (2018) in the following formulation.

We will see that, by restoring the original closure hypotheses under directed limits for the right class, we can even neglect the hereditary character of our cotorsion pair.
\bigskip

\textbf{Theorem.} \textsc{(CTCM, \cite{Saroch-tools})}
Let $\mathcal{C}=(\mathcal{A},\mathcal{B})$ be a cotorsion pair with $\lim \mathcal{B} = \mathcal{B}$. Then, $\mathcal{C}$ is of countable type, and $\mathcal{B}$ is definable.
\bigskip

Notwithstanding its inherent theoretical relevance, the CTCM also allows remarkable approaches to Enochs conjecture, as provided by Angeleri H\"ugel, Saroch and Trlifaj in \emph{'Approximations and Mittag-Leffler Conditions. The Applications.'} (2018), the twin paper to Saroch's work. Further developments appeared, for instance, in Saroch's Habilitation Thesis, \cite{Saroch-habilitation}.
\bigskip

The project is meant as an attempt to review the paper \cite{Saroch-tools} in an almost self-contained exposition grounding on the monograph \cite{Trlifaj&Goebel} by Trlifaj and G\"obel. It is then organized as follows.

The first section is devoted to \emph{Preliminaries}. Various formulations of the ML condition and stationarity are presented, together with their mutual relation.

Everything is explicitly stated and proved, except from \emph{Theorem} \cite{zimmermann},3.2 and 3.8, because it relies on heavy notation and is based on ideas that are distant from the spirit of the rest of the project.

Then, some set-theoretical prerequisites are introduced. In particular, we present the notion of filters, as well as a brief overview of the model theoretical approach, and we mention a 'revisited' approach to freeness conditions, following the work of Beke and Rosicky in \cite{Beke&Rosicky}.

Again, everything is made explicit, although the model theoretical digression on elementary cogenerators or Ziegler's spectrum are barely sketched, since each of them itself would require an independent review.
\smallskip

\smallskip

The second section is then devoted to developing the tools we will need to present a proof of the CTCM. In particular, we will survey various particular refinements of the notion of directed limits.

We will also introduce the Uniform Factorisation Condition, thus making the closure properties of stationary classes somewhat more explicit. Following the previous digression on some particularly well-behaved directed systems, we will generalize a criterion entailing the UFP to arbitrary cardinalities.

Thereafter, we will turn our attention to $C$-stationary classes of modules for $C$ a pure-injective, and leverage on the model theoretic formalism to provide for a 'canonical' class of pure-injective cogenerators, namely elementary cogenerators.

Finally, the last required device will consist of techniques to properly refine $\lambda$-dense systems of submodules of a given module, so as to form filtrations with good properties.

Our last endeavour will culminate in a 'Simplified Shelah's Compactness Principle', as proved by Saroch and Stovicek in \cite{Saroch&Stovicek}.
\bigskip

We will finally be ready to state and prove the current version of our main theorem, namely the CTMC, to which the third section will be dedicated. Finally, we will conclude by presenting an overview of the aforementioned applications of the CTMC to the Enochs conjecture, and by reporting an outlook on future research. In particular, we will mainly focus on $Add$ classes, from the viewpoint of the Enochs conjecture for small precovering classes, as well as from that of contramodules.

\newpage

\section{Preliminaries}
    
    In what follows, $Ring$ will always denote the category of rings with identity; $R$ will denote a general ring.
    
    We will adopt the following notation: \emph{f.} for finite, \emph{c.} for countable; \emph{g.} will stand for generated, while \emph{p.} for presented.
    Moreover, an abbreviation like $(<\kappa)-p.$ will mean 'less than $\kappa$ presented', and so on and so forth.
    
    Throughout the whole review, we will freely cite the results contained in the monograph \cite{Trlifaj&Goebel} by G\"obel and Trlifaj.

    \subsection{Purity, Definable classes and Coherent Functors}
    We state the following characterization of purity, as in \cite{Trlifaj&Goebel},2.19
    \begin{lemma}
    For a SES $\xymatrix{0 \ar[r] & A \ar[r] & B \ar[r] & B/A \ar[r] & 0}$ TFAE:
    \begin{enumerate}
        \item $A \subseteq_{*} B$.
        
        \item Any finite system of $R$-linear equations in $A$ that is solvable in $B$ must be solvable also in $A$.
        
        \item The given SES is the direct limit of a directed system of split SES's.
        
        \item The functor $(F\otimes \_)$ preserves the exactnes of the given SES for any $F$ (equiv. $F$ f.p.).
        
        \item The induced SES of character modules obtained by applying $[\_,\mathbb{Q}/\mathbb{Z}]_\mathbb{Z}$ is split-exact.
    \end{enumerate}
    \end{lemma}
    
    Let us now introduce a useful class of additive functors, as presented in the work of Crawley-Boevey \cite{coherent functors}.
    
    \begin{definition}
    An additive functor $F:Mod-R \rightarrow Mod-\mathbb{Z}$ is \textbf{coherent} provided that it commutes with arbitrary products and direct limits.
    \end{definition}
    
    \begin{definition}
    A class of modules is said to be \textbf{definable} provided that it is closed under arbitrary products, direct limits and pure submodules.
    
    Given any class of modules $\mathcal{G}$, we call \textbf{definable closure} $def(\mathcal{G})$ of $\mathcal{G}$ the smallest of its definable sup-classes.
    \end{definition}
    
    One can prove that definable classes are characterized as vanishing loci of sets of coherent functors. However, this is out of the scope of the present review. Thus, we will just provide a direct proof of the following useful, although weaker, property.
    
    \begin{lemma}
    The vanishing locus of a set of coherent functors is closed under definable closure and pure quotients.
    \end{lemma}

    \begin{proof}
    Notice that coherent functors preserve arbitrary products and direct limits, so we only need to check that the vanishing locus of a coherent functor is further closed under pure submodules. The latter property follows only by means of additivity and commutativity with $\lim$.

    Indeed, consider any pure-SES $E$ with middle-term in the vanishing locus of some coherent functor $F$. By the characterization $3.$ of purity, we can regard it as the colimit of a directed system of split-SES $(E_i)_I$. $F$ is additive, so it preserves the splitting property. Finally, it commutes with direct limits (that are exact in module categories), so that $F(E)=F(\lim_I E_i) = \lim_I F(E_i)$ is still a SES (and furthermore a pure-SES). Its middle term vanishes by assumption, and so the outer terms do by the exactness.
    \end{proof}
    
    \subsection{Mittag-Leffler condition}
    
    In this section we introduce the ubiquitous definitions of Mittag-Leffler condition and stationarity, and we will subsequently consider the properties of the objects satisfying these statements from various perspectives.
    
    \begin{definition}
    Given an inverse system $\mathcal{H}=(H_\alpha , h_{\alpha,\beta}: H_\beta \rightarrow H_\alpha | \alpha < \beta \in I)$ indexed by some set $I$, we say that it satisfies the \textbf{Mittag-Leffler condition} (in the following abbreviated as 'ML') provided that it holds:
    \vspace{-.4cm}
    
    \begin{equation*}
    \forall \alpha \in I \quad \exists \beta \in I \quad s.t. \quad \forall (\beta <) \gamma \in I, \quad Im(h_{\alpha,\beta}) = Im(h_{\alpha, \gamma})
    \end{equation*}
    
    \begin{itemize}
        \item A direct system $\mathcal{M}=(F_\alpha, f_{\beta,\alpha}:F_\alpha \rightarrow F_\beta | \alpha < \beta \in I)$ is said to be a \textbf{Mittag-Leffler directed system} if $[\_,B]_R$ induces a ML inverse system for each $R$-module $B$.
        
        \item A module $M$ is said to be a \textbf{Mittag-Leffler module} provided that each one of its directed presentations consisting of f.p.\ modules is ML.
    \end{itemize}
    
    \end{definition}
    
    \emph{Remark.}
    The ML condition is said to be \textbf{strict} whenever it further holds for the distinguished $\beta$'s that $Im(h_{\alpha, \beta}) = Im(h_\alpha)$, i.e. the sequence of images stabilizes at the limit morphism.

    \begin{definition}
    A directed system $\mathcal{M}=(F_\alpha , f_{\beta, \alpha} | \alpha < \beta \in I)$ indexed by some set $I$, is \textbf{$B$-stationary} provided that $[\_,B]_R$ induces a ML inverse system $\mathcal{H}^{B}$.
    
    We can then naturally extend the definition to modules (considering directed presentations consisting of f.p.\ modules, to save the analogy with ML modules).
    \end{definition}
    
    \emph{Remark.}
    In other words, the notion of $B$-stationarity specifies the module $B$ for which the ML condition is induced. If our objects induce strict ML properties, then they are said to be strict $B$-stationary.
    
    Moreover, notice that we can define $\mathcal{B}$-stationarity for any class $\mathcal{B} \subseteq Mod-R$, by requiring the previous property for each $B\in \mathcal{B}$.
    \bigskip
    
    We can also consider a relative version of the ML condition. For reference, see \cite{H-subgroups} by Angeleri H\"ugel and Herbera.
    
    In the absolute case, a module $M$ is ML iff for any choice of a directed presentation of f.p.\ modules and for any module $B$, the induced inverse system is $B$-stationary; similarly, $M$ is $B$-stationary iff it admits some presentation consisting of f.p.\ modules from which $[\_,B]_R$ induces a ML inverse system.
    
    However, in the relative case, the notion of stationarity turns out to present a remarkable formal analogy to the relative ML case. However, we postpone the discussion to the next section, so as to be able to employ the language of H-subgroups to shorten the notation.
    
    \subsubsection{H-subgroups}
    
    We will now present an equivalent approach in terms of H-subgroups, as adopted by Angeleri H\"ugel and Herbera in \cite{H-subgroups}. 
    
    \begin{definition}
    Given the modules $M,N,B$ and the $R$-morphism $v:M\rightarrow N$, we call \textbf{H-subgroup} the $End_R(B)$-submodule of $[M,B]_R$ given by $H_v (B) := Im([v,B]_R)$.
    \end{definition}
    
    With reference to the previous notation, consider any directed system of modules $\mathcal{M}$.
    
    Remark that, by the compatibility of the morphisms in $\mathcal{M}$, for each $B\in \mathcal{B}$ and $l\in I$, $Im(h^{B}_{l}) \subseteq Im(h^{B}_{jl}) \quad \forall l\leq j \in I$, so that it always holds $Im(h^{B}_{l})\subseteq \bigcap_{j\leq l} Im(h^{B}_{jl})$.
    
    In terms of H-subgroups, this amounts to $H_{f_{ji}}(B) \subseteq H_{f_{li}}(B)$ for each $i<j\leq l$, and $H_{f_{i}}(B) \subseteq H_{f_{li}}(B)$ for each $i\leq l$.
    
    Hence, the $B$-stationarity of the system is equivalent to the reverse inclusion $H_{f_{ji}}(B) \supseteq H_{f_{li}}(B)$ $(\exists l\geq j)$, and the strictness is achieved whenever it further holds $H_{f_{i}}(B) \supseteq H_{f_{li}}(B)$ $(\exists l\geq i)$.
    \bigskip
    
    We are now ready to remark the already mentioned formal analogy in the relative case. To introduce such a subtlety, however, we need the preliminary notion of domination.
    
    \begin{definition}
    Let $Q \in R-Mod$, $B\in Mod-R$. Given a co-angle $u:M\rightarrow N$, $v:M\rightarrow M'$ in $Mod-R$, we say that \textbf{$u$ $B$-dominates $v$ w.r.t. $Q$} provided that $Ker(u\otimes Q) \subseteq \bigcap_{h\in H_v(B)} Ker(h\otimes Q)$.
    
    We can then extend in the obvious way to classes $\mathcal{Q},\mathcal{B}$.
    \end{definition}
    
    In other words, we are trying to 'approximate' the hypothesis of the Factor Theorem to factorize $v\otimes Q$ by $u\otimes Q$ over $\mathcal{Q}$, namely $ker(u\otimes Q) \subseteq ker(v\otimes Q)$, by requiring the sought inclusion over all the right translations of $v$.
    
    We will now state the aforementioned analogy occurring in the relative setting via the following definitions-lemmas. The proofs, however, are omitted, because they are long and distant from the spirit of the rest of the dissertation.
    
    \begin{lemma}{(Relative stationary, \cite{H-subgroups},4.8)}
    In the usual notation, consider a class $\mathcal{S}$ consisting of f.p.\ modules, and take a module $M \in \lim \mathcal{S}$. TFAE:
    \begin{enumerate}
        \item $M$ is $B$-stationary.
        
        \item $\exists (\iff \forall )$ a directed presentation $\mathcal{M}$ consisting of f.p.\ modules such that, for each index $i\in I$, there is some $j\geq i$ for which $f_{ji}$ \textbf{$B$-dominates} $f_i$.
        
        \item For each f.p.\ $F$(w.l.o.g.\ $F\in \mathcal{S})$, and for each map $u:F\rightarrow M$, there is a factorization through $\mathcal{S}$ that is $B$-dominating. In other words, $\exists v:F\rightarrow S\in \mathcal{S}$ s.t. $u \in H_v(M)$ and $v$ $B$-dominates $u$.
    \end{enumerate}
    \end{lemma}
    
    \begin{lemma}{(Relative ML, by Facchini and Azumaya, \cite{H-subgroups},5.1; or adapt \cite{Trlifaj&Goebel},3.14)}
    Let $\mathcal{Q} \subseteq R-Mod$, $\mathcal{S}\subseteq Mod-R$ consisting of f.p.\, and take any $M\in \lim \mathcal{S}$. TFAE:
    
\begin{enumerate}
    \item $M$ is $\mathcal{Q}$-ML.
    
    \item $\forall$ directed presentation $\mathcal{M}$ consisting of f.p.\ modules, for any chosen index $i\in I$ there is some $j\geq i$ s.t. $f_{ji}$ dominates $f_i$ \textbf{w.r.t. $\mathcal{Q}$}.
    
    \item For each f.p.\ $F$(w.l.o.g.\ $F\in \mathcal{S})$, and for each map $u:F\rightarrow M$, there is a mutual factorization through $\mathcal{S}$. In other words, $\exists v:F\rightarrow S\in \mathcal{S}$ s.t. $u \in H_v(M)$ and $ker(u\otimes Q) = ker(v\otimes Q)$ $\forall Q\in \mathcal{Q}$ (i.e. the assumptions of the Factor Theorem are satisfied and one also has the converse factorization $v\otimes Q$ through $u\otimes Q$).
    
    \item For each $X\subseteq M$ with $\# X \leq \aleph_0$, there is some c.p.\ $\mathcal{Q}$-ML $N\in \lim \mathcal{S}$ together with some $v:N\rightarrow M$ s.t. $x \subseteq v(N)$ and $v\otimes Q$ mono for each $Q\in \mathcal{Q}$.
\end{enumerate}
    \end{lemma}
    
    Remark that condition $3.$ of the previous two statements provides a factorization as that entailed by Lenzing Lemma (see \cite{Trlifaj&Goebel},2.13) for limits of countably presented modules, that still satisfies some domination or factorisation properties after being tensored by the modules at stake.
    
    \begin{corollary}{(Local construction of $\mathcal{Q}$-ML, \cite{H-subgroups},5.2)}
    
    Let $\mathcal{Q}\subseteq R-Mod$, and consider a class $\mathcal{S} \subseteq Mod-R$ consisting of f.p.\ modules. Fix some $M \in \lim \mathcal{S}$, and denote by $\mathcal{C}$ the class of all the c.g.\ $N\leq M$ s.t. $N$ is $\mathcal{Q}$-ML and $\epsilon:N \xhookrightarrow{} M$ is preserved by $(\_ \otimes Q)$ for each $Q\in \mathcal{Q}$ (i.e. $\epsilon$ is $\mathcal{Q}$-pure).
    
    Then, $M$ is $\mathcal{Q}$-ML iff it is the directed union of the elements of $\mathcal{C}$.
    
    Moreover, provided that $R\in \mathcal{Q}$, $\mathcal{C}$ can be assumed to consist of c.p.\ modules from $\lim \mathcal{S}$.
    \end{corollary}
    
    \begin{corollary}{(c.g.\ iff c.p.\, \cite{H-subgroups},5.3)}\label{c.g. iff c.p.}
    Let $\mathcal{Q}\subseteq R-Mod$ containing $_{R}R$. Then, for each $\mathcal{Q}$-ML module $M$, the notions of c.g.\ and c.p.\ coincide.
    \end{corollary}
    
    \begin{proof}
    By the previous \emph{Corollary}, $M$ is the directed union of some family from $\mathcal{C}$, that consists of c.p.\ modules. Now, $M$ is c.g.\ by assumption, so we can extract a countable subchain whose directed union is $M$, and we are done.
    \end{proof}
    
\subsubsection{n-pointed modules}
    
    Finally, we will introduce ML conditions via $n$-pointed modules, as done by Zimmermann in \cite{zimmermann}.
    
    \begin{definition}
    We call $n$-pointed module any pair $(A,a)$ with $A\in Mod-R$ and $a\in A^n$. We make them into the category of $n$-modules $Mod-R_n$ by declaring as morphisms those arrows $h:(A,a)\rightarrow (M,m)$ whose action as elements of $[A,M]_R$ respects the pointed vectors via $h:a\mapsto m$.
    
    Given a directed system $\mathcal{M}:= (M_\alpha, u_{\beta,\alpha}|\alpha < \beta \in I)$ in $Mod-R$ with limit $(M,u_\alpha | \alpha \in I)$, it induces a directed system of $n$-pointed modules as follows: for each $\alpha \in I$ choose $x_\alpha \in M_\alpha ^{n}$ satisfying $u_{\beta,\alpha}:x_{\alpha} \mapsto x_{\beta}$, and pick up $m:=u_{\alpha} (x_{\alpha})$ $(\forall \alpha \in I)$. Then, $\big( (M_\alpha,x_\alpha),u_{\beta,\alpha} | \alpha <\beta \in I \big)$ is a direct system in $Mod-R_n$ with limit $(M,m)$.
    \end{definition}
    
    \emph{Remark.} The directed sequence of pointed vectors $\{x_\alpha\}_\alpha$ actually exists up to a choice of a well-ordered directed subsystem; such a choice can be performed as shown in the section on \emph{Direct Limits}.
    
    \bigskip
    
    We will state, but not prove, an auxiliary characterization of strict $B$-stationarity. (see \emph{Th.}\cite{H-subgroups},8.10 by Angeleri H\"ugel and Herbera, and \emph{Lemma}\cite{zimmermann},3.2 by Zimmermann.)
    
    Let's first introduce some notation: given $A,B \in Mod-R$ and $a\in A^n$, we define $\epsilon_a \in \big[ [A,B]_R,B^n \big]_{End(B)}$ by the rule $\epsilon_a: f \mapsto f(a)=(fa_i)_{i=1}^{n}$. Moreover, set $H_{A,a}(B):=Im(\epsilon_a)= \{ f(a) | f\in [A,B]_R \}$.
    
    \emph{Remark.} The following straightforward properties hold:
    \begin{itemize}
        \item If $a = (a_i)_{i=1}^{n}$ generates $A$, then $\epsilon_a$ is mono.
        
        \item For $v \in [A,M]_R$ s.t. $v:a\mapsto m$, $\epsilon_a (H_v(B))=H_{M,m}(B)$.
        
        \item For each $m\in M^n$, $H_{M,m}(B)=\epsilon_e (H_u(B))$, where $e$ is the canonical basis of $R^n$ and $u:R^n \rightarrow M$ acts as the $n$-scalar product on the $n$-vector $m$.
    \end{itemize}
    
    \begin{definition}
    For $B\in R-Mod-S$, $V\in S-Mod$, $M\in Mod-R$, consider the following morphism of abelian groups (see \emph{Lemma} \cite{Trlifaj&Goebel},2.16 about the properties of $Ext$ and $Tor$):
    \vspace{-.3cm}
    
    \begin{equation*}
        \begin{split}
            \nu \equiv \nu(M,B,V): M \otimes_R \big( _{S}[B,V]\big) & \longrightarrow {_{S}\big[ [M,B]_R, V  \big]} \\
            m\otimes \phi \quad & \longmapsto [\nu(m\otimes \phi):f\mapsto \phi(f(m))]
        \end{split}
    \end{equation*}
    \end{definition}
    
    \begin{lemma}
    Given the $R$-modules $M,B$, TFAE:
    \begin{enumerate}
        \item $M$ is $B$-stationary.
    
        \item For each $n<\omega$, $m \in M^n$, and for each direct system of $n$-pointed modules $\big( (M_\alpha,x_\alpha),  u_{\beta,\alpha} | \alpha \leq \beta \in I \big)$ with f.p.\ $M_\alpha$'s and colimit $(M,m)$, there is some $\beta \in I$ s.t. $H_{M,m}(B)=H_{M_\beta,x_\beta}(B)$.
    
        \item For each $n<\omega$, $m \in M^n$, there is some f.p.\ $R$-module $A$ and some $n$-tuple $a\in A^n$ s.t. $\exists h\in [(A,a),(M,m)]_{R,n}$ morphism of $n$-pointed modules s.t. $H_{M,m}(B)=H_{A,a}(B)$.
    \end{enumerate}

    Moreover, if $B\in S-Mod-R$, then also the following are equivalent:
    \begin{itemize}
        \item For each injective $_{S}V$, the canonical morphism $\nu(M,B,V)$ is mono.
        
        \item For each injective cogenerator $_{S}U$, the canonical morphism $\nu(M,B,U^{I})$ is mono for each set $I$.
    \end{itemize}
    \end{lemma}
    
    \emph{Remark.} We could weaken the extra condition on the injective module because any injective is a direct summand in some product of copies of an injective cogenerator, so that it suffices to consider character modules (i.e. $S$-modules of the form ${_{S}[M_R,{_S}U_R]}$ for ${_S}U$ injective cogenerator).
    
    \begin{corollary}
    Let $B\in Mod-R$. Then the following statements hold true:
    \begin{enumerate}
        \item The class of strict $B$-stationary modules is closed under pure submodules and pure extensions.
        
        \item A direct sum of modules is strict $B$-stationary iff each of these direct summands is such.
        
        \item $B$ is $\Sigma$-pure-injective iff every module in $Mod-R$ is strict $B$-stationary.
    \end{enumerate}
    \end{corollary}
    
    \begin{proof}
    $1).$ Given any pure SES $0 \rightarrow M \rightarrow_{\ast} N \rightarrow N/M \rightarrow 0$, $\nu$ induces a diagram of SES's (the exactness of the rows is given by purity):
    \begin{center}
        $\xymatrix{
        0 \ar[r] & {M\otimes (_{S}[B,V])} \ar[r] \ar[d]^{\nu_M} & {N\otimes (_{S}[B,V])} \ar[r] \ar@{^{(}->}[d]^{\nu_N} & {N/M \otimes (_{S}[B,V])} \ar[r] \ar[d]^{\nu_{N/M}} & 0 \\
        0 \ar[r] & {_{S}\big[ [M,B]_R, V  \big]} \ar[r] & {_{S}\big[ [N,B]_R, V  \big]} \ar[r] & {_{S}\big[ [N/M,B]_R, V  \big]} \ar[r] & 0
        }$
    \end{center}
    
    and we can conclude by the Snake Lemma.
    
    $2).$ Notice that for $M=\oplus_I M_i$ the canonical inclusion $\epsilon_i$ induces an inclusion $\epsilon^{\ast \ast}_i:{_{S}\big[ [M_i,B]_R, V  \big]} \rightarrow {_{S}\big[ [M,B]_R, V  \big]}$ by the left-exactness of the Hom-functor.
    
    The family of inclusions induces then an injective map by means of universality, namely
    \begin{equation*}
        \begin{split}
            E: \oplus_I \big(_{S}\big[ [M_i,B]_R, V  \big]\big) & \longrightarrow {_{S}\big[ [M,B]_R, V  \big]} \\
            (g_i)_I \quad & \longmapsto E(g_i)_I : f \mapsto \sum_I g_i \circ (f|M_i)
        \end{split}
    \end{equation*}
    
    Observe, then, that $\nu(\oplus_I M_i,B,V)$ is mono iff $\nu(M_i,B,V)$ is mono for each $i$.
    
    Indeed, the equivalence is obtained from the following commutative diagram:
    \begin{center}
        $\xymatrix{
    {M\otimes (_{S}[B,V])} \ar[rr]^{\nu} && {_{S}\big[ [M,B]_R, V  \big]} \\
    {\oplus_I M_i\otimes (_{S}[B,V])} \ar[rr]^{\oplus_I \nu_i} \ar@{=}[u] && {\oplus_I \big(_{S}\big[ [M_i,B]_R, V  \big]\big)} \ar@{^{(}->}[u]
    }$
    \end{center}
    
    $3).$ Is deduced by several results in \cite{zimmermann}, see \emph{Theorem} 3.8.
    \end{proof}
    
    We are finally ready to report a useful strict-stationarity criterion for SES, namely \emph{Lemma} \cite{Saroch-tools},4.4.
    
    \begin{lemma}\label{stationarity along ses}
    Let $B$ be an $R$-module and $0 \rightarrow N \rightarrow A \rightarrow M \rightarrow 0$ SES s.t. $M \in {^{\perp}(B^{(I)})^{cc}}$ for each set $I$.

    Then, $A$ strict $B$-stationary $\implies$ $N$ strict $B$-stationary.
    \end{lemma}

    \begin{proof}
    We have observed that $N$ is strict $B$-stationary iff $\nu(N,B,{\mathbb{Q}/\mathbb{Z}}^I)$ is mono for each set $I$.

    Since $[B,{\mathbb{Q}/\mathbb{Z}}^I]_{\mathbb{Z}} \simeq (B^{(I)})^{c}$, by \emph{Lemma} \cite{Trlifaj&Goebel},2.16 and the assumptions on $M$ on has that $Tor_1^R \big( M, [B,{\mathbb{Q}/\mathbb{Z}}^I]_{\mathbb{Z}} \big) \simeq 0$. Thus, we obtain the following commutative diagram of exact SES's (the injectivity of the first row corresponds to the condition on $Tor_1$):

    \begin{center}
        $\xymatrix{
        0 \ar[r] & {N \otimes_R \big[B,{\mathbb{Q}/\mathbb{Z}}^{I}\big]_\mathbb{Z}} \ar[d]^{\nu_N} \ar@{^{(}->}[rr] && {A \otimes_R \big[B,{\mathbb{Q}/\mathbb{Z}}^{I}\big]_\mathbb{Z}} \ar[d]^{\nu_A} \\
        & {\big[ [N,B]_R, {\mathbb{Q}/\mathbb{Z}}^{I} \big]} \ar[rr] && {\big[ [N,B]_R, {\mathbb{Q}/\mathbb{Z}}^{I} \big]}
        }$
    \end{center}

    The injectivity of $\nu_A$ provides that of $\nu_N$, as required.
    \end{proof}
    
    \subsubsection{Stationarity and Pure Injectivity}
    
    The following result relating stationarity and pure-injectivity was proved by Herbera as \emph{Proposition} \cite{Herbera},1.7.
    
    \begin{theorem}
    Let $B\in Mod-R$ be (locally) \texttt{pure-injective}. Then, for each directed system (of f.g.) $\mathcal{M}$, the inverse system induced by $[\_,B]_R$ satisfies $H_{f_i}(B) = \bigcap_{i\leq j} H_{f_{ji}}(B)$ for each $i\in I$.
    
    In particular, $\mathcal{M}$ is $B$-stationary iff $\mathcal{M}$ is strict $B$-stationary.
    \end{theorem}
    
    \begin{proof}
    \texttt{Assume} that our pure-injective module $B \in Mod-R$ is a dual module of the form $Q^c \equiv [Q,\mathbb{Q}/\mathbb{Z}]_\mathbb{Z}$ for some $Q\in R-Mod$.
    
    Then, we have the following diagram of canonical morphisms:
    
    \begin{center}
        $\xymatrix{
        [\lim M_i, Q^c]_R \ar[d]_{f_j^{\ast}} \ar[r]^{f^{\ast}_i} & [M_i, Q^c]_R \ar@{=}[d] \ar[r] & Coker(f_i^{\ast}) \ar@{.>}[d]_{\exists !} \ar@/^3pc/@{.>}[dd] \ar[r] & 0\\
        [M_j,Q^c]_R \ar[r]^{f_{ji}^{\ast}} & [M_i,Q^c]_R \ar[r] & Coker(f_{ji}^{\ast}) \ar[r] & 0 \\
        & & {\varprojlim_{(i\leq ) j} Coker(f_{ji}^{\ast})} \ar[u]^{\nu_j}
        }$
    \end{center}
    
    and, by universality, we obtain a morphism $\mu_i: Coker[f_i,Q^c]_R \rightarrow \varprojlim_{j} Coker(f_{ji}^{\ast})$.
    
    Observe that the thesis for $Q^c$ is equivalent to $\mu_i$ being an isomorphism. Indeed, fixed $i \in I$, by the compatibility of the family $\{ f_{ji}^{\ast} \}_j$ we would obtain the required isomorphism:
    
    $$\frac{[M_i,Q^c]}{H_{f_{i}}(Q^c)} = Coker(f_i^{\ast}) \simeq \varprojlim_j Coker(f_{ji}^{\ast}) = \frac{[M_i,Q^c]}{\cap_j H_{f_{ji}}(Q^c)}$$
    \smallskip
    
    \emph{CLAIM.} $\mu_i$ iso for each given $i\in I$.
    \smallskip
    
    \emph{PROOF.} By \emph{Lemma} \cite{Trlifaj&Goebel},2.2, there is an exact sequence $(\triangle)$:
    
    \begin{center}
        $\xymatrix{
        0 \ar[r] & {\cup_{j \geq i} Ker(f^{\ast}_{ji} \otimes Q)} \ar[r] & {M_i \otimes Q} \ar[rr]^{f^{\ast}_i \otimes Q} && {\lim_l (M_l \otimes Q)}
        }$
    \end{center}
    
    Dualization via $[\_,\mathbb{Q}/\mathbb{Z}]_\mathbb{Z}$ and the $(\otimes,Hom)$ adjunction (remark that directed $\lim$ of f.p.\ can be carried out of the first component of $Hom$ functors), then, yield the exact sequence:
    
    \begin{center}
        $\xymatrix{
        {\varprojlim [M_l,Q^{c}]_R} \ar[r] & [M_i,Q^c]_R \ar[r] & {\varprojlim_{j\geq i} [Ker(f^{\ast}_{ji} \otimes Q)]^{c}} \ar[r] & 0
        }$
    \end{center}
    
    Now, remark the following isomorphisms, again induced by $(\otimes,Hom)$ adjunction:
    
    $$(1) \qquad Coker[f^{\ast}_{ji},Q^c] \simeq [Ker(f^{\ast}_{ji}\otimes Q)]^c$$
    
    \vspace{-.4cm}
    
    $$(2) \qquad Coker[f^{\ast}_{i},Q^c] \simeq [Ker(f^{\ast}_{i}\otimes Q)]^c$$
    
    Thus, we obtain the commutative diagram:
    
    \begin{center}
        $\xymatrix{
        Coker[f^{\ast}_{i},Q^c] \ar[d]_{(2)}^{\simeq} \ar[rr]^{\mu_i} && {\varprojlim_{j\geq i} Coker[f^{\ast}_{ji},Q^c]_R} \ar[d]^{(1)+(\triangle)}_{\simeq} \\
        [Ker(f^{\ast}_{i} \otimes Q)]^c \ar[rr]_{(\triangle)}^{\simeq} && [\cup_{j\geq i} Ker(f^{\ast}_{ji} \otimes Q)]^c
        }$
    \end{center}
    
    that forces $\mu_i$ to be an isomorphism.
    \bigskip
    
    Now, let's consider the general case, i.e. take any pure-injective $B$, and set $Q:=B^c$.
    
    From the very definition, $ev:B\xhookrightarrow{}_{\ast} B^{cc}$ splits; call $h$ a retraction of $ev$.
    
    The first part of the proof provides for the inclusion
    \vspace{-.3cm}
    
    $$\cap_{j\geq i} H_{f_{ji}}(B) \overset{ev}{\xhookrightarrow{}} \cap_{j\geq i} H_{f_{ji}}(B^{cc}) = H_{f_{i}}(B^{cc})$$
    
    So that for each $u \in \cap_{j\geq i} H_{f_{ji}}(B)$ in the first member, there is some $g\in \varprojlim [M_l,B^{cc}]_R$ s.t. $ev \circ u = g\circ f^{\ast}_i$, that is $u = h\circ g \circ f^{\ast}_i \in H_{f_{i}}(B)$, as required.
    \end{proof}
    
    \emph{Remark.} Observe that the statement of the previous Theorem does not imply that the sequence $\big(H_{f_{ji}}(B)|i\leq j \in I\big)$ does ever stabilize; however, if it does before reaching the limit, then the image of the stabilizing morphism equals that of the corresponding canonical limit map.
    \bigskip
    
    The following \emph{Lemma} assures some useful closure properties of stationarity classes (slight modification of \emph{Cor.} \cite{H-subgroups},3.9).
    
    \begin{lemma}\label{closure properties}
    Given a class $\mathcal{G} \subseteq Mod-R$ and a $\mathcal{G}$-stationary $R$-module $M$, the following statements hold true:
    \begin{enumerate}
        \item w.l.o.g.\ $\mathcal{G}$ is closed under copies of pure submodules or pure quotients.
    
        \item TFAE:
            \begin{itemize}
                \item $M$ is $Add\mathcal{G}$-stationary.
        
                \item $M$ is $Prod \mathcal{G}$-stationary.
                
                \item $M$ is $def\mathcal{G}$-stationary.
        
                \item $M$ admits a direct presentation $\mathcal{M}$ consisting of f.p.\ modules s.t. we can define a strictly increasing function $s:I\rightarrow I$ witnessing the uniform factorization property over $\mathcal{G}$. In other words, $s$ satisfies the following property independently of $B\in \mathcal{G}$: each $c\in [M_i,B]_R$ that belongs to $H_{f_{s(i),i}}(B)$ is also in $H_{f_{j,i}}(B)$ for each $s(i) \leq j \in I$
            \end{itemize}
    
        \item $M$ is $Add(B),Prod(B),def(B)$-stationary for each $B\in \mathcal{G}$.
    \end{enumerate}
    \end{lemma}

    \begin{proof} $1)+3).$ $\mathcal{G}$-stationarity implies that there is a presentation of $M$ indexed by $I$ s.t. for each $B \in \mathcal{G}$ and for each index $i\in I$, there is some $j=j(B,i)$ s.t $H_{f_{ji}}(B) = H_{f_{li}}(B)$ for every $(j<)l \in I$.

    So, $B$ belongs to the vanishing locus of the coherent functors $(F_{l}^{B,i}(\_):=H_{f_{j(B,i),i}}(\_) / H_{f_{l,i}}(\_) | i\leq j \leq l < I)$ for each $i\in I$, so that also $def(B)$ lies there, and we obtain the $def(B)$-stationarity of $M$, that implies both $1)$ and $3)$.
    \smallskip

    As for $2)$, notice that, since $\cap$ and $\prod$ commute, $\prod B_j$-stationarity for each subset $(B_j)_J$ amounts to the uniformity of the factorisation, and hence to the $def(\mathcal{G})$-stationarity.

    Now, observe that also $\cap, \oplus$ commute, and that, in our setting, we can bring $\oplus$ out of the covariant Hom-functor. Indeed, for $(M_i)_I$ a presentation of $M$ consisting of f.p, $[M_i,\oplus_J B_j]_R \simeq \oplus_J [M_i, B_j]_R$ by \emph{Lemma} \cite{Trlifaj&Goebel},2.7. So, we can similarly replace $\prod$ by $\oplus$ in the previous equivalence.

    Finally, we can include also the directs summands of arbitrary $\prod,\oplus$ because they are in particular pure submodules: $N \leq_\oplus M$ is the direct limit of the countable directed tower induced by the splitting section after the retraction $M \rightarrow N \rightarrow M$.
    \end{proof}
    
    Let's now present a result by Angeleri H\"ugel and Herbera (see \emph{Proposition} \cite{H-subgroups},8.14) clarifying the relationship between stationarity and relative ML conditions.
    
    Recall that, given $M\in S-Mod$, $V$ an injective cogenerator of $S-Mod$, we denote the 'generalized' dual module of $M$ by $M^c := [M,V]_S$.
    \smallskip
    
    \begin{lemma}\label{ML vs stationarity}
    A Module $M$ is $C$-ML iff $M$ is strict $C^{c}$-stationary.
    
    Moreover, any strict $C$-stationary module is also $C^c$-ML. The converse holds whenever $C$ is pure-injective.
    \end{lemma}
    
    \begin{proof} The proof of the first equivalence differs from that provided by the authors, and uses Facchini and Azumaya's characterization of relative ML conditions, as provided in the section \emph{Preliminaries}.
    
    $M$ is $C$-ML iff given any $\mathcal{M}$ presentation of $M$ consisting of f.p.\ modules, for each index $i\in I$, there is some index $j\geq i$ s.t. $f_{ji}$ dominates $f_i$ w.r.t. $C$.
    
    In other words, it holds $Ker(f_{ji}\otimes C) \subseteq \bigcap_{dom(h') = M} Ker(h'\circ f_i \otimes C)$.
    
    Dualizing by some injective cogenerator $V$, and by exploiting the adjunction $(\_ \otimes C, Hom(\_,C))$ one obtains that our statement is equivalent to the following epi (for each arbitrary choice of $i$):
    
    $$\frac{[M_i,C^c]}{H_{f_{ji}}(C^c)}=Coker[f_{ji},C^c] \twoheadrightarrow \frac{M_i,C^c}{\sum H_{(h\circ f_{i})}(C^c)} = \frac{[M_i,C^c]}{H_{f_{i}}(C^c)}$$
    
    and the latter holds iff $H_{f_{ji}}(C^c) \subseteq H_{f_{i}}(C^c)$, that is the strict $C^c$-stationarity of $M$.
    
    Indeed, the epi is given by the following reasoning:
    
    Endowe $\mathcal{H}:=(h\in [M_i,N]_R | N\in Mod-R)$ with the partial order $h\leq h' \iff Ker(h\circ f_i \otimes C) \subseteq Ker(h'\circ f_i \otimes C)$.
    
    Notice that, for each choice of $h,h'$, they always have $id$ as an upper bound.
    
    Hence, $\mathcal{H}$ is a directed poset, and we can consider the inverse system
    \vspace{-.2cm}
    
    $$\mathcal{K}:= \{ K_h:=Ker(h\circ f_i \otimes C), e_{h'h}:K_{h'}\xhookrightarrow K_h|h\leq h' \in \mathcal{H} \}$$
    
    More explicitly, one may construct the following inverse system of exact sequences:
    
    \begin{center}
        $\xymatrix{
        0 \ar[r] & {Ker(h \circ f_i \otimes C)} \ar[r] & {M_i \otimes C} \ar[rr]^{h\circ f_i \otimes C} && {N \otimes C} \\
        0 \ar[r] & {Ker(h'\circ f_i \otimes C)} \ar[r] \ar@{^{(}->}[u]^{e_{h'h}} & {M_i \otimes C} \ar@{=}[u] \ar[rr]^{h'\circ f_i \otimes C} && {N' \otimes C} \ar@{.>}[u]_{"\mu_{h'h}"}
        }$
    \end{center}
    
    where the vertical dotted arrow $\mu_{h'h}$ is an essential (over $N'\otimes C$) extension of the canonical induced morphism
    
    $$Im(h'\circ f_i \otimes C) = \frac{M_i \otimes C}{K_{h'}} \overset{\mu_{h'h}}{\twoheadrightarrow} \frac{M_i \otimes C}{K_h} = Im(h \circ f_i \otimes C)$$
    
    Notice that $\varprojlim \mathcal{K} = \big( \bigcap_{\mathcal{H}} K_h,\iota_h | h\in \mathcal{H} \big)$ is the intersections of the kernels.
    
    After dualization, we obtain the directed system 
    $$\mathcal{K}^{c} = \{ C_h:= Coker[h\circ f_i,C^c], e_{h'h}^{\ast}:C_h \twoheadrightarrow C_{h'}|h\leq h' \in \mathcal{H} \}$$
    
    More explicitly, we have the directed system of sequences
    
    \begin{center}
        $\xymatrix{
        [N,C^c] \ar[d]^{{"\mu_{h'h}"}^{\ast}} \ar[rr]^{[h\circ f_i,C^c]} && [M_i,C^c] \ar@{=}[d] \ar[r] & Coker[h\circ f_i,C^c] \ar@{->>}[d]^{e_{h'h}^{\ast}} \ar[r] & 0 \\
        [N',C^c] \ar[rr]^{[h'\circ f_i,C^c]} && [M_i,C^c] \ar[r] & Coker[h'\circ f_i,C^c] \ar[r] & 0 \\
        }$
    \end{center}
    
    And the arrow $e_{h'h}^{\ast}$ acts, up to iso, as a quotient projection.
    Indeed, by the right exactness of the tensor product, our arrow is induced by the projection $\mu_{h'h}$:
    
    \begin{center}
        $\xymatrix{
        {[Ker(h\circ f_i \otimes C),V]} \ar[r]^{\simeq} \ar@{->>}[d]^{e_{h'h}^{\ast}} & {\frac{[M_i,C^c]}{Im[h\circ f_i,C^c]}} \ar@{=}[r] & {\frac{[M_i,C^c]}{\big(Im(h\circ f_i)\otimes C\big)^c}} \ar@{->>}[d]^{"\mu_{h'h}"} \\
        {[Ker(h'\circ f_i \otimes C),V]} \ar[r]^{\simeq} & {\frac{[M_i,C^c]}{Im[h'\circ f_i,C^c]}} \ar@{=}[r] & {\frac{[M_i,C^c]}{\big(Im(h'\circ f_i)\otimes C\big)^c}}
        }$
    \end{center}
    
    Therefore, one has the sought
    \vspace{-.3cm}
    
    $$\lim \mathcal{K}^{\ast} = \bigg( C:=  \frac{[M_i,C^c]}{\sum H_{(h\circ f_i)}(C^c)}, \hspace{.2cm} \nu_h \hspace{.2cm} | \hspace{.2cm} h\in \mathcal{H} \bigg)$$
    
    with canonical maps $\nu_h:C_h \rightarrow C$ acting as the quotient projections $f+H_{(h\circ f_i)}(C^c) \xmapsto{} f+H_{f_i}(C^c)$.
    \bigskip
    
    As for the second part, we will report the proof only of the needed (converse) implication.
    
    Let $C$ be a pure-injective module, and assume to be given a $C^c$-ML module $M$. Let's show that $M$ is also (strict) $C$-stationary.
    
    Let $_{S}V$ be an injective cogenerator of $S-Mod$, and write $B^c := [B,V]_S$ for the generalized dual of $B$.
    
    Notice that $_{\mathbb{Z}}V$ is an injective cogenerator of $Ab$ (tensor by $\_ \otimes_\mathbb{Z} S$ and apply the $\otimes,Hom$ adjunction).
    
    Thus, $[R,V]_\mathbb{Z}$ is an injective cogenerator of $Mod-R$ (tensor by $\_ \otimes_R R$ and use the adjunction).
    
    Consider, now, $T\in Ring$ s.t. $_{S}V_{T} \in S-Mod-T$, and let $U_T$ be an injective cogenerator of $Mod-T$. Up to taking powers, w.l.o.g.\ $V_T \leq U_T$.
    
    Also $_{R}{B^{c}}_{T} \in Mod-T$ for the right $T$-module structure induced by the covariant Hom; so, we can consider the double dual $B^{cc} = [B^{c},U]_T$.
    
    $M$ $B^c$-ML implies its (strict) $B^{cc}$-stationarity. Now, the morphism  $ev:B\xhookrightarrow{} B^{cc}$ is a pure mono, and hence further split by the pure-injectivity of $B$. Thus, $M$ is itself strict $B$-stationary, because $M$ is further strict $def(B^{cc})$-stationary, and $B \in def(B^{cc})$ (see the previous \emph{Lemma} for this closure property).
    \end{proof}

    \subsection{Miscellany of set- and model-theoretical results}
    
    \subsubsection{Filters and Ultrafilters}
    
    We will recall some notions on filters. We refer to the book from M. Prest \cite{Prest-new},3. Whenever an $I$-indexed product of modules occurs, we will freely consider its elements either as $I$-tuples or as maps $i \mapsto x_i \in M_i$ over $I$.
    
    \begin{definition}
    A \textbf{filter} on a cardinal $\lambda$ is a subset $\mathcal{F} \subseteq \mathcal{P}(\lambda)$ that contains $\lambda$ and is closed under finite intersections and super-sets. It is further \emph{proper} in case it does not contain the empty-set (i.e. it is not the whole power set).
    
    Given a filter $\mathcal{F}$ on $\lambda$, it is said to be an \textbf{ultrafilter} if, for each set in the power set of $\lambda$, $\mathcal{F}$ contains either that set or its complement.
    \end{definition}
    
    The following definition provides for a special and interesting class of filters on a poset.
    
    \begin{definition}
    The \textbf{associated filter} to a poset $I$ is the filter generated by the basis $(i\uparrow) := \{j\in I | j\geq i\}$.
    \end{definition}
    
    Now, let's turn to an application to module theory. Consider a family $(M_i)_I \subseteq Mod-R$ indexed by some poset $I$, and take a filter $\mathcal{F}$ on $I$.
    
    By the stated axioms, $\mathcal{F}$ induces an equivalence, say $\sim$, on $\prod_I M_i$:
    
    $$a \sim b \leftrightarrow \{ i \in I | a_i = b_i \} \in \mathcal{F}$$
    
    Indeed, it is clearly reflexive and symmetric, and the transitivity is given by the fact that, given $a \sim b, b\sim c$, $\{i | a_i = b_i \} \cap \{i | b_i = c_i \} \subseteq \{i | a_i = c_i\}$.
    \bigskip
    
    \emph{CLAIM.} Let $M:= \prod_I M_i$. The set-theoretic quotient $M / {\mathcal{F}} := M / \sim$ is a quotient module of $M$, called the \textbf{reduced product} of $(M_i)_I$ by $\mathcal{F}$.
    \smallskip
    
    \emph{Proof.} The set $Z := \{ a \sim 0 \}$ is a submodule of the product module: the closure under intersection and super-sets provides that under addition, while the closure under scalar products is entailed by the upward one. In other words, $\sim$ is compatible with the module operations, and factorising by $\mathcal{F}$ amounts to quotienting out by $Z$.
    \bigskip
    
    \emph{Remark.} If the underlying filter of a reduced product is an ultrafilter, then the obtained module is said to be an \textbf{ultrapower}.
    
    \begin{definition}
    With the previous notation, $\sum_{\mathcal{F}} M := Z$ is called the \textbf{$\mathcal{F}$-product} of $(M_i)_I$, and consists of all the elements of the products whose zero-set belongs to the filter.
    \end{definition}
    
    Let $\pi: \prod_I M_i \rightarrow M / \mathcal{F}$ be the canonical quotient projection.
    
    Then, we prove that the reduced product is actually the direct limit of the 'partial products' (see \cite{Prest-new},3.1)
    
    \begin{lemma}\label{reduced product}
    Let $(M_i)_I \subseteq Mod-R$, and consider a filter $\mathcal{F}$ on $I$. Denote by $M:= \prod_I M_i$.
    
    Any subset $J\subseteq I$ induces the partial product $M_J:=\prod_J M_i$, and hence the partial projections $\pi_{JK}:M_J \rightarrow M_K$. Consider, then, the directed system $(M_J,\pi_{JK} | J \leq K \in \mathcal{F}^{op})$.
    
    Then, the reduced product $M/\mathcal{F} \simeq \lim _{\mathcal{F}^{op}}\prod_J M_j$ is a directed colimit, and the quotient SES is pure-exact:
    
    $$0 \longrightarrow Z \longrightarrow_{*} M \longrightarrow M/\mathcal{F} \longrightarrow 0 $$
    \end{lemma}
    
    \begin{proof}
    Consider the directed system of the 'partial' SES's $(E_J, \phi_{JK}:=(\iota_{JK},id,\pi_{JK})|J\leq K\in \mathcal{F}^{op})$:
    
    $$
    \xymatrix{
    E_{J}: \ar[d]^{\phi_{JK}} & 0 \ar[r] & {\prod_{I\setminus J} M_j} \ar[d]^{\iota_{JK}} \ar[r] & M \ar@{=}[d] \ar[r] & {\prod_J M_j} \ar[d]^{\pi_{JK}} \ar[r] & 0 \\
    E_{K}: & 0 \ar[r] & {\prod_{I\setminus K} M_j} \ar[r] & M \ar[r] & {\prod_K M_j} \ar[r] & 0
    }$$
    
    The horizontal inclusions are virtually including all the elements whose zero-set contains $J$ (resp. $K$).
    
    Hence, the directed limit of the first column yields $Z$: it can be regarded as the co-directed limit on $\mathcal{F}$, that is the union of the co-directed poset consisting of those sets whose elemets have zero-set in $\mathcal{F}$.
    
    So, since $\lim$ is exact on the Grothendieck category $Mod-R$, the limit of the last column must yield the reduced product.
    
    Finally, the purity follows from the 3rd statement of the previous characterization.
    \end{proof}
    
    \subsubsection{pp-formulas}
    
    We will now briefly recall the parallel model theoretical approach to injectivity and Mittag-Leffler modules, as introduced by Prest in \cite{Prest-old},3 and then developed by Rothmaler in \cite{Rothmaler-old}, \cite{Rothmaler-new}. It translates the whole conceptual framework in terms of linear systems, thus stressing how module properties may be regarded as consequences of their inherent linear structure.
    
    In the work of Rothmaler \cite{Rothmaler-old}, ML modules are characterized as positively atomic models, while in \cite{Rothmaler-new} relative ML conditions are introduced. Moreover, it emerges a plain and elegant proof of the fact that the ML character of a directed system of ML modules is completely determined by that of all of its countable subsystems; the same is algebraically proved by Trlifaj and G\"obel in \emph{Lemma} \cite{Trlifaj&Goebel},3.11.
    
    The just recalled result, together with the several algebraic papers by Angeleri H\"ugel, Herbera, Trlifaj, Zimmermann suggests the theoretical equivalence of the two languages. Notice, however, that the model theoretical approach reveals sometimes to be extremely powerful, especially when fruitfully combined with the algebraic formalism, as it may be noticed in the work of Saroch and Stovicek, as well as in the following sections.
    
    An instance of this is the fundamental work of Ziegler relative to Ziegler's Spectrum, that allows us to characterize definable classes via their indecomposable pure-injectives, that do always constitute a skeletally small class, i.e. a proper \emph{set} up to quotienting by isomorphisms between its members. (Ziegler's spectrum is presented in a more intuitive way by Prest in \cite{Prest-new}.)
    
    Moreover, the model theoretical approach to module categories relies on the notion of pp-formula, i.e. on (projections of solutions of) homogeneous linear systems. The significant relevance of pp-formulas is due to a deep result by Braur, Monk et al. (see \cite{Prest-new},A1.1) stating that any formula in the language of $R$-modules has a solution set that is given by a boolean combination of the sets of solutions of pp-formulas of the same finite arity. Furthermore, it makes precise the generic idea that an inherently linear context as that of modules should be depicted solely by means of linear 'tools'. However, observe that, in general, we are not able to consider only free formulas. Indeed, the complete elimination of quantifiers characterizes Von Neumann regular rings (see \cite{Prest-new},2.3.24). Intuitively, this means that in order to be able to understand a linear framework we should generally regard it as immersed in some bigger linear environment.
    
    In what follows, we will mainly attempt an introduction, aiming at giving some context to the notion of definable class, so 'natural' from the model theoretical viewpoint. Finally, we will close this subsection with some comments on the notion of elementary cogenerators.
    
    \begin{definition}
    Let $R\in Ring$, and let's work in $R-Mod$. A \textbf{pp-formula} (i.e. positive primitive formula) over $R$ is an existentially quantified finite system (i.e. conjunction) of $R$-linear equations, say for instance
    
    $$\exists y_0,\dots,y_{l-1}: \quad \wedge_{i<n} \big( \sum_{j<m} r_{ij}x_j + \sum_{k<l} s_{ik}y_k = 0 \big) \qquad r_{ij},s_{ik} \in R$$
    \end{definition}
    
    We will abbreviate them as $\exists \bar{y}(B\bar{x} = A \bar{y})$ or, even more compactly, with the shorthand $A|B\bar{x}$.
    
    For $Mod-R$, we will write $A|\bar{x}B$, transposing columns to rows so as to represent the opposite multiplication. However, in what follows we will work only in $R-Mod$.
    \smallskip
    
    To grasp some intuition, notice that pp-formulas are precisely projections (induced by the existential quantification) of solutions of homogeneous $R$-linear systems.
    
    \begin{definition}
    Given a $n$-pp-formula $\phi$ (where $n$ corresponds to the number of free variables), we indicate by $\phi(N)$ the $n$-tuple in $N^n$ satisfying our formula. Notice that it is an additive subgroup of $N^n$.
    
    Furthermore, we define the following notions:
    
    \begin{itemize}
        \item A \textbf{pp-type} is a collection of formulas of the same finite arity.
        
        \item A \textbf{complete pp-type} in $N$ relative to the $n$-tuple $\bar{a}$ is a pp-type of the form: 
        $$tp_M^{+}(\bar{a}):= \{ \phi \quad | \quad \phi \text{ $n$-place pp-formula for $R$} \quad \wedge \quad \bar{a} \in \phi(N)  \}$$
        
        \item Given a pp-type $\Phi$, the \textbf{realizations} of $\Phi$ in $N$ are the elements of the set $\{ \phi(N) | \phi \in \Phi \}$.
    \end{itemize}
    \end{definition}
    
    We can define a partial order on pp-formulas via inclusion of their realizations, so as to obtain a pp-lattice:
    
    \begin{definition}
    Define a partial order relation over the class of pp-formulas with the same finite arity:
    \vspace{-.2cm}
    
    $$\phi \leq_{N} \psi \quad \leftrightarrow \quad \phi(N) \subseteq \psi(N)$$
    
    We have then an 'absolute' $\leq$ (i.e. independent of $N$) whenever the inclusion holds true for each $N$.
    
    Notice that the conjunction $\wedge$ and the sum $+$, defined as $\phi + \psi := \exists \bar{z}\big(\phi(\bar{x}-\bar{z}) \wedge \psi(\bar{z}) \big)$, induce meet and join operations (respectively).
    
    Call $\sim_{(N)}$ the equivalence relation induced by $\leq_{(N)}$. After quotienting out our class by $\sim$, it becomes a lattice for the previously defined relations. Such a lattice is referred to as the \textbf{pp-lattice} of all modules, and we write $pp_n(Mod-R)$.
    \end{definition}
    
    \emph{Remark.} Notice that, by quotienting out by $\sim$, we obtain a correspondence between each $n$-place formula class and the set of all its realizations.
    \bigskip
    
    \emph{Remark.} With a similar construction we can consider the lattice $pp_n(N)$ over $\{ \phi(N) \hspace{.1cm} | \hspace{.1cm} \text{$\phi$ $n-$place pp-formula} \}$ with meet $\cap$ and join the group-theoretic sum $+$. In this way, then, we obtain a sub-lattice of the lattice consisting of all the additive subgroups of $N^n$. Moreover, notice that our operations are plainly the restriction of those defined for general pp-formulas.
    
    Indeed, more generally, observe that the pp-lattices of $N$ are factor lattices of the corresponding (equivalent) lattices of all $R$-modules by $\sim_N$.
    \bigskip
    
    We can now introduce an equivalent definition of purity in model theoretic terms.
    
    \begin{definition}
    We call \textbf{pp-chain} any descending chain w.r.t. $\leq$ of pp-formulas of the same arity.
    
    Let $l$ associate to an $n$-tuple $\bar{a}$ its length, and consider any SES $0 \rightarrow K \overset{f}{\rightarrow} M \overset{g}{\rightarrow} N \rightarrow 0$. Then, the SES is called \textbf{pure-exact} if the following equivalent conditions are satisfied:
    \begin{itemize}
        \item $f$ pure-mono, i.e. $\forall \phi$ pp-formula, $\forall \bar{a} \in K^{l(\bar{a})}$, $f(\bar{a}) \in \phi(M) \implies \bar{a} \in \phi(K)$;
        
        \item g pure-epi, i.e. $\forall \phi$ pp-formula, $\forall \bar{b} \in \phi(N)$, $b$ has a $g$-preimage in $\phi(M)$.
    \end{itemize}
    \end{definition}
    
    \emph{Remark.} Notice that pp-chains induce descending chains in \emph{one} of the equivalent lattices of \emph{all} modules. Moreover, any pp-chain is closed (up to $\sim$) under meets.
    \bigskip
    
    \emph{Remark.} The definition of purity via pp-formula is exactly that of a pure-mono via matrix subgroups; also that of a pure-epi is the analogous factorisation condition.
    \bigskip
    
    The notion of definable class naturally emerges in this context as a subcategory that allows for a model theoretical description in terms of pp-types, our base structure.
    
    \begin{definition}
    A \textbf{definable subcategory} of $R-Mod$ is a full subcategory whose object class is the model class of a set of pp-implications.
    
    In other words, its class of objects is the class of all models satisfying the set of axioms given by the implications into some $n$-pp-type $\Phi$ together with the axioms of $R-Mod$.
    \smallskip
    
    In particular, given any class $\mathcal{L} \subseteq R-Mod$, we could consider the definable subcategory that it generates, namely $\langle \mathcal{L} \rangle := \{ C\in R-Mod \hspace{.1 cm} | \hspace{.1cm} \phi \leq_{\mathcal{L}} \psi \implies \phi \leq_{C} \psi \}$, that is the meet of all the sets of realizations containing $\mathcal{L}$ (by the compatibility of the order with $+$ on each given $C$, model of $R-Mod$).
    \end{definition}
    
    \emph{Remark.} Given a class $\mathcal{L}$, the definable closure $\langle \mathcal{L} \rangle$ is the \emph{minimal approximation of $\mathcal{L}$ among the classes that are expressible in model theoretic terms}.
    
    Prest proved that the previous property is equivalent to $\langle \mathcal{L} \rangle$ being the smallest sup-class of $\mathcal{L}$ that is closed under arbitrary products, direct limits and pure submodules; surprisingly, it can be constructed from $\mathcal{L}$ by successively considering the closure under the stated properties in the given order.
    \smallskip
    
    Moreover, notice that $\langle \mathcal{L} \rangle = \langle \mathcal{L}' \rangle$ implies that also $\leq_{\mathcal{L}} = \leq_{\mathcal{L}'}$ (by our construction), so that, in particular, $\leq_{\mathcal{L}} = \leq_{\langle \mathcal{L} \rangle}$.
    \bigskip
    
    Finally, let us present a powerful result that characterizes definable classes in terms of their pure-injectives (see again \cite{Prest-new},A1).
    \bigskip
    
    We say that two models of $R-Mod$ are elementarily equivalent provided that they satisfy the same sentences, i.e. formulas (made of pp-formulas, by the definition of our language) with no free variables.
    
    We can take as a more precise definition the following algebraic one; the equivalence of the two is a classical result of Freyne. 
    
    \begin{definition}
    Given $M,N$ $R$-modules, we say that they are \textbf{elementarily equivalent}, and we write $M \sim_{el} N$, provided that $M$ is a pure submodule of an ultrapower of $N$ and vice versa.
    \end{definition}
    
    Moreover, by a theorem of Eklof and Sabbagh, the pure injective hull of a module is elementarily equivalent to it; so, in particular, its double dual is such.
    
    \subsubsection{Freeness}
    
    The current section is concerned with generalizing the notion of 'freeness', that will be intuitively understood as the property characterizing the most 'elementary' and basic pieces of the theory at stake.
    
    For instance, free modules can be viewed as the basic pieces constructing, via direct limits, projective (and then flat) modules, as well as grounding any formulation of dimension via presentation or generation in $Mod-R$.
    
    This meta-mathematical approach will lead to several generalizations of Shelah's compactness principle, as pointed out and formalised by Beke and Rosicky in \cite{Beke&Rosicky}. Towards the end of our project, we will present a theorem stated by Saroch and Stovicek in \cite{Saroch&Stovicek} that can be regarded from this perspective.
    
    At first, let us approach such a viewpoint by the notion of \emph{locally $\mathcal{L}$-free} modules, where the localness is given by the following concept of $\lambda$-dense modules.
    
    \begin{definition}
    Let $\lambda$ be an infinite regular cardinal, and consider a module $M$. A system $\mathcal{S}$ consisting of $(<\lambda)$-presented submodules of $M$ satisfying the following conditions:
    \begin{enumerate}
        \item $\mathcal{S}$ is closed under unions of well-ordered ascending chains of length $<\lambda$.
        
        \item Each subset $X\subseteq M$ with $\# X <\lambda$ is contained in some $N\in \mathcal{S}$.
    \end{enumerate}
    
    is called a \textbf{$\lambda$-dense system} of submodules of $M$.
    \end{definition}
    
    \emph{Remark.} Notice that clearly $M$ is the directed union of these submodules in $\mathcal{S}$.
    
    \begin{definition}
    Let $\mathcal{F}$ be a set of $(<\lambda)-p.$ modules; call them \emph{free}. Denote by $Filt(\mathcal{F})^{<\lambda}$ the class of modules possessing a $(<\lambda)$-$\mathcal{F}$-filtration. A module $M$ is \textbf{locally $\mathcal{F}$-free} provided that it admits a $\lambda$-dense system of free-filtered submodules (i.e. from $Filt(\mathcal{F})^{<\lambda}$).
    \end{definition}
    
    \emph{Remark.} Remark that if $M$ is $(<\lambda)-p.$, then it is locally $\mathcal{F}$-free iff it is already in $Filt(\mathcal{L})^{<\lambda}$, i.e. if it is already 'free-filtered'.
    \bigskip
    
    \emph{Example.} A remarkable example of such an approach to freeness is the following characterization of the class of (flat) Mittag-Leffler modules.
    
    Indeed, \emph{Theorem} \cite{Trlifaj&Goebel},3.14 proves that the class $\mathcal{M}\mathcal{L}$ of Mittag-Leffler modules coincides with that of locally free modules over $\aleph_1$-pure-projectives (i.e. those modules admitting an $\aleph_1$-dense system of pure-projective submodules).
    
    Similarly, \emph{Theorem} \cite{Trlifaj&Goebel},3.19 replaces the class and the free objects in the previous result with $\mathcal{F}\mathcal{L}$ (the class of flat Mittag-Leffler modules) and locally free modules over $\aleph_1$-projectives (respectively). Incidentally remark that, by \emph{Lemma} \cite{Trlifaj&Goebel},3.18, a module is $\aleph_1$-projective iff it is flat and $\aleph_1$-pure-projective, so that the analogy has a very 'natural' origin.
    \bigskip
    
    Moreover, in what follows we will present a construction to refine $\lambda$-dense systems of a given well-behaved module compatibly with a $\kappa$-presenting SES. It heavily employs the machinery that will be introduced at the end of the next section, so we postpone its proof until then.
    
\section{The toolbox}
    
    \subsection{Direct Limits}
    
    \subsubsection{Well-ordered posets}
    
    Let's shift the focus in controlling the size of directed presentations from the dimension of the single terms to the length of the sequence.
    
    \begin{definition}
    For an infinite cardinal $\theta$, a module $M$ is $\theta$-presented provided that it admits a presentation consisting of f.p.\ modules indexed by a poset of cardinality at most $\theta$.
    \end{definition}
    
    Up to further subdividing each 'piece', the just stated definition is plainly equivalent to that considering a directed system consisting of $\theta-p.$ modules. However, it is better-behaved while dealing with ML modules and stationarity conditions, since they are formulated in terms of f.p.\ directed systems.
    
    Moreover, up to continuous modifications, it allows us to deal with well-ordered presentations, as performed in the extremely useful \emph{Lemma} by Iwamura reported in \cite{Saroch-tree}.
    
    \begin{lemma}\label{cofinal subsystem}
    Let $\theta$ be an infinite cardinal, and consider a $\theta$-presented module $M$. Then, $M$ is the directed limit of a \texttt{continuous} and \texttt{well-ordered} system of cardinality $cf(\theta)$ consisting of $(<\theta)-p.$ modules.
    \end{lemma}
    
    \begin{proof}
    Let $\mathcal{M}:= (M_i,f_{ji}|i<j\in I)$ be a f.p.\ directed presentation of $M$ of cardinality w.l.o.g.\ $\theta$.
    
    We propose the proof of \emph{Lemma} \cite{Trlifaj&Goebel},2.14, since it suggests a claim of independent interest.
    \bigskip
    
    \emph{CLAIM.} Any infinite directed set $(I,\leq)$ is a union of a well-ordered chain of its directed sub-posets, each of which has smaller cardinality than $I$.
    
    \emph{PROOF.} For $J\subseteq I$ finite, choose some upper bound of $J$ in $I$, say $j$. Set $J^{\ast}:= J \cup \{j\}$.
    
    Now, assume that $J\subseteq I$ is infinite. We will produce a directed poset $J\subseteq J^{\ast} \subseteq I$ of the same cardinality as $J$.
    
    To do so, let $J_0:=J$, and, given $J_n$, define $J_{n+1}$ by adding to $J$ the upper bounds of any pair of its elements (equiv. of any finite subset). Call $J^{\ast}:= \cup_n J_n$. At each step $\# J_{n+1} = \# \mathcal{P}_f (J_n) = \#J_{n}$, so that $\# J^{\ast} = \# (\cup_{n<\omega} \# J) = \# J$.
    
    We can then inductively define the sought chain $(I_\alpha | \alpha < \theta)$, for a cardinal $\theta$ s.t. $I=(i_\alpha | \alpha < \theta)$.
    
    Set $I_0:=\emptyset$. For $\alpha^+$ successor ordinal, set $I_{\alpha+1}:=(I_{\alpha}\cup \{i_\alpha\})^{\ast}$. Finally, make the chain continuous by positing $I_\alpha:= \cup_{\beta <\alpha}I_{\beta}$ for limit ordinals.
    
    Then, clearly $I=\cup_{\alpha < \theta} I_{\alpha}$.
    
    \bigskip
    
    Now, consider any chain of cardinals $\big(\beta_\gamma | \gamma < cf(\theta)\big)$ realizing the cofinality of $I$. Let's extract a cofinal sub-chain of sub-posets of $I$ indexed by $\beta_\gamma$. Since cardinals are all limit ordinals and $\#I_{\alpha +1} = \# I_\alpha$, we can argue by induction that $\#I_{\beta_{\gamma}} = \beta_{\gamma}$.
    
    Define $M_\gamma := \varinjlim \big( M_i,f_{ji}|i\leq j \in I_\gamma \big)$, and call $\big(m_{\delta,\gamma}:M_\gamma \rightarrow M_\delta | \gamma <\delta <cf(\theta)\big)$ the canonical maps induced by universality. Notice that each $M_\gamma$ is $(<\theta)-p.$
    
    $M=\varinjlim \big( M_\gamma,m_{\delta,\gamma}|\gamma <\delta < cf(\theta)\big)$ is then the sought system; the continuity is given by the equality $\varinjlim_{j} \varinjlim_{I_{\gamma_j}} M_i = \varinjlim_{ I_{\cup_j \gamma_j}} M_i$.
    \end{proof}
    
    \subsubsection{$\lambda$-Continuous Limits}
    
    In this section we will generalize the following 'finite-dimensional' argument to arbitrary cardinality.
    
    Take a directed system $\mathcal{M}$ of modules, and pick some module $N$. Consider a family of maps $(g_i:N\rightarrow M_i)_i$ that is compatible with $(f_{ji})_{j\geq i}$. Call $g := f_i \circ g_i = f_j \circ g_j$, and assume that $Ker (g)$ is f.g.\ Then, also $K_i := g_i (Ker(g)) \leq M_i$ is f.g.\ for each $i\in I$.
    
    Each of the finitely many generators of $K_i$ is sent to $0\in \lim M_l$ by $f_i$, so that there must be some index $j$ s.t. $f_{ji}\circ g_i(Ker(g)) \subseteq M_j$ is already 0 (by \emph{Lemma} \cite{Trlifaj&Goebel},2.2). 
    
    Now, they are finitely many, so we can consider the supremum of such $j$'s, call it $j_{max}$, and conclude that $g_i (Ker(g))$ must already vanish in $M_j$.
    
    Let's call $\aleph_0$-continuity the fact that $I$ is closed under suprema of finite sets. We want to generalize this notion to arbitrary cardinals $\lambda$:
    
    \begin{definition}
    A poset is $\lambda$-continuous if it is closed under suprema of well-ordered ascending chains of length $<\lambda$.
    
    We say that a directed system of modules is \textbf{$\lambda$-continuous} if it is indexed over a $\lambda$-continuous directed poset and $M_{sup(i_\alpha)}=\lim_{\alpha < \tau} M_{i_\alpha}$ for each $(i_\alpha |\alpha < \tau) \subseteq I$ well-ordered ascending chain of length $<\lambda$.
    \end{definition}
    
    In other words, we want to be able to move continuously along the various lattice branches of the indexing poset, and to translate this continuous 'paths' into  paths along modules at the level of the directed system.
    
    Remark that via further specifications of the branch/path structure of the system one may obtain powerful results, as for instance in the case of tree modules (see the generalization proposed by Saroch in \cite{Saroch-tree}).
    \smallskip
    
    Moreover, notice that this stronger continuity notion fulfills our purpose, since it allows us to trivially generalize our argument to a map $g$ that factors through $(f_{ji})_{i\leq j}$ and s.t. $Ker(g)$ is $(<\lambda)-g.$
    
    Incidentally, observe that the factorization condition is always satisfied whenever $N$ is f.p.\ (see \emph{Lemma} \cite{Trlifaj&Goebel},2.7-2.8).
    \smallskip
    \smallskip
    
    \emph{Lemma} \cite{Saroch&Stovicek},1.1 by Saroch and Stovicek shows that continuity conditions are not restrictive, up to modifications of the considered system.
    
    \begin{lemma} \label{continuous colimits}
    Let $M$ be a module, and consider any infinite regular cardinal $\lambda$. Then, $M$ is the directed limit of a $\lambda$-continuous directed system of $(<\lambda)-p.$ modules.
    \end{lemma}
    
    \begin{proof}
    For a presentation of $M$, $R^{(X)} \overset{f}{\rightarrow} R(Y) \rightarrow M \rightarrow 0$, consider the set
    \vspace{-.2cm}
    
    $$I := \{ (x,y) \in \mathcal{P} (X) \times \mathcal{P} (Y) \quad | \quad \# x + \# y < \lambda \quad \wedge \quad f(R^{(x)}) \subseteq R^{(y)} \}$$
    
    that is trivially a directed and $\lambda$-continuous poset under inclusion on both components.
    
    Now, for each $(x,y)\equiv i \in I$, set $M_i := Coker[(f|R^{(x)}):R^{(x)} \rightarrow R^{(y)}]$. Then, $\big( M_i | i\in I \big)$ endowed with the maps induced by the universal property of cokernels is a $\lambda$-continuous directed system with limit $M$.
    \end{proof}
    
    Moreover, the notion of $\lambda$-continuity extends also to filters.
    
    \begin{definition}
    A filter $\mathcal{F}$ on a poset $I$ is $\lambda$-continuous if it admits meets (for the order induced by the reversed inclusion, so intersections) of sets of cardinality $<\lambda$.
    \end{definition}
    
    Remark that $\lambda$-continuity of posets induces that of its associated filter.
    
    \begin{lemma}
    Let $I$ be a $\lambda$-continuous poset. Then, any of its subsets of cardinality $<\lambda$ has an upper bound in $I$. In particular, the associated filter $\mathcal{F}_I$ is $\lambda$-complete.
    \end{lemma}
    
    \begin{proof}
    Enumerate $J=\{ j_\alpha | \alpha < \tau \}$ for some $\tau < \lambda$. Let's construct by induction a well-ordered ascending chain $(k_\alpha | \alpha < \tau) \subseteq I$ s.t. $k_0=j_0$, and $k_\alpha \geq \{ j_\alpha, sup_{\beta <\alpha} k_\beta \}$. The successor step is given by the directedness and the limit one by the $\lambda$-continuity + directedness of $I$. It has supremum $k$ in $I$ by the $\lambda$-continuity, and $k$ is an upper bound for $J$ in $I$, as required.
    
    Now, in particular, consider a subset of $\mathcal{F}_I$ of cardinality $<\lambda$, say $X=(x_\alpha | \alpha < \eta)$ for some $\eta < \lambda$. Then, by our construction, each $x_\alpha$ contains some subset of $I$ of the form $(i_\alpha \uparrow)$. The previous argument yields some $k\in I$ s.t. $k \geq J \equiv (i_\alpha | \alpha < \eta)$, and $(k\uparrow) \subseteq \bigcap X$, so that the latter is in $\mathcal{F}_I$.
    \end{proof}

    \subsubsection{Uniform Factorisation Conditions}
    
    \begin{definition}

    Given a directed system $\mathcal{M}:=(M_i,f_{ji}|i<j\in I)$ and a class $\mathcal{B} \subseteq Mod-R$, we say that $\mathcal{M}$ satisfies the \textbf{uniform factorisation property} provided the existence of a strictly increasing function $s: I \rightarrow I$, independent of $B\in \mathcal{B}$, with the following property: $\forall c\in [M_i,B]_R, \quad c \in H_{f_{s(i),i}}(B) \implies c\in H_{f_{j,i}}(B) \quad \forall s(i) \leq j \in I$.
    \smallskip 
    
    In other words, we are requiring that for each $i\in I$, the index $j=j(B,i)$ whitnessing the ML character of $\mathcal{M}$ does not depend on $B$, but only on $i$. 
    \end{definition}
    
    \emph{Remark.} Observe that the functor $H_v:B \mapsto H_v(B)$ defines a \emph{co-sieve} of $End_R(\_)$-modules.
    
    Indeed, given $v \in [M,M']_R$, $H_v(B) \equiv [M,B] \circ v \leq [M',B]_R$, and for each module $D$, $\forall g \in H_v(D)$, $\forall f:D \rightarrow B$, it clearly holds $f\circ g \in H_v(B)$.
    \smallskip
    
    Moreover, fixed $i_0\in I$, the following functor $J_{i_0}$ defines a Grothendieck topology of cosieves: 
    \begin{equation*}
        \begin{split}
            B & \mapsto J_{i_0}(B):=\{ H_{f_{ji_0}}(B) | (i_0<)j\in I \} \\
            [h:B\rightarrow B'] & \mapsto J(h):= [h_{\ast}:H_{f_{ji_0}}(B) \rightarrow H_{f_{ji_0}}(B')]
        \end{split}
    \end{equation*}
    
    The Uniform Factorisation Condition means precisely that $\mathcal{B}$ has 'thick' elements for the considered topology $J_i$, meaning that the intersection of all the covering sets for each $B$ is eventually stationary after the same index $j=j(i_0)$, independently of $B\in \mathcal{B}$.
    \bigskip
    
    The following reformulation of results 2.3 and 2.5 in \cite{Saroch&Stovicek} (as in \cite{Saroch-tools},4.2) provides for sufficient conditions to obtain the uniform factorisation property.
    
    \begin{lemma}\label{uniform stationarity}
    Assume that $\mathcal{B}$ is a \texttt{filter-closed} class of modules, and that $M \in {^{\perp}\mathcal{B}}$.
    
    Then, $M$ is strict $\mathcal{B}$-stationary, and the indices witnessing such stationarity are uniform over $\mathcal{B}$ (for each considered presentation).
    \end{lemma}
    
    \begin{proof}
    Let $M$ be orthogonal to $\mathcal{B}$, and consider any f.g.\ directed presentation of $M$, say $\mathcal{M}$. Recall that, by the compatibility of the morphisms in $\mathcal{M}$, for each $B\in \mathcal{B}$ and $l\in J$, $Im(h^{B}_{l}) \subseteq Im(h^{B}_{jl}) \quad \forall l\leq j \in J$, so that it always holds $Im(h^{B}_{l})\subseteq \bigcap_{j\leq l} Im(h^{B}_{jl})$. Hence, the strict $B$-stationarity amounts to the remaining inclusion.
    
    Let's show by contradiction the stronger statement about the uniformity of the index, namely that such an inclusion holds uniformly w.r.t. $B\in \mathcal{B}$.
    
    In other words, assume as absurdum that:
    \vspace{-.2cm}
    
    \begin{equation*}
        \mathcal{S}:= \{ l\in I \quad | \quad \big( \forall i \in (l \uparrow)\big) \quad (\exists B_{i}\in \mathcal{B} \quad s.t. \quad \exists g_{i} \in [M_i,B_i]_R \quad : \quad f_{il}^{\ast} (g_i) \not \in H_{f_l}(B_i) \} \neq \emptyset
    \end{equation*}
    
    Complete our data to a $J$-indexed family by setting $B_{i}, g_{i}=0$ for $i<l$, and consider
    
    $$\{[t_{ji}:M_i \rightarrow B_j]:= g_j f_{ji}\}_{i\leq j}, \quad \{ t_{ji}=0 \}_{i>j}$$
    
    \begin{minipage}{.70\textwidth}
    Set $B:= \prod_J B_j$. By means of universality, we canonically obtain the morphism $t$ as in the aside diagram.
    
    Let $\Phi: \oplus_J M_i \rightarrow \lim_J M_i \equiv M$ be the canonical pure-epic.
    \end{minipage}
    \begin{minipage}{.30\textwidth}
    \raggedleft{
    $\xymatrix{
    \oplus_{k\in J} M_k \ar[r]^{t} & B \equiv \prod_{J} B_j \ar@{->>}[d]^{\pi_j}\\
    M_i \ar@{^{(}->}[u]^{\nu_i} \ar[r]^{h_{ji}} & B_j 
    }$
    }
    \end{minipage}
    
    \emph{CLAIM.} Denote by $\mathcal{F}_J$ the associated filter to the directed poset $J$. Then, $t(Ker \Phi) \subseteq \sum_{\mathcal{F}_J} B$.
    \smallskip
    
    \emph{Proof.} First, notice that $\forall i\leq j \in J, \quad \{ k\in J | t_{ki} = t_{kj}f_{ji} \} \in \mathcal{F}_{J}$, since it contains each $k\geq j$ by the compatibility of $\{t_{ji}\}_{i,j}$ induced by that of $\{f_{ji}\}_{j\geq i}$.
    
    Now, explicitly, $Ker \Phi = \langle y:= x-f_{ji}x | x\in M_i, i\leq j \in J \rangle$, and $i,j$ are fixed for each $y$. Consider any $k\geq j \in J$. Then, $\pi_k \circ t(y) = \pi_k \circ t \circ \nu_i (x) - \pi_k \circ t \circ \nu_j \circ f_{ji} (x) = (t_{ki}-t_{kj}\circ f_{ji})(x)$. Hence, each element $y\in Ker \Phi$ has zero-set $\zeta (ty) \in \mathcal{F}_J$, and $Ker\Phi \subseteq \sum_{\mathcal{F}_J}B$.
    
   \vspace{.2cm}
    
    Then, the hypothesis of the Factor Theorem are satisfied, and there is a well-defined morphism $u: M \rightarrow B/\mathcal{F}_{J}$, that in turn induces a morphism of SES's:
    
    \begin{center}
    $\xymatrix{
    0 \ar[r] & {\sum_{\mathcal{F}_{J}} B}  \ar[r] & B \ar[r]^{\rho} & B/\mathcal{F_{J}} \ar[r] & 0 \\
    0 \ar[r] & Ker \Phi \ar@{.>}[u]^{"t"} \ar[r] & \oplus_J M_i \ar[u]^{t} \ar@{->>}[r]^{\Phi} & {\lim_J M_i} \ar@{.>}[u]^{\exists u} \ar[r] & 0
    }$
    \end{center}
    
    Now, $\mathcal{B}$ is filter-closed, so that the $\mathcal{F}$-product of $B$ is still in the class, and $Ext^{1}(M,\sum_{\mathcal{F}_J}B)=0$. Hence, the map $\rho_{\ast}:[M,B] \rightarrow [M, B/\mathcal{F}_J]$ is epic, and there is some $g \in [M,B]$ through which $u$ factors, say $u = \rho \circ g$.
    
    Thus, for each $i\in J$, $\rho g f_i = \rho g \Phi \nu_i = (u \Phi) \nu_i = \rho (t \nu_i)$, that implies $[gf_i-t\nu_i]:M_i \rightarrow \sum_{\mathcal{F}_J}B$ takes values in the kernel of $\rho$.
    
    Since $M_i$ is f.g.\ for each $i\in J$, there is a $(i\leq) j \in J$ s.t. $\pi_k \circ g\circ f_i = \pi_k \circ t \circ \nu_i$ $(\forall j\leq k\in J)$.
    
    Indeed, by \emph{Lemma} \ref{reduced product}, the reduced product can be seen as a direct limit, so, for any generator $x_\alpha \in M_i$, $[gf_i-t\nu_i](x_\alpha) = 0 \in B/\mathcal{F}_J$ implies that $\exists A_\alpha \in \mathcal{F}_J$ s.t. it already vanishes in $\prod_{a\in A_\alpha} M_a$. $M_i$ is f.g.\, so that we can take the maximum $A_{max}:= max\{A_\alpha\}_\alpha \in \mathcal{F}_J$. All the images of the generators vanish in $\prod_{a\in A_{max}} M_a$, and we can pick a suitable $j\in A_{max}$, since the sets in $\mathcal{F}_J$ contain some $(j\uparrow)$.
    
    In particular, $\pi_k \circ t \circ \nu_i = t_{ki} = g_k \circ f_{ki}$, so that we obtain $\pi_k \circ g \circ f_i = g_k \circ f_{ki}$ $(\forall j\leq k\in J)$.
    
    It holds also in the special case $i=l, k=j$, that yields the factorisation $(\pi_j g) \circ f_l = g_j \circ f_{jl}$, that is the sought contradiction.
    \end{proof}
    
    \textbf{Construction.} Notice that the previous proof suggests the following construction:
    
    Fixed $l\in I$, for each $j \geq l$ consider the inverse system
    \vspace{-.2cm}
    
    $$\big(H_i:=H_{f_{ji}}(B_j), \quad u_{ki}:=f^{\ast}_{ik} \quad | \quad k \leq i (\leq j) \in J)$$
    
    and pick the element in its inverse limit $x \in \varprojlim_{i\in J} H_{f_{ji}}(B_j)$ s.t. $x_j = g_j$. Embed it as a map in $\prod_{i \in J} [M_i,B_j]_R$ by putting $0$'s in each component $i > j$.
    
    Do the same for each $j \in J$, and consider the resulting product sequence over $J \times J$, where the second coordinate indexes the codomain.
    
    In this way we obtain the maps $t_{ji}$'s.
    
    Call $\triangle_J$ the lower triangle on $J^2$, i.e. $\{ (j,i) \in J^2 \quad | \quad i\leq j \}$, and make it into a poset via the order $(j,i') \leq (j,i) \in J^2 \quad \leftrightarrow \quad i'\leq i \leq j \in J$.
    
    Now, fix $j_0 \in J$ and remark that the inverse system
    \vspace{-.2cm}
    
    $$\big( H_{f_{ji}}(B_{j_0}),\quad u_{(j,i),(j,i')}:= f^{\ast}_{i,i'} \quad | \quad (j,i') \leq (j,i) \in \triangle_J \big)$$
    
    induces a canonical quotient diagram by means of universality:
    
    \begin{center}
        $\xymatrix{
        0 \ar[r] & {\varprojlim_{i \leq j} H_{f_{ji}}(B_{j_0})} \ar[r] \ar@{->>}[d] & {\prod_{i\leq j} H_{f_{ji}}(B_{j_0})} \ar[r]^{\sigma_0} \ar@{->>}[d]^{\rho_0} & {\prod_{i\leq j}H_{f_{ji}}(B_{j_0})} \ar[r] \ar@{->>}[d] & 0 \\
        0 \ar[r] & {\varprojlim_{i\leq j_0} H_{f_{j_0i}}(B_{j_0})} \ar[r] & {\prod_{i\leq j_0} H_{f_{j_0i}}(B_{j_0})} \ar[r]^{\sigma^{\leq j_0}_0} & {\prod_{i\leq j_0}H_{f_{j_0 i}}(B_{j_0})} \ar[r] & 0
        }$
    \end{center}
    
    where $\sigma:= \pi_{i}-f^{\ast}_{ji}\pi_{j}$ is the canonical map in the presentation of the inverse limit, and the vertical maps are obtained by truncation at $j_0$.
    
    Observe that the embedding $\epsilon$ of $\prod_{i\leq j_0} H_{f_{j_0i}}(B_{j_0})$ into $\prod_{i\leq j} H_{f_{ji}}(B_{j_0})$ by completing the sequences with $0$'s provides a section of $\rho_0$ that, however, does not extend to a section of the whole diagram of SES.
    
    Taking the product over all the $j_0$'s in $J$, $t$ is plainly the image under $\epsilon$ of \emph{any choice} of a map in $\prod_{j_0 \in J} \varprojlim_{i \leq j_o} H_{f_{j_oi}}(B_{j_0})$.
    
    It always holds that $t(Ker \Phi) \subseteq \sum_{\mathcal{F}_J}B$, and our hypotheses on $\mathcal{B}$ provide the factorization (and hence the strict stationarity).
    \bigskip
    
    \emph{Remark.} In case $J$ is countable, then the associated filter is just the Frech\'et filter, and the filter product amounts to the direct sum; thus we can weaken the condition of $\mathcal{B}$ being filter-closed to the closure under countable direct sums.
    \bigskip
    
    Moreover, notice that the following generalization (see \emph{Lemma} \cite{Saroch&Stovicek},2.3) states that, under some non-restrictive continuity assumption on $\mathcal{M}$, we could achieve the uniform factorisation property up to a cofinal $\lambda$-continuous subposet of $I$.
    
    \begin{lemma}\label{generalized uniformity}
    Let $\mathcal{B} \subseteq Mod-R$ be a \texttt{filter-closed} class. Consider any uncountable regular cardinal $\lambda$ and a $\lambda$-continuous directed system of modules $\mathcal{M}$ consisting of $(<\lambda)-g.$ modules with colimit $M$.
    
    If $M\in {^{\perp}\mathcal{B}}$, then there is a $\lambda$-continuous cofinal subposet $J\subseteq I$ s.t. (for each $j\in J$) $[M_j,B]_R \subseteq H_{f_j}(B)$ uniformly w.r.t. $B\in \mathcal{B}$. (Hence the equality holds.)
    \smallskip
    
    We will denote this by $[M_J, \mathcal{B}] = H_{f_J}(\mathcal{B})$.
    \end{lemma}
    
    \begin{proof}
    Our absurdum is again that the set
    
    $$\mathcal{S}:= \{ i\in I \quad | \quad \exists B_i \in \mathcal{B}, \quad \exists g_i \in [M_i,B_i]_R \quad s.t. \quad g_i \not \in H_{f_i}(B_i) \}$$
    
    has non-trivial intersection with each other cofinal $\lambda$-continuous subposet of $I$.
    
    Repeat the previous construction until the point in which we leverage on the finite assumptions proper of the previous setting.
    
    Notice that, carrying on basically the same argument as in the finite case and applying the observation pointed out in the \emph{Preliminaries}, one gets that the set $K_i := \{ k\in I \quad | \quad \pi_k\circ g \circ f_i = \pi_k \circ h \circ \nu_i \} \in \mathcal{F}_I$.
    
    Perform the same substitution and define
    
    $$J := \{i\in I \quad | \quad \big(\forall k\in (i\uparrow)\big) (\pi_k \circ g \circ f_i = g_k \circ f_{ki}) \}$$
    
    In particular, evaluating at $k = i$, $g_i \in H_{f_i}(B_i)$ for each $i\in J$. In other words, $J\cap \mathcal{S} = \emptyset$.
    
    So, we are left to show that $J$ is $\lambda$-continuous and cofinal to $I$.
    
    The continuity of $\mathcal{M}$ forces that of $J$ (because compatibility forces the defining formulas to behave continuously w.r.t. $i$ as well).
    
    Hence, let's show the cofinality part. Fixed $i$, $K_i \in \mathcal{F}_i$ yields (since $\mathcal{I}$ is generated by the sets of the form $(j\uparrow)$) an increasing function $s:I\rightarrow I$ s.t. $(s(i) \uparrow) \subseteq K_i$, i.e. the formula $\big(\forall k \geq s(i)\big)(\pi_k\circ g \circ f_i = g_k \circ f_{ki})$ holds true.
    
    Now, fix any $i' \in I$, and consider the sequence $(j_n)_{n < \omega}$ defined by $j_0 := i'$ and $j_{n+1}:= s(j_n)$. Call $j := sup (j_n)$.
    
    One has, then, both $j\geq i'$ and $j\in J$: by the $\aleph_1$-continuity of $\mathcal{M}$, $\big(\forall k\geq sup j_{n+1} \big)(\pi_k\circ g \circ f_{sup j_n} = g_k \circ f_{k,sup j_n})$.
    \end{proof}
    
    Finally, let's apply the previous results to derive some closure properties for $Ext^1$ conditions.
    
    \begin{lemma}\label{closure of Ext conditions}
    Let $\mathcal{G}\subseteq Mod-R$, and take any $M \in {^{\perp}\mathcal{G}}$ c.p.\ TFAE:
    \begin{enumerate}
        \item $M\in {^{\perp}D}$ for any copy of a pure submodule of a product of modules from $\mathcal{G}$; i.e. w.l.o.g.\ $\mathcal{G}$ is closed under pure-submodules of its products.
        
        \item $M$ is $D$-stationary for any $D\in def(\mathcal{G})$, the definable closure of $\mathcal{G}$, i.e. w.l.o.g.\ $\mathcal{G}$ is definable.
        
        \item $M \in {^{\perp}D}$ for any copy of a countable directed sum from $\mathcal{G}$; i.e. w.l.o.g.\ $\mathcal{G}$ is closed under countable direct sums.
        
        \item $M$ is $\mathcal{G}$-stationary.
    \end{enumerate}
    \end{lemma}

\begin{proof}
$1) \implies 3).$ Direct sums are pure submodules in the corresponding direct products.

\bigskip
$3) \implies 2).$ It is almost provided by the previous \emph{Lemma}:

For $M$ c.p.\, we can regard it as the countable direct limit of a directed system of f.p.\ modules, and provided that $\mathcal{G}$ is closed under countable direct sums (as assumed in $2)$), our \emph{Lemma} assures the existence of a strictly increasing function $s: \omega \rightarrow \omega$, independent of $B\in \mathcal{G}$, witnessing the uniform factorization property, i.e. $\forall c\in [M_i,B]_R, \quad c \in H_{f_{s(i),i}}(B) = H_{f_{i}}(B) \implies c\in H_{f_{j,i}}(B) \quad \forall s(i) \leq j \in \omega$.
\smallskip

\emph{CLAIM.} The previous factorization holds uniformly over the all of $def(\mathcal{G})$. (See \emph{Lemma} \cite{Saroch&Stovicek},2.6.)
\smallskip

\emph{PROOF.} By the work of Crawley-Boevey (see \emph{Preliminaries}), the vanishing locus of a set of coherent functors is closed under definable closure, so we aim at proving that the previous factorisation property can be written in terms of some set of coherent functors.

To this extent, as in \emph{Prop.} \cite{Saroch&Stovicek},2.7, consider the family of functors $(F_{j}(\_):=H_{f_{s(i),i}}(\_) / H_{f_{j,i}}(\_) | i\leq s(i) \leq j < I)$ given by the family of SES's $\xymatrix{
0 \ar[r] & H_{f_{j,i}}(\_) \ar[r]^{\subseteq} & H_{f_{s(i),i}}(\_) \ar[r] & F_j(\_) \ar[r] & 0
}$. They are clearly coherent functors, because arbitrary products and direct limits of the argument commute with the SES structure for each $j$ (the latter since the $M_i$'s are f.p.). Moreover, a class of modules admits the previous factorisation via the function $s$ iff it is in the vanishing locus of our family, so we are done.

Finally, notice that $def(\mathcal{G})$ is clearly filter-closed. Indeed, definable classes are closed under products and pure submodules, and hence also under filter products, since the latter sits in the first position of pure-exact SES's with middle terms direct products from $\mathcal{G}$ (see \emph{Lemma} \ref{reduced product}).
\bigskip

$2) \implies 1).$ Write $D\subseteq_{\ast} \prod_k B_k$ for some modules in $\mathcal{G}$. Any c.p.\ module can be seen as the direct limit of a countable chain of f.p.\ $(M_i | i<\omega)$. Consider, then, the ML (by assumption) system induced by $[\_,D]_R$.

Consider the pure-SES $0\rightarrow D \rightarrow_\ast \prod B_k \rightarrow \prod B_k / D \rightarrow 0$.

It induces via $[M,\_]_R$ the exact sequence $[M,\prod B_k]_R \rightarrow [M,\prod B_k / D]_R \rightarrow Ext^1 (M,D) \rightarrow 0 = Ext^1(M,\prod B_k)$, where the vanishing of the last term is because $Ext^1$ commutes with arbitrary products in the second component, and because $M\in {^{\perp}\mathcal{G}}$.

Moreover, the given SES also yields the exact sequence $\varprojlim [M_i,\prod B_k]_R \rightarrow \varprojlim [M_i,\prod B_k / D]_R \rightarrow 0 = \varprojlim^{1}[M_i,D]_R \rightarrow 0 = \varprojlim^1 [M_i,\prod B_k / D]_R$.

And the last vanishings are implied by the equivalence between the ML property and the vanishing of the derived functor $\lim^1$.

Bringing the inverse limits inside the contravariant Hom-functor, the 5-Lemma implies that the latter is isomorphic to the former sequence, and it forces $Ext^1(M,D)\simeq 0$, as required.
\bigskip

$4) \implies 2).$ Provided the uniform factorisation condition, it is just condition $2).$ of a result by Angeleri H\"ugel and Herbera, namely \emph{Cor.} \cite{H-subgroups},3.9 (see \emph{Lemma} \ref{closure properties}).
\end{proof}

\subsection{Stationarity and Pure-Injectivity}

Another application of the previous tools is the following \emph{Lemma} relating stationarity conditions w.r.t. a $\Sigma$-pure-injective module.

\begin{lemma}
Let $B\in Mod-R$. TFAE:
\begin{enumerate}
    \item $B$ $\Sigma$-pure-injective.

    \item All modules are strict $B$-stationary.
    
    \item All modules are $B$-stationary.
    
    \item All countably presented modules are $B$-stationary.
\end{enumerate}
\end{lemma}

\begin{proof}
$1) \iff 2).$ It is \emph{Theorem} \cite{zimmermann},3.8.
\bigskip

$4) \implies 3).$ The following result, namely \emph{Proposition} \cite{H-subgroups},3.10, holds true:
\smallskip

Let $(M_i,f_{ji}|i<j\in I)$ be a directed system in $Mod-R$, and let $B\in Mod-R$. TFAE:
\begin{itemize}
    \item For each countable chain $(i_n)_{n<\omega} \subseteq I$, $(F_{i_n},f_{i_{n+1},i_n}|n<\omega )$ is $B$-stationary.
    
    \item For each countable chain $(i_n)_{n<\omega} \subseteq I$, $\big( H_{f_{i_{n+1},i_n}}(B) |n<\omega \big)$ is stationary.
\end{itemize}

Each of them imply that the original direct system is $B$-stationary.
\smallskip

\emph{PROOF.} The equivalence of the two statements is just a matter of definitions. Let's prove the implication. Assume by contradiction that there is some $i \in I$ s.t., for each $(i\leq )j \in I$, there is some $i\leq l$ for which $H_{f_{ji}}(B) \not \supseteq H_{f_{li}}(B)$. We will inductively construct a countable chain that does not satisfy the first of the equivalent conditions.

Set $i_0:= i$. Given $(i_k)_{k\leq n}$ s.t. $H_{f_{i_k,i_0}}(B)$ strictly contains $H_{f_{i_{k+1},i_0}}(B)$, the assumption applied to $j=i_n$ yields some $i_0 \leq l \in I$ s.t. $H_{f_{i_n,i_0}}(B)$ strictly contains $H_{f_{l,i_0}}(B)$. Consider any $i_{n+1} \in I$ s.t. $i_{n+1} \geq i_n,l$, so that $H_{f_{i_{n+1},i_0}}(B) \subseteq H_{f_{l,i_0}}(B)$ and the latter is a proper subset of $H_{f_{i_n,i_0}}(B)$, as required.
\bigskip

Now, take any module $M$ and consider any of its presentations in f.p.\ modules. By the previous Proposition, such a directed system is $B$-stationary iff each of its countable subsystems is such, and the latter condition, since $M$ is arbitrary, is properly equivalent to condition $4)$.
\bigskip

$3) \implies 2).$ We proved that $B$-stationarity can be strengthened to $def(B)$-stationarity (see \emph{Lemma} \ref{closure properties}).

Notice that $M \sim_{el} N$ implies that $M$ is a pure submodule of an ultrapower of $N$ (that in turn is a pure-epi image of $N^I$ for some set $I$ by our result on the presentation of reduced product), and we also showed that $def(B)$ is closed under pure-epi images. Hence, we can assume $def{B}$ to be closed under elementary equivalence.

Now, a result by Eklof and Sabbagh states that the pure-injective envelope of a module is elementarily equivalent to it, so that we obtain also the $\mathcal{P}\mathcal{E}(B)$-stationarity of $Mod-R$.

Finally, as proved in the \emph{Preliminaries}, stationarity and strict stationarity w.r.t. a pure-injective module are equivalent notions, and hence $Mod-R$ is strict $\mathcal{P}\mathcal{E}(B)$-stationary. But then notice that the embedding $\iota:=[B\xhookrightarrow{} \mathcal{P}\mathcal{E}(B)]$ induces $\iota_{\ast}: H_v(B) \xhookrightarrow{} H_v(\mathcal{P}\mathcal{E}(B))$ ($\forall v$), that entails the strict $B$-stationarity (for the same index function) via the commutative diagram:

\begin{center}
    $\xymatrix{
    H_{f_{s(i)i}}(B) \ar@{^{(}->}[rr]^{\iota_\ast} && H_{f_{s(i)i}}(\mathcal{P}\mathcal{E}(B)) \\
    H_{f_{i}}(B) \ar@{^{(}->}[rr]^{\iota_\ast} \ar@{^{(}->}[u]^{\subseteq} && H_{f_{i}}(\mathcal{P}\mathcal{E}(B)) \ar@{=}[u]
    }$
\end{center}

Thus, $B$ must be $\Sigma$-pure-injective.
\bigskip

Finally, to complete the circle, notice that e.g.\ $2) \implies 4)$.
\end{proof}

\emph{Remark.} Notice that the implication $4)\implies 3)$ provides an alternative proof of the aforementioned countable character of ML conditions.

\subsubsection{Elementary Cogenerators}

We will now introduce the notion of elementary cogenerators.

\begin{definition}
A \texttt{pure-injective} module $E$ is an \textbf{elementary cogenerator} if every pure injective direct summand of any module elementarily equivalent to $E^{\aleph_0}$ is a direct summand of a direct product of copies of $E$.
\end{definition}

The 'cogenerating' character of these modules will be clarified by the following \emph{Lemma}.

To grasp some intuition, indeed, they will supply for a 'canonical' class of pure-injective cogenerators of sufficiently rich classes (very close to be definable), where 'canonical' is understood from the model theoretical perspective. And the 'cogenerating' character amounts to the fact that each element of such rich classes can be purely embedded in some power of the corresponding elementary cogenerator. In this way, then, we will then also obtain a class of 'canonical' pure-injective cogenerators of $Mod-R$.
\bigskip

\emph{Remark.} Prest proved that each module $M$ admits an elementary cogenerator that is elementarily equivalent to it (see \emph{Cor.} \cite{Prest-new},9.36).

\begin{lemma}\label{canonical cogenerator}
Let $\mathcal{B}$ be a class of modules that is closed under \texttt{direct products} and \texttt{direct limits}.

Then, $\mathcal{B}$ contains a pure-injective module $C$ s.t. each module $B\in \mathcal{B}$ embedds purely in the product of some copies of $C$. In other words, for each $B\in \mathcal{B}$ there exists some set $I$ s.t. $B\subseteq_{\ast} C^I$.
\smallskip

Moreover, if $\mathcal{B}$ contains a generator, then $C$ is a cogenerator for $Mod-R$.
\end{lemma}

\begin{proof}
$\mathcal{B}=\lim \mathcal{B}$ provides for the closure under direct summands, and hence also under elementary equivalence to pure-injectives.

Indeed, any pure-injective module $A$ elementarily equivalent to some $B \in \mathcal{B}$ is actually a direct summand of an ultrapower of $B$. The latter is a direct limit (see \emph{Lemma} \ref{reduced product}), so that it is still in $\mathcal{B}$, and also $A$ must be in our class.

Notice, in particular, that $\mathcal{B}$ is closed under taking pure-injective hulls and double duals (see the section on pp-formulas in the \emph{Preliminaries}).

By the previous \emph{Remark}, for each module in $\mathcal{B}$ we can consider a (pure-injective) elementary cogenerator within the class that is elementarily equivalent to the given module.

So now, by the Ziegler's spectrum \emph{Theorem}, we can choose a representative \texttt{set} $\mathcal{S}$ of elementary cogenerators for $\mathcal{B}$.

Call $C:= \prod \mathcal{S}\in \mathcal{B}$ its direct product, that is still pure-injective. Moreover, for each $B\in \mathcal{B}$, there is some set $I$ s.t. $B\subseteq_{\ast} C^I$ up to iso.

Indeed, $B^{cc} \in \mathcal{B}$ and there is some elementary cogenerator $E\in \mathcal{S}$ s.t. $B^{cc} \leq_{\oplus} (B^{cc})^{\aleph_0} \sim_{el} E^{\aleph_0}$.
The latter holds because $B^{cc} \leq_{\oplus} (B^{cc})^{(\aleph_0)} \subseteq_{\ast} (B^{cc})^{\aleph_0}$, and hence $B^{cc} \leq_{\oplus} (B^{cc})^{\aleph_0}$ by the pure-injectivity; moreover, $(B^{cc})^{\aleph_0} \sim_{el} B^{\aleph_0} \sim_{el} E^{\aleph_0}$.
Now, since $E$ is an elementary cogenerator, $B^{cc} \leq_{\oplus} E^I \leq_{\oplus} C^I$. Hence, $B\subseteq_{\ast} B^{cc}$ implies $B\subseteq_{\ast} C^I$.

Finally, clearly if $B\in \mathcal{B}$ is a cogenerator for $Mod-R$, then $C$ cogenerates $Mod-R$ as well.
\end{proof}

Given a pure-injective cogenerator $C$, the following \emph{Lemma} provides us with a 'rich' (i.e. $\lambda$-dense) supply of c.p.\ (strict) $C$-stationary modules. We will consider only the $\aleph_1$ case, although (as Saroch incidentally remarks) it can be generalized to higher cardinals $\lambda$ by adapting the theorem characterizing relative ML modules to the case of $\lambda$-directed colimits.

\begin{lemma}\label{rich C-stationary}
Let $C$ be a pure-injective module which cogenerates $Mod-R$, and consider any (strict) $C$-stationary module $M$. Then, there is an $\aleph_1$-dense system $\mathcal{L}$ of (strict) $C$-stationary submodules of $M$ s.t. $Ext^1 (M/N,C)=0$ (i.e.
$[M,C]_R \rightarrow [N,C]_R$ is epic) for each $N$ directed union of modules from $\mathcal{L}$.
\end{lemma}

\begin{proof}
Let $Q:= C^c$. We will use several results from the paper \emph{"Mittag-Leffler conditions on modules"} by Angeleri H\"ugel and Herbera, all recalled in the previous sections.

A module is strict $C$-stationary iff $Q$-ML (see \emph{Lemma} \ref{ML vs stationarity}). Then, by the characterization of relative ML provided by Facchini and Azumaya (as in \emph{Theorem} \cite{H-subgroups},5.1/4), a presentation of $M$ consisting of f.p.\ modules forms a set $\mathcal{S}$ s.t. $M\in \lim \mathcal{S}$ and for each $X\subseteq M$ with $\#X \leq \aleph_0$, there is a c.p.\ $Q$-ML $N\in \lim \mathcal{S}$ and there exists a map $v:N\rightarrow M$ s.t. $X\subseteq v(N)$ and $v \otimes Q$ is mono.

Consider a collection of countable generating sets for $M$ (e.g.\ induced by the presentation) whose union generates the whole of $M$. Starting from such a family, consider the corresponding c.p.\ $Q$-ML modules containing them, as provided by an application of the previous property. Iterate this procedure to successive generating sets, so as to obtain an $\subseteq$-directed set $\mathcal{F}$ consisting of c.p.\ $Q$-ML submodules of $M$ satisfying the axioms of $\aleph_1$-density.

Moreover, remark that the map $v$ is actually a mono, since $v\otimes Q = [[v,C]_R,\mathbb{Q}/\mathbb{Z}]_\mathbb{Z}$ is a mono obtained by a double dualization of $v$ via a cogenerator and an injective cogenerator, respectively.

We can, now, gradually extend $\mathcal{F}$ by adding countable directed unions (that are necessarily well-ordered).

Notice that for each newly added c.p.\ $N\leq M$, the latter inclusion is preserved by $(\_ \otimes Q)$, because $Tor_1 (\_,Q)$ commutes with $\lim$.

Hence, the limit is still ML, since it is the directed union of some c.p.\ $Q$-ML submodules of $M$ whose inclusion is preserved by $(\_ \otimes Q)$ (see \emph{Cor.} \cite{H-subgroups},5.2).

After a countable number of iterations, $\mathcal{F}$ becomes the desired $\aleph_1$-dense system $\mathcal{L}$ of strict $C$-stationary submodules of $M$ (because all $Q$-ML).

Finally, remark that any directed union of modules from $\mathcal{L}$ satisfies our condition on the vanishing of $Ext^1$, since it is implied by the injectivity of the involved tensors, which is in turn preserved by taking directed unions (that are directed limits).
\end{proof}

Moreover, let us present a similar statement for a general cardinal $\lambda$. Apart from their independent interest, these two existence results will be heavily used in tandem in the proof of the main theorem.

\begin{lemma}\label{initial step}
Let $\mathcal{B}$ be a \texttt{filter-closed} class of modules that \texttt{cogenerates} $Mod-R$. Then, for each $\lambda$ uncountable regular cardinal and for each $M\in {^{\perp}\mathcal{B}}$, there is some $\mathcal{C}_\lambda$ $\lambda$-dense system of submodules of $M$ s.t. $M/N \in {^{\perp}\mathcal{B}}$ for each $N \in \mathcal{C}_\lambda$.
\end{lemma}

\begin{proof}
The proof grounds on several Lemmas.

By \emph{Lemma} \ref{continuous colimits}, let's consider a $\lambda$ continuous presentation of $M$, say $\mathcal{M}$, consisting of $(<\lambda)-p.$ modules.

By \emph{Lemma} \ref{generalized uniformity}, there is some $J\subseteq I$ $\lambda$-continuous and cofinal s.t. $[M_J,\mathcal{B}] = H_{f_J}(\mathcal{B})$.

Being cofinal, $J$ is still directed, and we can consider w.l.o.g.\ $\mathcal{M}$ to be indexed over $J$, because we are interested in its limit.

Since $\mathcal{B}$ is cogenerating, $f_j$ must be a mono for each $j\in J$: for each $j$ there is an inclusion of $M_j$ in some direct sum of copies of elemets of $\mathcal{B}$ (the Fr\'echet-filter-product of $\mathcal{B}$); the latter is still in the filter-closed $\mathcal{B}$, so that the inclusion must factor through $f_j$.

Hence, w.l.o.g.\ $M_j \subseteq_{f_j} M$, as well as $M_i \subseteq_{f_{ji}} M_j$. Then, call again $\mathcal{M}:= \{M_j | j\in J\}$, and write $\bar{\mathcal{M}}$ for its closure under the directed union of well-ordered (not necessarily ascending) subchains of length $<\lambda$.

Remark that, although $\mathcal{M}$ is $\lambda$-continuous, by adding all the unions of well-ordered (not necessarily ascending) subchains we could actually obtain a strictly bigger class.

However, it still remains $\lambda$-continuous as well as directed, since each subset of $\bar{\mathcal{M}}$ of cardinality $<\lambda$ has an upper bound in $\bar{\mathcal{M}}$.

Thus, we can regard $\bar{\mathcal{M}}$ as a $\lambda$-continuous directed system indexed by $\bar{\mathcal{M}}$ itself.

Let's apply again \emph{Lemma} \ref{generalized uniformity} to extract a new $\lambda$-continuous cofinal subposet $\mathcal{C}_\lambda$ out of $\bar{\mathcal{M}}$ s.t. (for each $N\in \mathcal{C}_\lambda$) $[N,\mathcal{B}]_R$ factors through some canonical colimit morphism $\nu_N$.

The morphisms $(\nu_{M_j} | j\in J)$ induce by means of universality a canonical morphism $\nu: M = \lim_J \mathcal{M} \rightarrow \lim \bar{\mathcal{M}}$.

Hence, in particular, we have factorizations $[M_j,\mathcal{B}]_R = H_{\nu_{M_j}}(\mathcal{B)}$ through $\lim \bar{\mathcal{M}}$, that can be extended to factorizations $[M,\mathcal{B}]_R = H_{\nu}(\mathcal{B)}$.

In other words, for each $N \in \mathcal{C}_\lambda$, we have an epimorphism $[M,\mathcal{B}]_R \twoheadrightarrow{} [N,\mathcal{B}]_R$ that forces $Ext^1 (M/N, \mathcal{B}) = 0$, because $Ext^1 (M, \mathcal{B})$ already vanishes by assumption.

We are then left to show that $\mathcal{C}_\lambda$ is actually $\lambda$-dense.

$\mathcal{C}_\lambda$ is $\lambda$-continuous in $\bar{\mathcal{M}}$, so it is clearly closed under unions of well-ordered ascending chains of length $<\lambda$.

Moreover, the cofinality part yields $\cup \mathcal{C}_\lambda = M$. Then, consider any $X \subseteq M$ of cardinality $<\lambda$; for a generic element $x\in X$ there is some $N\in \mathcal{C}_\lambda$ s.t. $x\in N$. Consider the collection of those submodules, say $\mathcal{N}$, and observe that there must be some upper bound $\bar{N} \in \mathcal{C}_\lambda$ for $\mathcal{N}$. For such a module it holds $X\subseteq N$, as required.
\end{proof}

\emph{Remark.} Notice, however, that the generalized \emph{Lemma} is effectively weaker than the previous one. Indeed, non-providing a concrete algebraic construction, the $Ext^1$-vanishing condition restricts to those directed unions that remain in the $\lambda$-dense system; moreover, unless the class $\mathcal{B}$ is hereditary, we cannot directly infer that the obtained dense system consists of $C$-stationary modules for e.g.\ the 'canonical' elementary cogenerator for $\mathcal{B}$.

Extending the previous procedure to arbitrary unions would amount to the fact that $M/U\mathcal{N}$ is $C$-stationary for any $\mathcal{N} \subseteq \mathcal{C}_\lambda$, but the class of $C$-stationary modules need not be closed under inverse limits.

\subsubsection{Refinements of dense systems}

Let us now present the aforementioned construction to refine $\lambda$-dense systems to make them compatible with presentation SES's.

\begin{lemma}\label{refinement of dense systems}
Let $M$ be a $\kappa-p.$ module with presentation $0\rightarrow K \rightarrow R^{(\kappa)} \overset{f}{\rightarrow} M \rightarrow 0$.
    
Consider a cotorsion pair $(\mathcal{A},\mathcal{B})$, where $\mathcal{B}$ is a class \texttt{closed under $\lim$} (and hence also filter-closed). Assume, further, that $M\in \mathcal{A}$.

Choose any uncountable regular cardinal $\lambda$, and consider a $\lambda$-dense system $\mathcal{C}_{\lambda}$ consisting of submodules of $M$ s.t. $M/N\in {^{\perp}\mathcal{B}}$ for each $N\in \mathcal{C}_\lambda$. Let's call it simply $\mathcal{C}$.
\smallskip

Then, we can assume $\mathcal{C} \subseteq \mathcal{A}$ provided that $\mathcal{C}$ consists of c.p.\ modules.
\end{lemma}

\begin{proof} We will almost explicitly construct a refinement of $\mathcal{C}$ that is in $\mathcal{A}$. Notice that we will need the c.p.\ hypothesis only at the end of our proof.
\bigskip

\emph{CLAIM.} w.l.o.g.\ $\mathcal{C}$ is compatible with the given presentation, meaning that we can assume that each $N\in \mathcal{C}$ is of the form $f(R^{(x_n)})$ for some $x_n \subseteq \kappa$ with $\# x_n < \lambda$.
\smallskip
    
\emph{PROOF.} Consider $\mathcal{C} \cap \{f(R^{(x)}) | x\subseteq \kappa,  \#x <\lambda \}$:
    
For each $N\in \mathcal{C}$, let $x_N$ be an $R$-generating system for $f^{-1}(N)$ (that might not be itself free). We then have an epi $R^{(x_N)} \twoheadrightarrow f^{-1}(N) \overset{f}{\twoheadrightarrow} N$, so that $x_N$ is a generating set for $N$, and w.l.o.g.\ $\# x_N < \lambda$. Since $M = \cup C$, we can assume that $R^{(\kappa)} = \cup_{N\in \mathcal{C}} R^{(x_N)}$, and we are done.
\bigskip
    
\emph{CLAIM.} Let $C$ be the pure-injective cogenerator provided by \emph{Lemma} \ref{canonical cogenerator}. Then, $K$ is strict $C$-stationary.
\smallskip
    
\emph{PROOF.} The right class of a cotorsion pair is always closed under direct products by the properties of $Ext^1$; so, whenever it is further closed under directed limits, one obtains also the closure under double duals, as shown in the first part of the proof of \emph{Lemma} \ref{canonical cogenerator}. Hence, for each $B\in \mathcal{B}$, it holds $(B^{(I)})^{cc}\in \mathcal{B}$ for each set $I$; so, in particular, $M\in {^{\perp}(B^{(I)})^{cc}}$, and we can apply \emph{Lemma} \ref{stationarity along ses}, to obtain the strict $C$-stationarity of $K$ whenever $R^{(\kappa)}$ has such property.

Finally, $R^{(\kappa)}$, being free, is ML and a fortiori $C^c$-ML, so that we can conclude by \emph{Lemma} \ref{ML vs stationarity}.
\bigskip

Therefore, we can apply \emph{Lemma} \ref{rich C-stationary} and consider a $\lambda$-dense system $\mathcal{L}$ of strict $C$-stationary submodules of $K$ in $^{\perp}C$ with quotients $K/L \in {^{\perp}C}$ for each $L\in \mathcal{L}$.
\smallskip

\emph{Remark.} We will need only those properties of the system that would have been entailed by the weaker construction of \emph{Lemma} \ref{initial step}. However, without the hereditariness of the cotorsion pair we cannot infer that $K\in \mathcal{A}$. Hence, we needed to prove the $C$-stationarity of $K$ and then apply the stronger \emph{Lemma} \ref{rich C-stationary}.
However, as remarked by the author, the last result can be extended to any uncountable regular cardinal $\lambda$ (up to refine our notion of direct limit), so that here we still do not need to specialize $\lambda = \aleph_1$.
\bigskip

Now, define $\mathcal{K}:= \{ Ker(f|R^{(x_N)}) | N\in \mathcal{C} \}$, and observe that $\mathcal{K} \cap \mathcal{L}$ is still $\lambda$-dense.

In other words, we considered a system $\mathcal{K}$ induced by $\mathcal{C}$ in the form of a $\lambda$-system of SES's, namely 
$$\{ 0 \rightarrow Ker(f|R^{(x_N)}) \rightarrow R^{(x_N)} \rightarrow f(R^{(x_N)}) \equiv N \rightarrow 0 \quad | \quad N\in \mathcal{C} \}$$

or (equivalently) indexed by $Ker(f|R^{(x_N)}) \in \mathcal{K}$.

From this we extract a $\lambda$-dense subsystem, this time indexed by $\mathcal{K} \cap \mathcal{L}$, that induces a new $\lambda$-dense system of 'well-behaved' submodules of $M$. Call it $\mathcal{C}'=\{N\in \mathcal{C}| Ker(f|R^{(x_N)}) \in \mathcal{K} \cap \mathcal{L} \}$.
\bigskip

\emph{CLAIM.} $\mathcal{C'} \subseteq {^{\perp}C}$.
\smallskip

\emph{PROOF.} Given any $N\in \mathcal{C}'$, $Ext^1(N,C)=0$ iff each morphism $h:Ker(f|R^{(x_N)}) \rightarrow C$ can be extended to some $h':R^{(x_N)}\rightarrow C$, and this holds by our construction.

Indeed, $[R^{(x_N)},C]_R \twoheadrightarrow{} [K,C]_R$ whenever $Ext^1(M,C)=0$ (that holds by assumption), and $[K,C]_R \twoheadrightarrow{} [Ker(f|R^{(x_N)}),C]_R$ provided that $K/Ker(f|R^{(x_N)}) \in {^{\perp}C}$, which is given by our construction of $\mathcal{L}$.
\bigskip

\emph{CLAIM.} $\mathcal{C}' \subseteq  \mathcal{A}={^{\perp}\mathcal{B}}$.
\smallskip

\emph{PROOF.} $\mathcal{C}\subseteq {^{\perp}C}$ by the \emph{Lemma}. Moreover, $N = Ker[M \twoheadrightarrow{} M/N]$, for $M \in {^{\perp}(C^I)^{cc}}$ and $M$ strict $C$-stationary, so that we can apply \emph{Lemma} \ref{stationarity along ses} to obtain that also $N$ is strict $C$-stationary.

Hence $\mathcal{C}$ consists of strict $C$-stationary modules in $^{\perp}C$.

\texttt{Assume} further that \texttt{$\mathcal{C}$ consists of c.p.\ modules}, so as to meet the assumptions of \emph{Lemma} \ref{closure of Ext conditions}. Then, each $N\in \mathcal{C}$ is also in $^{\perp}D$ whith $D$ ranging (up to iso) over all the pure submodules of direct powers of $C$.

However, by our construction of $C$ notice that the latter is the general form of the elements in $\mathcal{B}$, so that $N\in {^{\perp}\mathcal{B}}$. In other words, $\mathcal{C}' \subseteq \mathcal{A}$, as required.
\end{proof}

\begin{lemma}\label{cofinal refinement}
Given a cotorsion pair $(\mathcal{A},\mathcal{B})$, consider a $\lambda$-dense system $\mathcal{C}_\lambda$ with $M/N \in \mathcal{A}$ for each $N\in \mathcal{C}_\lambda$. Then, $\mathcal{C}_\lambda \cap \mathcal{A}$ is again $\lambda$-dense and satisfies the $Ext^1$-vanishing property of quotients, provided that it is cofinal to $\mathcal{C}_\lambda$ (i.e. it is such that for each $N\in \mathcal{C}_\lambda$ there is some $L\in \mathcal{C}_\lambda \cap \mathcal{A}$ containing $N$).
\end{lemma}

\begin{proof}
Clearly the class of all this $L$'s 'transitively' satisfies the local part of $\lambda$-density. The closure under well-ordered ascending $(<\lambda)$-chains goes as follows. Consider any well-ordered ascending chain in $\mathcal{C}_\lambda \cap \mathcal{A}$, say $(N_\alpha |\alpha < \lambda)$.

For any successor ordinal $\alpha +1$, $N_{\alpha +1}/N_{\alpha} \in {^{\perp}\mathcal{B}}$: we required $M/{N_\alpha} \in {^{\perp}\mathcal{B}} = \mathcal{A}$, so that the restriction map induces an epi $[M,\mathcal{B}]_R \twoheadrightarrow{} [N_\alpha^i,\mathcal{B}]_R$, that in turn factorizes into $[M,\mathcal{B}]_R\rightarrow[N^i_{\alpha+1},\mathcal{B}]_R \rightarrow [N_\alpha^i,\mathcal{B}]_R$, thus yielding the required epi $[N^i_{\alpha+1},\mathcal{B}]_R \twoheadrightarrow{} [N_\alpha^i,\mathcal{B}]_R$. In other words, our ascending chain is an $\mathcal{A}$-filtration of its limit (i.e. its union), and $\mathcal{A}$ is closed under $\mathcal{A}$-filtrations by Eklof Lemma.

Finally, still $M/L \in \mathcal{A}$ by our assumption.
\end{proof}

We close this section with a simplified version of Shelah's Singular Compactness Theorem, as stated by Saroch and Stovicek (see \emph{Theorem} \cite{Saroch&Stovicek},3.2).

\begin{theorem}\label{compactness principle}
Let $\kappa$ be a singular cardinal, and consider a $\kappa-g.$ module $M$. Consider some cardinal $\mu$ s.t. $cf(\kappa) \leq \mu < \kappa$, and assume to be given, for each $\mu < \lambda < \kappa$, a $\lambda$-dense system $\mathcal{C}_{\lambda}$ consisting of submodules of $M$.

Then, there is a filtration of $M$, $(M_\alpha | \alpha \leq cf\kappa)$ and a continuous strictly increasing chain of cardinals $(\kappa_\alpha | \alpha < cf\kappa)$ that is cofinal to $\kappa$ s.t. $M_\alpha \in \mathcal{C}_{\kappa_{\alpha}^+}$ for each $\alpha < cf\kappa$.
\end{theorem}

\begin{proof}
Fix any chain $(\kappa_\alpha | \alpha <cf\kappa)$ with $\mu < \kappa_0$.

Consider an increasing chain $(X_\alpha | \alpha < cf\kappa)$ of subsets of $M$ whose union generates $M$, and s.t. $\# X_\alpha = \kappa_\alpha$ for each $\alpha < cf \kappa$.

The localness property of dense systems translates this into an increasing chain of submodules of $M$, namely $(N_\alpha^0 | \alpha < cf\kappa)$ s.t. $N_\alpha ^0 \in \mathcal{C}_{\kappa_{\alpha}^+}$ and $X_\alpha \cup (\bigcup_{\beta < \alpha} N_{\beta}^{0}) \subseteq N_{\alpha}^0$ ($\alpha < cf\kappa$). Furthermore, fix a generating set $Y_\alpha ^0$ for each $N_\alpha ^0$ together with an enumeration $Y_\alpha ^0 := \{ y_{\alpha,\gamma}^0 | \gamma < \kappa_\alpha \}$.

We will refine this initial chain along our dense systems in a gradual manner, that resembles Cantor's diagonal procedure.

More exactly, we will construct countably-many increasing chains $\mathcal{N}_n:=(N_\alpha^n | \alpha < cf\kappa)$ of submodules of $M$. To better control these chains, we will progressively fix the corresponding generating sets $Y_\alpha ^n$, together with an enumeration $Y_\alpha ^n:= \{ y_{\alpha,\gamma}^n | \gamma < \kappa_\alpha \}$; this will be performed requiring our chains to satisfy the following compatibility properties:
\begin{enumerate}
    \item $N_\alpha^n \in \mathcal{C}_{\kappa_\alpha ^+}$
    
    \item $N_\alpha ^n \quad \supseteq \quad \{ y_{\delta,\gamma}^{n-1} \quad  | \quad \alpha \leq \delta < cf\kappa \quad \wedge \quad \gamma < \kappa_\alpha \} \quad \cup \quad  \bigcup_{\beta < \alpha} N_\beta ^n$
\end{enumerate}

One can spatially visualize our chains in a $cf\kappa \times \omega$ rectangular grid with rows indexed by $\alpha$ and columns indexed by $n$.

So, the first property means that we are refining our sequence jumping around the various given dense systems so that the $\alpha$-th row is made of modules in $\mathcal{C}_{\kappa_\alpha^+}$.

As for the second property, it ensures the compatibility of our enumerations; in other words, it entails that each column is increasing, and that the limit along the diagonal of the grid contains all the generators $y_{\delta,\gamma}^{n}$'s.

Moreover, the weak inequality $\alpha \leq \delta$ in the second property implies that also rows are increasing, since each $N_\alpha ^n \supseteq N_\alpha ^{n-1}$, that is generated by $\{ y_{\alpha,\gamma}^{n-1} | \gamma < \kappa_\alpha \}$.

Let's now construct our chains via a double induction. For each $n<\omega$, perform induction on $\alpha$ as follows.

Assume to have already constructed $(\mathcal{N}_l | l\leq n)$ as well as $(N_\beta ^{n+1} | \beta < \alpha)$.
$$X^{n+1}_\alpha := \{ y_{\delta,\gamma}^{n} \quad  | \quad \alpha \leq \delta < cf\kappa \quad \wedge \quad \gamma < \kappa_\alpha \} \quad \cup \quad (\cup_{\beta < \alpha} Y^{n+1}_\beta)$$
has cardinality $\kappa_\alpha$, so let $N_\alpha ^{n+1} \in \mathcal{C}_{\kappa_{\alpha}^+}$ contain $X^{n+1}_\alpha$; fix one of its generating sets $Y_{\alpha}^{n+1}$ containing $X^{n+1}_\alpha$ and having again cardinality $\kappa_\alpha^+$.

Consider the limit along the rows, ($\alpha < cf\kappa$) $M_\alpha := \bigcup_{n<\omega} N_\alpha^n \in \mathcal{C}_{\kappa_\alpha^+}$ (since the latter system is closed under unions of well-ordered ascending chains). Remark that $\cup_{\alpha < cf \kappa} M_\alpha = M$, because it holds $X_\alpha \subseteq N_\alpha ^0 \subseteq M_\alpha$, and $\cup X_\alpha$ generates $M$.

In our picture, spatially place the limits in the extra column indexed by $\omega$.

All the chains considered so far need not be continuous; however, we claim that the limit chain along the $\omega$-th column is such.

To show this, fix any limit ordinal $\alpha < cf\kappa$; it plainly holds $M_\alpha \supseteq \bigcup_{\beta <\alpha} M_\beta$ (we have term-wise inclusion for each $n$).

For the opposite inclusion we will exploit the Cantor diagonal construction. Indeed, fixed $n < \omega$, for each $\beta < \alpha$ one has $\{ y^{n-1}_{\alpha,\gamma} | \gamma < \kappa_\beta \} \subseteq N_{\beta}^n$ by the second property. Hence, $Y_\alpha ^{n+1} \subseteq \bigcup_{\beta<\alpha} N_{\beta}^{n}$, that implies $N_{\alpha}^{n-1} \subseteq \bigcup_{\beta<\alpha} N_\beta^n$.

In other words, each point in the grid is included in the union of the part of the successive column right above it. Thus, by moving along the increasing direction of $n$ (w.l.o.g.\ downwards) from $N_\alpha^n$, we remain within the union of the triangle above the anti-diagonal with one vertex 'almost' at $N_\alpha^n$.

More precisely, this amounts to the required condition $M_\alpha \subseteq \bigcup_{\beta<\alpha} M_\beta$.

Finally, modify $M_0:=0$, and set $M_{cf\kappa}:=M$, so as to satisfy the axiom of a filtration.
\end{proof}

\section{The main Theorem}

\subsection{The Countable Telescope Conjecture}

We will apply the preliminary digression on refinements of dense systems to the well-known Countable Telescope Conjecture.

\begin{theorem}[CTCM]
Let $\mathcal{C}=(\mathcal{A},\mathcal{B})$ be a cotorsion pair with $\lim \mathcal{B} = \mathcal{B}$. Then, $\mathcal{C}$ is of countable type, and $\mathcal{B}$ is definable.
\end{theorem}

\begin{proof}
Let $C$ be the 'canonical' pure-injective cogenerator constructed for $\mathcal{B}$ as in \emph{Lemma} \ref{canonical cogenerator}. Indeed, $\mathcal{B}$ is clearly closed under direct products, being the right class of a cotorsion pair.

Call $\mathcal{A}_0 := \mathcal{A}^{\leq \omega}$ the class of all the c.p.\ modules in $\mathcal{A}$. We will proceed by induction on cardinals $\kappa$ to prove that $\kappa-p.$ modules in $\mathcal{A}$ are in $Filt(\mathcal{A}_0)$. (w.l.o.g.\ $\kappa$ is uncountable.)

Observe that we can reduce the induction procedure to regular uncountable cardinals, because they constitute a proper class.
\bigskip

\emph{IDEA.} Consider a $\kappa-p.$ module $M\in \mathcal{A}$. For each \texttt{uncountable regular cardinal} $\lambda \leq \kappa$, we will construct by transfinite induction a $\lambda$-dense system $\mathcal{C}_\lambda \subseteq \mathcal{A}$ of submodules of $M$ s.t. $M/N \in \mathcal{A}$ for each $N\in \mathcal{C}_\lambda$. We will use them to build the required $\mathcal{A}_0$-filtration of $M$.

Indeed, assume to have such a family $\{\mathcal{C}_\lambda | \lambda<\kappa \}$. We will build an $\mathcal{A}$-filtration of $M$ consisting of $(<\kappa)-p.$ modules. Then, notice that, by the induction premise, we will be able to assume that each term of the $\mathcal{A}$-filtration is further $\mathcal{A}_0$-filtered, so that also $M$ will be such.
\bigskip

Now, let's implement the \texttt{induction procedure}. The initial step is plainly provided by \emph{Lemma} \ref{initial step} (whose assumptions are satisfied by our construction together with the completeness of $\mathcal{C}$).
\smallskip

\begin{itemize}
    \item \textsc{Case 1: $\lambda = \aleph_1$.}
        Let $\mathcal{C}_{\aleph_1}$ be an $\aleph_1$-dense system in $M$ as provided by \emph{Lemma} \ref{initial step}, and consider a presentation SES of $M$: $0\rightarrow K \rightarrow R^{(\kappa)} \rightarrow M \rightarrow 0$. By \emph{Lemma} \ref{refinement of dense systems} we can refine our system compatibly with the presentation, so as to obtain w.l.o.g.\ $\mathcal{C}_{\aleph_1} \subseteq \mathcal{A}$, as required.
        
    \item \textsc{Case 2: $\lambda$ regular limit cardinal.}
        Let $\mathcal{C}_\lambda$ be a $\lambda$-dense system as provided by \emph{Lemma} \ref{initial step}. We will construct a cofinal subsystem into $\mathcal{C}_\lambda \cap \mathcal{A}$ in the spirit of \emph{Lemma} \ref{cofinal refinement}.
    
        Pick any $N_0 \in \mathcal{C}_\lambda$, and let $\mu_0$ be a regular cardinal s.t. $N_0$ is $(<\mu_0)-p.$ Consider the dense system $\mathcal{C}_{\mu_0} \subseteq \mathcal{A}$ given by the induction premise, and chose some $N_0 ' \in \mathcal{C}_{\mu_0}$ containing $N_0$. Now, there is some $N_1 \in \mathcal{C}_\lambda$ containing $N_0 '$; since the latter is still $(<\lambda)-p.$ we can choose a new $\mu_1$ and iterate.
    
        In this way we construct a countable chain $N_0 \leq N_0' \leq N_1 \leq N_1' \leq \dots$ whose directed union is $N$, and s.t. it possesses two cofinal subsystems, namely $(N_n |n<\omega)$ and $(N_n' | n<\omega)$. The first one yields $N\in \mathcal{C}_\lambda$, while the second one entails $N\in \mathcal{A}$ (it provides an $\mathcal{A}$-filtration of $N$ that is in $\mathcal{A}$ by Eklof Lemma). Call again $\mathcal{C}_\lambda$ the obtained cofinal subfamily of $\mathcal{C}_\lambda \cap \mathcal{A}$.
    
    \item \textsc{Case 3: $\lambda$ successor cardinal.} (Recall that successor cardinals are always regular.)

        \begin{itemize}
            \item \textsc{Case 3.1: $\lambda$ successor of a regular cardinal.}
                Let $\nu$ be a regular cardinal s.t. $\lambda = \nu +1$. Construct the usual $\lambda$-dense system $\mathcal{C}_\lambda$, and pick up any $N_0\in \mathcal{C}_\lambda$. w.l.o.g.\ $N_0$ is $(<\nu)-p.$, and consider the corresponding $\nu$-dense system $\mathcal{C}_\nu \subseteq \mathcal{A}$ of $M$ given by the inductive premise.
                
                We can inductively build a chain $\mathcal{N}_0 := \big( N^0_{\alpha} | \alpha <\nu_0 \big) \subseteq \mathcal{C}_\nu \subseteq \mathcal{A}$ whose union contains $N_0$.
                
                Indeed, let $X^0 := \{ x^0_\alpha | \alpha <\nu \}$ be a generating set for $N_0$, and take an initial subset $X_0$; by the localness part of density, there is some module $N^0_0 \in \mathcal{C}_\nu$ containing $X_0$.  Now, at each successor step $\alpha +1$, call $Y_\alpha$ a generating set of $N^0_\alpha$, and consider $X_{\alpha+1} := Y_\alpha \cup \{x_{\alpha +1}\}$; take a corresponding $X_{\alpha+1} \subseteq N^0_{\alpha +1} \in \mathcal{C}_{\nu}$. At each limit step make the chain continuous by taking the union of the preceding modules.
                
                $\mathcal{C}_\nu$ is closed under union of well-ordered ascending chains, so $\cup \mathcal{N}_0 \in \mathcal{C}_\nu \subseteq \mathcal{A}$.
                
                Choose some $N_1 \in \mathcal{C}_\lambda$ containing such a union, and - from the choice of a generating set $X^1$ of cardinality $\nu$ - build a similar chain $\mathcal{N}_1 = \big( N^1_\alpha | \alpha < \nu_1 \big)$ containing $N_1$; iterate countably many times, and set $N := \cup_{i<\omega} N_i$.
                
                Remark that $\nu_i = \nu$ for each $i<\omega$, because $\nu$ is a regular cardinal and otherwise we would write $X^i = \cup_{\alpha < \nu_i} X^i_\alpha$ as a union of less than $\nu$ subsets of strictly smaller cardinality. Indeed, the $X^i_\alpha$'s should be thought as a set of $\nu$ pointer-variables for the possibly repeated generators $x^i_\alpha$'s; this enables us to carry on the regularity argument, up to allowing for possible repetitions in the chains $\mathcal{N}_i$'s.
                
                Thus, we could also assume $N^{i}_{\alpha} \subseteq N^{i+1}_{\alpha}$ at each step by requiring $X^{i+1}_{\alpha +1} = Y^{i+1}_{\alpha} \cup X^i_{\alpha+1} \cup \{x^{i+1}_{\alpha +1}\}$.
                
                One could visualize the various chains as a rectangular $\omega \times \nu$ grid with rows indexed by $i$ and columns indexed by $\alpha$.
                
                Add an extra row in position $\omega$ given by the ascending chain $\big( \cup_{i<\omega} N^i_{\alpha} | \alpha < \nu \big)$ consisting of the unions along the columns. Since each term is in $\mathcal{C}_\lambda \subseteq \mathcal{A}$, it forms an $\mathcal{A}$-filtration of $N$ as in \emph{Lemma} \ref{cofinal refinement}, and so $N \in \mathcal{A}$ itself.
                
                The $N$'s constitute, then, the sought cofinal dense subsystem of $\mathcal{C}_\lambda \cap \mathcal{A}$; call it again $\mathcal{C}_\lambda$.
                
            \item \textsc{Case 3.2: $\lambda$ successor of a singular cardinal.}
            Let $\lambda = \nu^+$, for $\nu$ singular cardinal of cofinality $cf(\nu) = \mu<\nu$. Choose a strictly increasing continuous chain $(\nu_\alpha | \alpha <\mu)$ of infinite cardinals which is cofinal to $\nu$ and that starts with $\nu_0 > \mu$.
            
            Consider the usual $\lambda$-dense system $\mathcal{C}_{\lambda}$ given by our \emph{Lemma}, and pick any $N_0$ from there.
            
            We will perform a construction similar to the above case, in the sense that we will again move around a rectangular grid, although this time of dimension $\omega \times \mu^+$. However, since $\nu$ is no longer regular and we can apply our \emph{Lemma} only with regular cardinals, to 'access' $\nu$ we need a Cantor diagonal argument going through the various regular cardinals $\nu_\alpha$'s.
            
            Thus, we cannot construct a chain at a time, and we will employ \emph{Theorem} \ref{compactness principle} to obtain the required $\mathcal{N}_n := (N_{\alpha}^n | \alpha < \mu) \subseteq \mathcal{C}_{\nu_\alpha}$.
            
            Indeed, w.l.o.g.\ $N_0$ is $\nu-g.$ with generating set $X$, and consider an ascending chain $(X_\alpha | \alpha < \mu)$ of subsets of $X$ s.t. $\# X_\alpha = \nu_\alpha$.
            
            Inductively construct the chain $\mathcal{N}_0$ as in the \emph{Theorem}, and carry on the inductive enlargement.
            
            This will provide for the countably-many $\mu$-indexed increasing chains $\mathcal{N}_i = (N_\alpha ^i | \alpha < \mu) \subseteq \mathcal{C}_{\nu_\alpha^+}$.
            
            However, this time perform the diagonal procedure on a grid with rows indexed by $\alpha \leq \mu$ (and \emph{not} $\alpha < \mu$) and columns by $i < \omega$, where we let each $N^{i}_\mu$ to be some module of $\mathcal{C}_{\lambda}$ which contains the union $\cup \mathcal{N}_i$ along the already constructed $i$-th column. We choose $N^{i}_\mu$ at the beginning of the $(i+1)$-th step, together with a compatible generating set $Y^i_\mu := \{ y_{\gamma}^i | \gamma < \lambda \}$.
            
            So, for example, (in the notation of the proof of \emph{Theorem} \ref{compactness principle} adapted to the current one) the set $X^{i+1}_{\alpha}$ becomes
            $$X^{i+1}_\alpha := \{ y_{\delta,\gamma}^{i} \quad  | \quad \alpha \leq \delta \leq \mu \quad \wedge \quad \gamma < \kappa_\alpha \} \quad \cup \quad (\cup_{\beta < \alpha} Y^{i+1}_\beta)$$
            
            And the rest is analogous.
            
            Define $N:=\cup_{i<\omega} N_i$ as the union of the $\mu$-th row, and observe that it is in $\mathcal{C}_\lambda$.
            
            Moreover, it equals the union performed first along the rows and then along the columns, namely $N= \cup \mathcal{H}$ for $\mathcal{H}:=(\cup_{i<\omega} N_{\alpha}^i | \alpha \leq \mu)=(\cup_{i<\omega} N_{\alpha}^i | \alpha < \mu)$.
            
            Where the latter equality holds because, by our modification of the sets $X^{i+1}_\alpha$'s, it is irrelevant to include or not $\mu$ in the union (it is actually the diagonal procedure).
            
            The latter increasing chain is continuous, and has each term in $\mathcal{C}_{\nu_\alpha^+} \subseteq \mathcal{A}$ (by the general induction premise). Thus, as in \textsc{Case 3.1}, it provides for an $\mathcal{A}$-filtration of $N \in \mathcal{C}_\lambda$, which turns out to live in $\mathcal{C}_\lambda \cap \mathcal{A}$.
        \end{itemize}
\end{itemize}

Thus, having constructed the families, let's provide the sought \texttt{$\mathcal{A}$-filtration of $M$} consisting of $(<\kappa)-p.$ modules.

If $\kappa$ is regular, then $\mathcal{C}_\kappa$ already contains the filtration for $M$.

Otherwise, in our setting, \emph{Theorem} \ref{compactness principle} applied to $M$ and $(\mathcal{C}_\lambda | \lambda <\kappa ) \subseteq \mathcal{P}(\mathcal{A})$ yields a $cf\kappa$-filtration of $M$ of the form $(M_\alpha | \alpha < cf\kappa)$ s.t. $M_\alpha \in \mathcal{C}_{\kappa_\alpha^+}\cap \mathcal{A}$, and it turns out to be an $\mathcal{A}$-filtration (because the $Ext^1$-vanishing property of the quotients is the same for each dense system, so we can follow again the proof of \emph{Lemma} \ref{cofinal refinement}).

Hence, since by our induction premise each $M_\alpha$ is $\mathcal{A}_0$-filtered, also $M$ is such.

Indeed, given the $\mathcal{A}_0$-filtrations $(M_{\alpha,\beta} | \beta < \eta_\alpha)$ relative to $M_\alpha$, modify them by setting $M_{\alpha+1,\beta}' := M_{\alpha+1,\beta}+M_\alpha$ and merge them together into a $(\sum_{\alpha < cf\kappa} \eta_\alpha)$-filtration of $M$. (Up to the elimination of some terms, we can further assume that the chain is strictly increasing.)

Moreover, it is clearly an $\mathcal{A}_0$-filtration, since the successive quotients are not altered by our modification, namely $M_{\alpha+1,\beta+1}'/M_{\alpha+1,\beta}' \simeq M_{\alpha+1,\beta+1}/M_{\alpha+1,\beta}$.

Thus, $(\mathcal{A},\mathcal{B})$ is $\mathcal{A}_0$-deconstructible. To prove its countable type, we are left to show that modules in $\mathcal{A}_0$ have \texttt{c.p.\ }\texttt{first syzygies}.

Take $M_0 \in \mathcal{A}_0$ and consider a resolution $0\rightarrow \Omega^1(M_0) \rightarrow P_0 \rightarrow M_0 \rightarrow 0$.

$M_0\in \mathcal{A}_0 \subseteq {^{\perp}\mathcal{B}}$ and $P_0$ is projective (and hence ML and strict $C$-stationary, by \emph{Lemma} \ref{ML vs stationarity}), so that the assumptions of \emph{Lemma} \ref{stationarity along ses} are satisfied, and also $\Omega^1(M)$ is strict $C$-stationary. Thus, it is strict $def(C)$-stationary (see \emph{Lemma} \ref{closure properties}) and in particular strict $\mathcal{B}$-stationary (by the construction of $C$).

In particular, $\Omega^1(M_0)$ is strict $_{R}R^c$-stationary ($R^c \in \mathcal{B}$, being injective), and hence $_{R}R$-ML (see \emph{Lemma} \ref{ML vs stationarity}/1).

Now, being $M_0$ c.p., $\Omega^1(M_0)$ is c.g.\ by a classical pullback argument, and hence also c.p, because the two notions agree for $_{R}R$-ML modules (see \emph{Cor.} \ref{c.g. iff c.p.}).

Finally, let's turn to the \texttt{definability of the class $\mathcal{B}$}.

Under our closure assumptions, it equates the closure of $\mathcal{B}$ under pure submodules. We will show this via a slight modification of the proof of \emph{Theorem} \cite{Trlifaj&Goebel},13.41.
\bigskip

\emph{CLAIM.} Let $\mathcal{B}$ be a class of modules closed under $\prod$, $\lim$, and consider any c.p.\ module $A\in {^{\perp}\mathcal{B}}$. Then, $A^{\perp}$ is closed under copies of pure submodules of products in $\mathcal{B}$.

\smallskip

\emph{PROOF.} By \emph{Lemma} \ref{uniform stationarity}, $A\in {^{\perp}\mathcal{B}}$ is strict $\mathcal{B}$-stationary. Choose a countable tower presentation of $A$ consisting of f.p.\ modules that witnesses such a stationarity. Notice that it further has the uniform factorisation property, so that, as remarked in \emph{Lemma} \ref{closure properties}, it is moreover $def(\mathcal{B})$-stationary.

Consider $B' \subseteq_{\ast} \prod B_\alpha$ as in the statement. In particular, our tower is $B'$-stationary. In other words, the induced inverse system $([A_i,B']_R,h^{\ast}|i<\omega)$ satisfies the ML condition.

Now, consider the pure-SES $0 \rightarrow B' \rightarrow_{\ast} \prod B_\alpha \overset{\rho}{\rightarrow} \prod B_\alpha / B' \rightarrow 0$

Since each $A_i$ is f.p.\, by the properties of purity $([A_i,\_]_R | i<\omega)$ yields an inverse tower of SES's, namely
$$0 \rightarrow [A_i,B'] \rightarrow [A_i,\prod B_\alpha] \overset{(\rho_i)_\ast}{\rightarrow} [A_i,\prod B_\alpha / B'] \rightarrow 0$$
with transition morphisms induced by the tower $(h_i|i<\omega)$.

Since the system in the first member satisfies the Mittag-Leffler property, passing to the inverse limit we obtain the SES:
$$0 \rightarrow \varprojlim_{i<\omega} [A_i,B']_R \rightarrow \varprojlim_{i<\omega} [A_i,\prod B_\alpha]_R \overset{\rho_\ast}{\rightarrow} \varprojlim_{i<\omega} [A_i,\prod B_\alpha / B']_R \rightarrow 0$$

And, bringing the $\varprojlim$ inside the first component of the Hom, the latter is properly
$$0 \rightarrow [A,B'] \rightarrow [A,\prod B_\alpha] \overset{\rho_\ast}{\rightarrow} [A,\prod B_\alpha / B'] \rightarrow 0$$
so that, in particular, $\rho_\ast = [A,\rho]$ is epi.

Since also $Ext^1(A,\prod B_\alpha)=0$ by the properties of $Ext^1$ and our assumption, we conclude that $Ext^1(A,B')=0$, as required.

\end{proof}

\subsection{Applications}

In this section we will sketch the main applications of the CTCM as presented in \cite{Saroch-applications}, the twin paper \emph{'Approximations and Mittag-Leffler Conditions. The Applications'} by Angeleri-H\"ugel, Saroch and Trlifaj. Remarkably, we will be able to state and prove a weaker version of the Enochs' Conjecture for cotorsion pairs with right class closed under direct limits.

In what follows, we will just hint at the various underlying ideas, briefly highlighting how to apply the previously introduced tools, but without claiming to be exhaustive.
\smallskip

We begin our detour by recalling the notion of a Bass module over a c.p.\ class $\mathcal{D}$, which is a module that can be written as the direct limit of a tower in $\mathcal{D}$.

The terminology comes after the famous Bass 'P' Theorem, in which such a class was introduced to produce a counterexample of a c.p. flat $R$-module that is not projective over a non-right perfect ring $R$.

Subsequently, in \cite{tree modules} Slavik and Trlifaj proved, by introducing the construction of tree modules, that Bass modules play a crucial role in deconstructibility issues.

In the rest of this section, we will present various deep results, that are essentially a consequence of refining procedures of suitable dense systems. We will then be able to employ Bass modules to characterize the cotorsion pairs with both classes closed under direct limits among those with a limit-closed right class.

This will turn out to yield a weaker version of the Enochs' Conjecture, as well as a deconstructibility criterion for the left classes at stake. Finally, we will be able to describe the kernels of such cotorsion pairs in a manner that is formally analogous to the particularly interesting case of tilting cotorsion pairs.
\smallskip

In what follows, we will almost always work with a cotorsion pair $\mathcal{C}=(\mathcal{A},\mathcal{B})$ s.t. $\mathcal{B} = \lim \mathcal{B}$ is closed under direct limits.

Denote by $\bar{\mathcal{A}}$ the class of all pure-epi images of $\mathcal{A}$. 

We will need the following powerful technical result. 

\begin{lemma}
Let $(\mathcal{A},\mathcal{B})$ be a cotorsion pair s.t. $\mathcal{B}=\lim \mathcal{B}$. Let $\mathcal{B}'$ be the class of all pure-injective modules in $\mathcal{B}$. Then, $^{\perp}\mathcal{B}' = \bar{\mathcal{A}}$
\end{lemma}

Moreover, one may prove the following \emph{Lemmas}.

\begin{lemma}{(\cite{Saroch-applications},3.3)}
    Let $(\mathcal{A},\mathcal{B})$ be a cotorsion pair s.t. $\mathcal{B}=\lim \mathcal{B}$. Then, $\bar{\mathcal{A}} = \lim \bar{\mathcal{A}}^{\leq \omega}$.
\end{lemma}

\begin{proof}
    \emph{(Sketch.)} The proof provides a refinement of a suitable $\aleph_1$-dense system, as in \emph{Lemma} \ref{refinement of dense systems}.
\end{proof}
    
\begin{lemma}{(\cite{Saroch-applications},3.4)}
Let $(\mathcal{A},\mathcal{B})$ be a cotorsion pair. Under the relaxed assumption on countable limits $\mathcal{B} = \lim_\omega \mathcal{B}$, then $\bar{\mathcal{A}}^{\leq \omega} \subseteq \lim \bar{\mathcal{A}}^{< \omega}$ whenever $\mathcal{C}$ is further complete.
\end{lemma}    

\begin{proof}
    \emph{(Sketch.)} The proof is reminiscent of El Bashir's generalization of Enoch's Theorem to Grothendieck categories: expand a f.p.\ countable presentation $\mathcal{M}$ of any $M\in \bar{\mathcal{A}}^{\leq \omega}$ to a direct system of special $\mathcal{A}$-precovers $\pi_n:A_n \rightarrow M_n$. Then, the limit of the precovers is a pure-epi; so, by the finiteness property of $\mathcal{M}$, up to dropping some terms, each transition map $f_{n+1,n}$ factors through the successive precover $\pi_{n+1}$ via some $\nu_{n}\in [M_n,A_{n+1}]$.
    
    Thus, $(M,f_n \circ \pi_n|n<\omega)=\lim(A_n,g_{nm}|m\leq n<\omega)$ for $g_{nm}=\nu_{n-1} \circ f_{n-1,m} \circ \pi_{m}$, so that $M\in \lim_\omega \mathcal{A}$.
\end{proof}

\begin{corollary}
In particular, provided that actually  $\mathcal{B} = \lim \mathcal{B}$, the CTMC implies that $\mathcal{C}$ is of countable type, and we obtain that $\bar{\mathcal{A}}^{\leq \omega}$ consists of Bass modules over $\mathcal{A}^{\leq \omega}$. Moreover, if further $\mathcal{C}$ is of finite type, then it suffices to consider Bass modules over $\mathcal{A}^{< \omega}$.
\end{corollary}

\begin{proof}
\emph{(Sketch.)} Suitably refine the previous system of special $\mathcal{A}$-precovers via Hill Lemma, so as to obtain the required Bass module; notice that any countable tower is in particular $\aleph_1$-continuous.
\end{proof}

In other words, we can sum up the outcome entailed by these deep results as the following Theorem.

\begin{theorem}{(Characterization of $\bar{\mathcal{A}}$.)}
Any pure-epi image of a module from $\mathcal{A}$ is the direct limit of a direct system of Bass modules from $\mathcal{A}^{\leq \omega}$.
\end{theorem}

As a Corollary, for a cotorsion pair $(\mathcal{A},\mathcal{B})$ with $\mathcal{B}=\lim \mathcal{B}$ and a class $\mathcal{G}$, the $\mathcal{G}$-stationarity of pure-epi images of modules from $\mathcal{A}$ can be checked on Bass modules over $\mathcal{A}^{\leq \omega}$. More precisely,

\begin{corollary}{(\cite{Saroch-applications},3.5)}
Let $(\mathcal{A},\mathcal{B})$ be a cotorsion pair with $\mathcal{B}=\lim \mathcal{B}$. Consider any class of modules $\mathcal{G}$.

If all Bass modules over $\mathcal{A}^{\leq \omega}$ are $\mathcal{G}$-stationary, then also $\bar{\mathcal{A}}$ consists of $\mathcal{G}$-stationary modules.
\end{corollary}

\begin{proof}
\emph{(Sketch.)} By the previous \emph{Theorem}, we can write the modules from $\bar{\mathcal{A}}$ as direct limits of direct systems of Bass modules from $\mathcal{A}^{\leq \omega}$. By \ref{ML vs stationarity}, all the considered Bass modules turn out to be $\mathcal{G}^c$-ML. Since ML conditions can be checked on countable chains, $\bar{\mathcal{A}}^{\leq \omega}$ consists of $\mathcal{G}^c$-ML iff $\mathcal{G}$-stationary.
\end{proof}

Hence, we can deeply characterize the cotorsion pairs with both classes closed under direct limits among those that satisfy the closure requirement on the right class.

\begin{theorem}{(\cite{Saroch-applications},3.6)}\label{Bass modules and closure}
Let $\mathcal{C}=(\mathcal{A},\mathcal{B})$ be a cotorsion pair s.t. $\mathcal{B}=\lim \mathcal{B}$. TFAE:
\begin{enumerate}
    \item $\mathcal{C}$ is cogenerated by a (discrete) pure-injective module.
    
    \item $\mathcal{A}$ is closed under pure-epi, pure-submodules.
    
    \item $\mathcal{C}$ is closed (i.e. $\mathcal{A} = \lim \mathcal{A}$).
    
    \item $\mathcal{A}$ contains all Bass modules over $\mathcal{A}^{\leq\omega}$.
\end{enumerate}
\end{theorem}

Its proof is long and heavily relies on the ideas about refinements of dense systems involved in the proof of the CTCM and on the following technical deconstructibility result, namely \emph{Proposition} \cite{Saroch&Stovicek},5.12.

\begin{lemma}
Let $(\mathcal{A},\mathcal{B})$ be a complete cotorsion pair with $\mathcal{B}$ closed under direct limits. Then, there is some pure-injective module $E$ s.t. $^{\perp}E = \bar{\mathcal{A}}$. Moreover, we can assume $E$ to be discrete, i.e. to be the product of some indecomposable pure-injectives.
\end{lemma}

Thence, we have the following corollary to the previous theorem, \cite{Saroch-applications},3.7.

\begin{corollary}
Let $\mathcal{C}=(\mathcal{A},\mathcal{B})$ be a cotorsion pair of finite type. Then, $\mathcal{C}$ is closed iff $\mathcal{A}$ contains all Bass modules over $\mathcal{A}^{<\omega}$.
\end{corollary}

Now, let's turn our attention to cotorsion pairs of countable type. Let $\mathcal{L}$ be the class of all locally $\mathcal{A}^{\leq \omega}$-free modules.

As an application of the tree modules construction, in \cite{tree modules} Slavik and Trlifaj proved that the class $\mathcal{L}$ is non-deconstructible provided that there is some Bass module over $\mathcal{A}^{\leq \omega}$ that does not belong to $\mathcal{A}$.

Thus, by applying this result to the previous theorem, we have the following characterization of the deconstructibility character of locally free $\mathcal{A}^{\leq \omega}$-free modules.

\begin{theorem}{(Deconstructibility, \cite{Saroch-applications},4.1)}
Let $(\mathcal{A},\mathcal{B})$ be a cotorsion pair s.t. $\mathcal{B}$ is complete, and let $\mathcal{L}$ be as before. Then, $\mathcal{L}$ is deconstructible iff $\mathcal{A}$ contains all Bass modules over $\mathcal{A}^{\leq \omega}$ iff $\mathcal{C}$ is closed.
\end{theorem}

\emph{Remark.} The just stated result has interesting applications to tilting classes; in particular, it states that for a tilting module $T$, the class of all locally $T$-free modules is (pre)covering iff $T$ is $\Sigma$-pure-split. (See \cite{Saroch-applications},4.2-4.3.)
\smallskip

Moreover, another application of the classical refining procedure as in \ref{refinement of dense systems}, yields the following characterization of locally $\mathcal{A}^{\leq \omega}$-free modules, \cite{Saroch-applications},4.4.

\begin{theorem}
Let $(\mathcal{A},\mathcal{B})$ be a cotorsion pair s.t. $\mathcal{B}$ is complete. Then $M$ is locally $\mathcal{A}^{\leq \omega}$-free iff $M$ is a  $\mathcal{B}$-stationary pure-epi image of a module from $\mathcal{A}$.

In particular, $\mathcal{L}$ is closed under pure submodules.
\end{theorem}

\begin{proof}
\emph{(Sketch.)} $(\implies)$ $M \in \mathcal{L}$ possesses an $\aleph_1$-dense system of submodules in $\mathcal{A}^{\leq \omega}$ whose union is $M$ itself. Moreover, by the criterion on the Uniform Factorisation Property, $\mathcal{A}$ consists of strict $\mathcal{B}$-stationary modules, so in particular a directed union of modules from $\mathcal{A}$ is such.

Conversely, in our setting we can apply the classical refining procedure to obtain an $\aleph_1$-dense system of submodules of $M$ in $\mathcal{A}^{\leq \omega}$. The last claim relies on the initial technical lemma: $\bar{\mathcal{A}} = {^{\perp}C}$ for an elementary cogenerator for $\mathcal{B}$.
Then apply \emph{Lemma} \ref{ML vs stationarity} and \emph{Lemma} \ref{closure properties}.
\end{proof}

We are now ready to state the most important consequence of the CTCM, namely a weak version of Enochs' Conjecture.

The latter result lies at the heart of representation theory of modules by means of left and right approximations, and can be considered to be by far the most important open one.

Indeed, given a complete cotorsion pair $\mathcal{C}=(\mathcal{A},\mathcal{B})$ (i.e. special preenveloping iff special precovering; with the equivalence due to Salce's Lemma), Enochs proved in a famous Theorem (e.g.\ see \cite{Trlifaj&Goebel} or \cite{Xu}) that $\mathcal{C}$ is perfect (i.e. covering iff enveloping) whenever it is closed.

Now, there is a rich abundance of complete cotorsion pairs, since a variant of Quillen's small object argument states that a cotorsion pair is complete whenever it is generated by a set; hence it is not a really restrictive condition among those of interest. So, the nutshell of the question seems to lie in the implication \emph{closed} $\implies$ \emph{perfect}.

Rewriting the issue from the categorical viewpoint introduced in the \emph{Introduction}, in usual contexts we can always consider the non-empty categories of all the special preenvelopes/precovers of a given module. Then, the closure property of the left class allows us to inductively construct a minimal object for the subobject preorder.

Moreover, notice that, whenever such categories are trivial, our research question is even meaningless.

Accordingly with this logical flow, Enochs' Conjecture inquires whether such a theorem admits a converse, and it posits an affirmative answer.

As already mentioned, such a problem is still open; however, the tools developed to show the CTCM allowed Angeleri H\"ugel, Saroch and Trlifaj to prove a particular case, under the further assumption that $\mathcal{B} = \lim \mathcal{B}$. Whence the name 'Weak Enochs Conjecture'.

Before stating the Theorem, let us incidentally notice that, whenever it further holds $\mathcal{B}=\lim \mathcal{B}$, one has that the left class $\mathcal{A}$ is cogenerated by all the pure-injective modules of $\mathcal{B}$, and the latter cogenerating class can be taken to be a set by an application of Ziegler's Spectrum theorem. Hence, the completeness requirement is automatically satisfied and, as already noticed, due to the properties of the cotorsion pair at stake, the converse implication is a meaningful inquire and does exquisitely involve closure and perfectness. 

\begin{theorem}{(Weak Enochs Conjecture; see \emph{Theorem} \cite{Saroch-applications},5.2)}
Let $\mathcal{C}= (\mathcal{A},\mathcal{B})$ be a cotorsion pair s.t. $\mathcal{B} = \lim \mathcal{B}$. Let $\mathcal{L}$ denote the class of all locally $\mathcal{A}^{\leq \omega}$-free modules. TFAE:
\begin{enumerate}
    \item $\mathcal{C}$ is closed.
    
    \item every module (in $\mathcal{B}$) has an $\mathcal{A}$-cover, i.e. $\mathcal{C}$ is perfect.
    
    \item $Ker(\mathcal{C})$ is closed under $\lim_\omega$.
    
    \item $(\lim \mathcal{A})^{\leq \omega}$ consists of strict $\mathcal{B}$-stationary modules.
    
    \item $Ker(\mathcal{C})$ consists of $\Sigma$-pure-split modules.
    
    \item $\mathcal{L}$ coincides with the class $\bar{\mathcal{A}}$ of all pure-epi images from $\mathcal{A}$.
    
    \item Every module (in $(\lim \mathcal{A})^{\leq \omega}$) has an $\mathcal{L}$-(pre)cover.
\end{enumerate}
\end{theorem}

\emph{Remark.} By \emph{Theorem} \ref{Bass modules and closure}, we could formulate an extra equivalent condition involving Bass modules over $\mathcal{A^{\leq \omega}}$. The fact that these modules enter the picture when considering precovering issues is a deeply-rooted feature.

To grasp some intuition, consider the following instance. In \cite{Saroch-tools}.3.3 Saroch proved that the nicely-behaved class $\mathcal{F}\mathcal{L}$ is precovering iff the underlying ring is right perfect (so surprisingly it is not precovering in general). The latter is due to the role played by Bass modules, which induce some deconstructibility issues for the class of Flat Mittag-Leffler modules, as shown by Slavik and Trlifaj in \cite{tree modules} in the more general context of classes of locally free modules.

Here, the 'locally free' behaviour is given by our characterization of the class $\bar{\mathcal{A}}$, that relies on accurate refinements of dense systems of modules (our basic free pieces).
\bigskip

Moreover, \emph{Lemma} \cite{Saroch-applications},5.4 sheds some more light on the kernel of a cotorsion pair with complete right class.

\begin{lemma}
Let $\mathcal{C}=(\mathcal{A},\mathcal{B})$ be a cotorsion pair with $\mathcal{B}=\lim \mathcal{B}$. Then, there is some module $K$ s.t. $Ker(\mathcal{C}):= \mathcal{A} \cap \mathcal{B} = Add(K)$.
\end{lemma}

The module $K$ is determined by an inductive procedure, and, under further requirements on the kernel, one obtains additional equivalent conditions for the Weak Enochs Conjecture. Namely,

\begin{theorem}{(\cite{Saroch-applications},5.5.)}
Let $\mathcal{C}=(\mathcal{A},\mathcal{B})$ be a cotorsion pair with $\mathcal{B}=\lim \mathcal{B}$. Let $K$ be a module s.t. $Ker(\mathcal{C})=Add(K)$. TFAE:

\begin{enumerate}
    \item $\mathcal{A}$ is closed.
    
    \item $\mathcal{A}$ is covering (is covering for $\mathcal{B}$).
    
    \item Every module in $Ker(\mathcal{C})$ has a semiregular endomorphism ring.
    
    \item $K$ has perfect decomposition.
    
    \item $K$ is $\Sigma$-pure-split.
    
    \item $Ker(\mathcal{C})$ is a covering class (is covering for $\mathcal{B}$).
\end{enumerate}
\end{theorem}

The just stated \emph{Theorem} is of particular meta-mathematical interest for the following main reasons.

Firstly, notice that the properties of $Ker(\mathcal{C})$ are a strong formal analogue of those of the kernel of tilting cotorsion pairs. So, we can notice one more time the significant predictive relevance of such a class as a proxy for well-behaved contexts.

Secondly, our result highlights the importance of $Add$ classes. These are very difficult to understand in the general setting, and, at the current time, their study represents a considerable open research branch. We will, then, involve in a deeper discussion in the next section, which aims at presenting some outlooks of current and future research.
\bigskip

Finally, we conclude this section by reporting an interesting remark that emerges from the proof of the previous \emph{Theorem}.
\smallskip

\begin{corollary}
Each $Ker(\mathcal{C})$-cover of $\mathcal{B}$ is actually an $\mathcal{A}$-cover.
\end{corollary}

\begin{proof}
\emph{(Sketch.)} Fix any $B\in \mathcal{B}$. Consider any $Ker(\mathcal{C})$-cover of $B$, say $\psi$, and take a special $\mathcal{A}$-precover $\phi$ of $B$ (hence also a $Ker(\mathcal{C})$-precover of $\mathcal{B}$; to see this consider the induced LES in cohomology to compute $Ext^1(\mathcal{A},\_)$). Then, by a classical result of Xu (\cite{Xu},1.2.7), we have the decomposition $dom(\phi) \simeq X \oplus Y$ s.t. $(\phi|X):X\rightarrow B$ is a $Ker(\mathcal{C})$-cover of $B$ and $Y\subseteq Ker (\phi)$. Now, $(\phi|X) = h \psi$ for a unique iso $h$, so that the $\mathcal{A}$-precovering property implies that $\psi$ factorises the morphisms $\mathcal{A}\rightarrow M$ via $\mathcal{A} \rightarrow X\oplus Y \overset{pr_1}{\rightarrow} X \simeq_h dom(\psi)$. Finally, it is automatic that any endomorphism $f \in End(dom\psi)$ must be an automorphism, because $\psi$ is a  $Ker(\mathcal{C})$-covering.
\end{proof}

\section{Outlooks for future research}

\subsubsection{Enochs Conjecture for small precovering classes}

The most recent results by Saroch presented in his Habiliation Thesis \cite{Saroch-habilitation} and expanded in \cite{Saroch-small} are part of the aforementioned attempts to understand $Add$ classes. In illustrating them, we will follow the terminology of the author and call \textbf{small precovering classes} those of the form $Add(M)$ (or equiv. $Add(\mathcal{S})$ for some set $\mathcal{S}$, up to writing $M:= \oplus \mathcal{S}$). The main theorem provides for a criterion linking the sufficiency of the perfect decomposition property of the module $K$ s.t. $Add(K) = Ker(\mathcal{C})$ to a covering property of the limit closure of its $Add$ class.

\begin{theorem}
Under some set-theoretical 'local incompactness' principle $(\ast)$, a module $M$ s.t. $Add(M)$ is covering in $\lim Add(M)$ must have perfect decomposition. In particular, $Add(M)$ is closed under $\lim$.
\end{theorem}

The proof is really long and complicated, so we will just sketch the idea, highlighting the role of the set theoretical assumption.
\smallskip

Firstly, we need to introduce the definition of \textbf{local direct summand}, that is a direct sum $\oplus_I N_i \leq N$ s.t. each of its \emph{finite} sub-sums is an actual direct summand, although the whole sum is not.
\smallskip

A classical example to foster intuition is provided by an infinite family of injectives $\{ I_\alpha | \alpha < \kappa \}$ over a non-Noetherian ring, regarded as a subset of its product, namely $\oplus_{\kappa} I_\alpha \leq \prod_\kappa I_\alpha$. Indeed, the sum of any finite subfamily coincides with its product (we are in an abelian category), which is itself actually injective; hence, its injectivity forces the inclusion $\oplus_{j\leq n} I_{\alpha_j} \leq \prod_\kappa I_\alpha$ to be further split, as required. On the other hand, whenever they are not isomorphic, clearly $\oplus_\kappa I_\alpha$ is not a direct summand in $\prod_\kappa I_\alpha$, although it is always a pure submodule (see \cite{Trlifaj&Goebel},2.21).
\smallskip

We can characterize modules admitting a perfect decomposition as those $M$ s.t. every local direct summand in a module from $Add(M)$ is an actual direct summand.
\smallskip

Now, let us formulate the set-theoretical hypothesis. We will need some preliminary definitions.

Let $\lambda$ be a limit ordinal. We say that a subset $C\subseteq \kappa$ is \emph{closed} in $\lambda$ if for each ordinal $\alpha < \kappa$, $sup(C\cap \alpha) = \alpha \neq 0 \implies \alpha \in C$ (i.e. it is an 'attractive' subspace of the ordinal space $\lambda$).

Moreover, $C$ is \emph{unbounded} in $\lambda$ provided that for each $\alpha <\lambda$ there is some $\beta \in C$ s.t. $\beta > \alpha$ (i.e. it is cofinal to $\lambda$).

A closed and unbounded subset of a limit ordinal is called a \textbf{club}, it amounts to an unbounded attractive subset of the corresponding ordinal space.

Consider, now, a cardinal $\kappa$ with uncountable cofinality (to avoid trivialities). A subset $S \subseteq \kappa$ is \textbf{stationary} in case it intersects each club in $\kappa$ in a non-trivial way. Stationarity can be seen as a notion of 'density' (it intersects all unbounded attractive subspaces), and hence of being 'big'.

Finally, we say that $E \subseteq \kappa^+$ is \textbf{non-reflecting} if $E\cap \alpha$ is non-stationary in $\alpha$ for each limit ordinal $\alpha < \kappa^+$. Intuitively, this time we are considering a subspace of the ordinal space that is non-dense when restricted to each proper limit ordinal of our compact ordinal space, and hence it is 'locally non-dense'.

We are now ready to state our set-theoretical assumption:

\begin{center}
    $(\ast)$: There is a proper class of cardinals $\kappa$ s.t. \emph{each} stationary subset $E'\subseteq \kappa^+$ has a non-reflecting stationary subset $E$.
\end{center}

In other words, we are requiring that we can choose arbitrarily large cardinals $\kappa^+$ s.t. each of its big dense (i.e. stationary) subsets $E'$ is locally non-dense; more precisely, we can choose a locally non-dense subset (i.e. non-reflecting) $E\subseteq E'$ that is big dense (i.e. stationary) in $E'$.

In this sense we can call it an 'incompactness' principle: it allows us to have properties that are valid for an arbitrarily big family of objects, but in a locally scattered way within it.

In concluding our digression, notice that $(\ast)$ can be proved to be consistent with ZFC.

Now let's turn our attention to the Theorem.
\bigskip
\begin{proof}
\emph{(Sketch.)} Pick a cardinal $\mu$ s.t. $M$ is $\mu-p.$ And assume by contradiction that $M$ does not admit a perfect decomposition or, equivalently, that there is some local direct summand $K:=\oplus_I N_i \leq N\in Add(M)$ that is not an actual direct summand.

Let's take a $K$ that is minimal for such a property, meaning that $\oplus_J N_i \leq_\oplus N$ whenever $\# J < \# I$. In this way we are releasing the arbitrary constraint on the biggest cardinal which witnesses the 'local character' w.r.t. $\oplus$ of local direct summands (that is $\aleph_0$ in the definition).

By Walker's Lemma, w.l.o.g.\ $N_i$ is $\mu-p.$ for each $i\in I$; moreover, after some work with the theory of locally split epics, w.l.o.g.\ $\aleph_0 \leq \# I \leq \mu$.

Now, if the cardinal $\# I$ is singular, we can non-trivially decompose the underlying set into $I = \coprod_J I_j$ with e.g.\ $\# J = cf (\# I)$; by the minimality of $I$, $K=\bigoplus_J (\oplus_{I_j} N_i)$ has terms which are direct summands in $N$, so that, after the relabelling $N_j':= \oplus_{I_j} N_i'$, $K=\oplus_J N_j '$ is still a local direct summand in $N$ with an indexing set of regular cardinality.

Hence, we can assume to be working with a set of an infinite cardinality $\lambda \leq \mu$ s.t. $K$ is $\mu-p.$.

And again by Walker's Lemma we can assume that also $N$ is $\mu-p.$ (up to enlarging $N$ to make it a direct sum of less than $\mu$ copies of $M$).

We can then carry on an inductive homogenisation procedure based on Eilenberg's trick.

Take some arbitrarily big cardinal $\kappa > \mu$, and set $H:=(K \oplus N)^{\kappa} \in Add(M)$; it is $\kappa-p.$

Define a splitting filtration $\mathcal{F}:= (M_\alpha | \alpha \leq \lambda)$ by $M_{\alpha +1}:= M_\alpha \oplus N_\alpha \oplus H$ and extend it by continuity.

Then, it satisfies the following properties:
\begin{itemize}
    \item By Eilenberg's trick, $M_\lambda = K \oplus H^{(\lambda)} \leq N \oplus H^{(\lambda)} \simeq H$; and the inclusion is non-split by the definition of $K$.
    
    \item $M_\alpha \simeq H$ ($\forall \alpha \leq \lambda$), $M_{\alpha+1}/M_{\alpha} \simeq H$ ($\forall \alpha <\lambda$) again by Eilenberg's trick.
    
    \item By the minimality of $\# I = \lambda$, $M_\alpha \leq_\oplus N\oplus H^{(\lambda)}$ ($\forall \alpha <\lambda$), whilst $M_\lambda$ does not split (as already noted).
\end{itemize}

We aim at extending $\mathcal{F}$ to a new filtration $\mathcal{H}=(M_\alpha | \alpha \leq \kappa^+)$ so as to obtain an arbitrary big module $Z:= \cup M_\alpha$.

Apply the set-theoretic assumption $(\ast)$; it posits the
existence of a proper class of arbitrarily large cardinals $\kappa \geq \mu$ s.t. one can fix a subset $E \subseteq \kappa^+$ that is non-reflecting and stationary. So, up to enlarging our chosen cardinal, we can assume to have such a set $E$.

Moreover, notice that we can assume $E$ to consist of limit ordinal numbers of cofinality $\lambda$ (this will be needed in the successor step of the induction).

Indeed, consider the stationary set $E':=\{\alpha < \kappa^+ | cf(\alpha) = \lambda \}$, and intersect it with the club of limit ordinal numbers in $\kappa^+$; we obtain again a stationary set, that we call $E'$ with a slight abuse of notation. Finally, we can apply our assumption $(\ast)$ to extract the sought stationary set $E$ out of it.

So, chosen $E$, we require our filtration to satisfy the following properties:

\begin{enumerate}
    \item $0<\alpha<\kappa^+ \implies M_\alpha \simeq H \in Add(M)$
    
    \item $\beta <\alpha <\kappa^+ \hspace{.3cm} \wedge \hspace{.3cm} \beta \not \in E \implies M_\alpha = M_\beta \oplus G$ for some copy $G$ of $H$ in $M_\alpha$.
    
    \item $\beta \in E \implies M_\beta \leq_{\not \oplus} M_{\beta+1}$.
\end{enumerate}

$\mathcal{F}$ actually satisfies the stated properties, so we will inductively construct $\mathcal{H}$ by prolonging $\mathcal{F}$. Then, assume to have the $M_\beta$'s ($\forall \beta < \alpha$ for some $\lambda < \alpha <\kappa^+$), and distinguish the classical cases.

\begin{itemize}
    \item $\lambda < \alpha < \kappa^+$ \textsc{limit:} by the non-reflecting property of $E$, $E\cap \alpha$ is non-stationary, so that we can find a club $S\subseteq \alpha$ s.t. $S \cap E = \emptyset$. Define $M_\alpha = \bigcup_{S} M_s$.
    
    \item $\alpha$ \textsc{successor:}
        \begin{itemize}
            \item $\alpha = \beta^+$, $\beta \not \in E$: $M_\alpha := M_\beta \oplus H$ works.
            
            \item $\alpha = \beta^+$, $\beta \in E$: $E$ consists of limit ordinals of cofinality $\lambda$, so $\beta$ occurs in the previous case, and we obtain a club $S \subseteq E^c \cap \beta$ inducing a splitting filtration $(M_\gamma | \gamma \in S)$, that may be assumed to be $\lambda$-indexed (because $cf(S) = cf(\beta) = \lambda$), say $(M_{\gamma_\delta} | \delta < \lambda)$.
            
            Our relabelling provides for a family of isomorphisms $\iota_\delta:M_\delta \overset{\sim}{\rightarrow} M_{\gamma_\delta}$ that in turn induces an isomorphism $\iota:M_\lambda \overset{\sim}{\rightarrow} M_\beta \equiv M_{\gamma_{\lambda}}$ which respects the labelling (i.e. it restricts to each $\iota_\delta$). Then, following a classical construction (see \cite{Beke&Rosicky} for some more context), we can define the new term as the pushout of the iso $\iota$ and of the non-split inclusion $\subseteq:K\otimes H^{(\lambda)} \xhookrightarrow{} N\otimes H^{(\lambda)}, i.e. $ $M_\alpha := PO(\iota,\subseteq)$.
            
            Here we actually need the properties of $\lambda$ to take care of property $2)$.
        \end{itemize}
\end{itemize}

By the properties of $E$ our induction procedure produces an actual filtration. (One should check the continuity, that is where the full properties of $S$ are needed.)

Now, set $Z:= M_{\kappa^+}$, and notice that it clearly is in $\lim Add(M)$, being the directed union of the $\kappa^+$-directed system $\mathcal{H} \subseteq Add(M)$. Furthermore, it is arbitrarily large, since its dimension depends on the arbitrarily chosen $\kappa$.

Moreover, it locally resembles a module in $Add(M)$, meaning that there is a cofinal sub-filtration that lives in $Filt(Add(M))$.

Now, one can show that the canonical epi $\zeta:\oplus_{\kappa^+} M_\alpha \twoheadrightarrow{} Z$ is a special $Add(M)$-precover, and, by some classical results, that it is also split, so that $Z\in Add(M)$.

Hence we can apply again Walker's Lemma, and $Z$ does not decompose into modules of strictly smaller cardinality. This brings to a contradiction.

Indeed, we can write $Z=\oplus_{\gamma < \kappa^+}P_\gamma$ for some $\kappa-p.$ direct summands, and, by taking the $\alpha$-partial direct sums, induce a filtration $\mathcal{H}':=(M_\alpha':= \oplus_{\gamma<\alpha} P_\gamma |\alpha <\kappa ^+)$.

Now, a zig-zag argument with the generating sets of the modules reached by the bounching between $\mathcal{H}$ and $\mathcal{H}'$ shows that $T:=\{ \alpha <\kappa^+ | M_\alpha = M_\alpha ' \}$ is a club. Since $E$ is stationary, $T\cap E$ is non-trivial. Thus, there exists some ordinal $\beta$ s.t. $M_\beta \leq_\oplus Z$, that induces a splitting $M_\beta \leq_\oplus M_{\beta+1}$, against property $3)$.
\end{proof}

Hence, in an extremely simplified way, our Theorem aims at explicitly constructing a special $Add(M)$-precover of an arbitrarily big module in $\lim Add(M)$ that is very pathological from the viewpoint of decomposition, and that is presented by a filtration which fails to be an $Add(M)$-filtration at enough limit ordinals of cofinality $\# I$, where the latter is a regular cardinality for which there is a minimal local direct summand in some module of $Add(M)$.

In particular, $(\ast)$ allows us to assume that such a 'failure set' can be taken to be locally scattered within a big dense subset of our arbitrary cardinal.

Moreover, notice that in our construction we do not really need to be able ex absurdum to produce an arbitrarily long filtration $\mathcal{H}$, or equivalently an arbitrarily big pathological module $Z$ (as remarked by the author while sketching the general idea guiding his proof). Indeed, we just need to produce one that is 'big enough', meaning that it is $\kappa-p.$ for some cardinal $\kappa > \mu$.
\bigskip

\emph{Remark.} It seems intuitively reasonable that it should be possible to remove the additional set-theoretical hypothesis. However, at the moment it is not clear how to approach the issue.

Perhaps, one interesting strategy could be to effectively put into practice Saroch's idea of considering a 'sequence' of big $\kappa$-p. and pathological modules indexed by increasing sequences of sets into the proper class posited by $(\ast)$, which consists of arbitrarily big good cardinals $\kappa$. Then, one could attempt with this a generalization of Erdos Probabilistic Approach, provided that one manages to have some good and well-defined notion of probability on ordinal spaces.
\bigskip

This lengthy discussion is motivated by the deep connections of the just presented Theorem with Enochs' Conjecture.

More explicitly, leveraging on the classical work about modules with perfect decomposition, Saroch proved the following result.

\begin{corollary}{(Enochs' Conjecture for small precovering classes)}
Let $M\in Mod-R$ and suppose that the previous theorem is true; i.e. assume either:
\begin{itemize}
    \item the incompactness assumption $(\ast)$ holds true, or
    
    \item $\exists n<\omega$ s.t. $M$ is a direct sum of $\aleph_n-p.$ modules.
\end{itemize}

Then, TFAE:
\begin{enumerate}
    \item $M$ has a perfect decomposition.
    
    \item $Add(M)$ is closed under direct limits.
    
    \item $Add(M)$ is a covering class of modules.
    
    \item $Add(M)$ is covering for $\lim Add(M)$.
    
    \item For each cardinal $\kappa$, the ring $End_R(M^{(\kappa)})$ is semiregular.
\end{enumerate}
\end{corollary}

$1) \implies 2) \implies 3) \implies 4)$ and $1)\implies 5) \implies 3)$ were already known; however, our theorem allows us to close the circle and prove $4)\implies 1)$. Hence, in particular, under our set theoretical assumptions, $2) \iff 3)$, that amounts to a verification of Enochs Conjecture for small precovering classes.

Notice that the condition is fairly sharp, since it turns out to be also necessary for a c.g.\ $K$, as subsequently proved in the same paper.

In relation to our previous digression on $(\ast)$, Saroch extends that result to direct sums of $\aleph_\omega-g.$ modules, while a full generalization of it would be inherently equivalent to Enochs Conjecture.

Finally, observe that $Add(\mathcal{S})$ is covering for many particular, although interesting, classes $\mathcal{S}$ of modules, so that a generalization that relies on no further set-theoretical assumptions is supported by some 'empirical' evidence (as shown again in \cite{Saroch-applications}), although it still remains an open problem.

\subsubsection{Contramodules}

Another promising field of research is triggered by the notion of contramodules, firstly introduced by Eilenberg and More in their memoir of the AMS \cite{Eilenberg&Moore} establishing foundations of relative homological algebra, and then deeply studied by Bazzoni and Positselski.

Remarkable are, with this respect, the recent characterizations of the class $Add$ of contramodules appearing in \cite{Positselski&Trlifaj}, by Positselski and Trlifaj.

In what follows, we will introduce the notion of contramodules, closely following  \cite{Positselski&Trlifaj}, and completing the exposition with references to the lecture notes \cite{Bazzoni-contramodules} by Bazzoni.
\bigskip

We will be concerned with \emph{left linear} topological rings, i.e. those $S\in Ring$ endowed with a topology for which a base of neighborhoods of $0$ consists of open left ideals.

Its completion is defined to be the projective limit $\hat{S} := \varprojlim_{_{S}I \leq S} S/I$ over its open left ideals; notice that our definition can be proved to reduce to the expected one for valuation domains, in what it guarantees the convergence of all Cauchy sequences in $\hat{S}$ for the topology induced by the limit.

We in turn endow $\hat{S}$ with  the induced limit topology, which is left linear and makes it into a separated and complete topological ring. (See \cite{Bazzoni-contramodules} for more details.)

Now, consider any complete, separated, left linear topological ring $\mathcal{S}$. Denote by $[X]\mathcal{S}$ the finite formal linear combinations of the elements of $X$ with coefficients in $\mathcal{S}$.
\smallskip

Let's define the notion of contramodules; it will allow us to leverage on our topological additional structure to make sense of (not only formal) infinite linear combinations.
\smallskip

Denote by $[[X]]\mathcal{S}$ the limit $\varprojlim_{_{\mathcal{S}}J\leq \mathcal{S}} [X]\mathcal{S}/J$ over the open left ideals. They are properly all the infinite formal linear combinations of elements of $X$ with coefficients that converge to $0$ for the limit topology, i.e.\ those elements which constitute a sequence eventually belonging to each open ideal of $\mathcal{S}$.

Consider the map $T_\mathcal{S} :X\mapsto [[X]]\mathcal{S}$ in $Set$; it is clearly covariantly functorial, and caries a \textbf{monad} structure, as expressed by the following natural maps:
\vspace{-.3cm}

\begin{equation*}
    \begin{split}
    \epsilon_X : & X \longrightarrow [[X]]\mathcal{S} \hspace{4cm} \phi_X:[[[[X]]\mathcal{S}]]\mathcal{S} \longrightarrow [[X]]\mathcal{S} \\
    & x \longmapsto x \cdot 1_{\mathcal{S}} \hspace{5cm} \sum_{\alpha} \big(\sum_{\beta} x_{\beta,\alpha} \cdot  r_\beta^\alpha \big) \cdot r_\alpha \longmapsto \sum_{\beta} x_{\beta,\alpha} \cdot \big( \sum_{\alpha} r_\alpha r^{\alpha}_{\beta} \big)
    \end{split}
\end{equation*}

The first one is the 'point measure' map, since it uniformly embeds $X$, whereas the second map is called the 'opening of parentheses', in what it introduces the associativity of infinite linear combinations, and  pictorially allows us to adopt the conventions of Einstein notation.

Notice that the previous presentation is meaningful by the hypotheses about the topology of $\mathcal{S}$, that characterizes converging sequences as eventually converging (i.e. it is described by the cofinite filter). Then, we are allowed to refer to infinite sums, regarding them as the limit (for the topology of $\mathcal{S}$) of their finite sub-sums.

We are now ready to define contramodules on a left linear, separated and complete ring $\mathcal{S}$.

\begin{definition}
Given a left linear, separated and complete ring $\mathcal{S}$, an $\mathcal{S}$-contramodule is an algebra over the monad $T_\mathcal{S}$.

More explicitly, it is a pair $(C,\pi_C:[[C]]\mathcal{S} \rightarrow C)$ in $Set$ s.t. the following composition diagrams hold:
\vspace{-.4cm}

\begin{equation*}
    \begin{split}
    \textrm{Contraassociativity:} \qquad [[[[C]]\mathcal{S}]]\mathcal{S} & \doublerightarrow{[[\pi_C]]\mathcal{S}}{\phi_C} [[C]]\mathcal{S} \overset{\pi_C}{\rightarrow} C \\
    \textrm{Contraunit:} \hspace{1.9cm}  C & \overset{\epsilon_C}{\rightarrow} [[C]]\mathcal{S} 
    \overset{\pi_C}{\rightarrow} C
    \end{split}
\end{equation*}

The map $\pi_C$ is called \emph{contraaction}, and 'interprets' within $C$ the previously defined formal infinite linear combinations.

Define morphisms of $\mathcal{S}$-contramodules to be those maps of sets that induce  morphism squares in the comma category $(TSet \downarrow Set)$.

We can then consider the category of all $\mathcal{S}$-contramodules, $Contra-\mathcal{S}$.
\end{definition}

\emph{Remark.} Notice that, given any $S\in Ring$ endowed with the discrete topology, modules over the monad of $Set$ $T_S:X\mapsto [X]S$ are exactly right $S$-modules for the obvious contraaction summing up finite linear combinations.
\bigskip

Moreover, each $\mathcal{S}$-contramodule has a 'natural' underlying $\mathcal{S}$-module structure, given by the composition
\vspace{-.2cm}

$$[C]\mathcal{S} \xhookrightarrow{} [[C]]\mathcal{S} \overset{\pi_C}{\rightarrow} C$$

Hence, we obtain a forgetful functor $for: Contra-\mathcal{S} \xhookrightarrow{} Mod-\mathcal{S}$. One can prove that $for$ is exact, faithful and is a right-adjoint to $\triangle$. (See again \cite{Bazzoni-contramodules}, also for the definition of $\triangle$.)
\bigskip

Such properties allow us to describe the projectives in $Contra-\mathcal{S}$. Namely, call \textbf{free right $\mathcal{S}$-contramodules} those induced by the 'opening of parenthesis', that are of the form $([[X]]\mathcal{S},\phi_X)$ for $X\in Set$.

The arbitrary direct sum of free objects is induced by the disjoint union of sets, and the latter, in turn, induces that in $Contra-\mathcal{S}$ (Consider a free presentation of the family of contramodules at stake, and notice that $\oplus$ preserves cokernels).

We can, then, characterize the projective objects in $Contra-\mathcal{S}$ as exactly the direct summands of the free ones.
\bigskip

Moreover, given any $M\in Mod-R$, it is interesting to consider the action of the topological ring $End_R(M)$ endowed with the \emph{finite topology} given by the annihilators of f.g.\ submodules $Ann(f.g.):=\{Ann(F) | F \leq M \hspace{.1cm} \textrm{f.g.}\}$.

$M$ has then a canonical left $End(M)$-module structure that makes it into a discrete module.
\bigskip

We will now define the notion of the contratensor product.

Given a topological ring $\mathcal{S}$ as before, a discrete left $\mathcal{S}$-module $M$ and a right $\mathcal{S}$-contramodule $C$, define:

\begin{equation*}
    \begin{split}
        t_C:[[C]]\mathcal{S} \otimes_\mathbb{Z} M & \rightarrow C \otimes_Z M \\
        (\sum_{c\in C} cs_c) \otimes b & \mapsto \sum_{c\in C} (c\otimes s_c b)
    \end{split}
\end{equation*}

Such a map allows us to 'transfer' infinitely many scalars from the left object to the right one.

Observe that it is well-defined since $Ann(b)\in Ann(f.g.)$ is open, and hence all but a finite number of the scalars $s_c$ annihilate $b$.

\begin{definition}
With reference to the previous notation, the \textbf{contratensor product} $C \odot_\mathcal{S} M$ is the coequalizer $Coeq(t_C,\pi_C \otimes_\mathbb{Z} M)$.
\end{definition}

In other words, it is a universal among those generalizations of $\otimes_\mathbb{Z}$ that allow us to consider an unambiguous notion of 'scalar action' of $\mathcal{S}$, in what we can let it act via the contramodule structure on the left argument as well as via the module scalar rule on the right one.

Therefore, we obtain similar properties to those of $\otimes$, that is constructed as the universal object representing the bifunctor $(\_ \prod \_)$. More precisely, one has the following morphisms in $Ab$:
\begin{itemize}
    \item There is a natural epimorphism in $Ab$: $C \otimes_{\mathcal{S}} M \twoheadrightarrow{} C \odot_\mathcal{S} M$.
    
    \item $[[X]]\mathcal{S} \odot_\mathcal{S} M \simeq [X]M = M^{(X)}$.
    
    \item $\odot_\mathcal{S}$ is right exact and preserves coproducts in both its arguments.
\end{itemize}

Our analogy brings us to state the following adjunction theorem. We incidentally remark the inherent relevance of the tools we have just developed, by mimicking an internalization procedure guided by the aim of meaningfully interpret infinite sums within the language of modules.

\begin{theorem}{(\cite{Positselski&Trlifaj},8.1)}
Let $M\in Mod-R$, and consider the finite topological ring $\mathcal{S}:=(End_R(M),Ann(f.g.))$. The restriction of the adjunction:
\vspace{-.2cm}

\begin{equation*}
    [M,\_]_R:Mod-R \longrightarrow Contra-\mathcal{S} \qquad ; \qquad (\_ \odot_\mathcal{S} M): Contra-\mathcal{S} \longrightarrow Mod-R
\end{equation*}

induces an equivalence of categories $Add(M) \sim (Contra-\mathcal{S})_{proj}$, that characterizes $Add(M)$ in terms of the projective objects of the induced contramodule category.
\end{theorem}

Notice that in the previous result we canonically endow $[M,N]_R$ with the $\mathcal{S}$-contramodule structure induced by the discrete left $\mathcal{S}$-module structure via composition of linear morphisms.
\bigskip

Moreover, we are able to provide a characterization of the limit closure of $Add(M)$.

\begin{theorem}{(\cite{Positselski&Trlifaj},8.2)}
With the previous notation, $\lim Add(M)$ coincides with the class of all the $R$-modules of the form $F \otimes_\mathcal{S} M$, where $F\in \lim (Contra-\mathcal{S})_{proj}$.
\end{theorem}

Now, remark that for $\mathcal{S}$-contramodules $G,C$ of which the first is f.g.\ free, the forgetful functor induces an isomorphism $Contra-\mathcal{S}[G,C] \simeq Mod-\mathcal{S}[G,C]$ between the groups of morphisms (employing Mac Lane's notation), given by the presentability of $\oplus$.

In particular, if we denote classes of c.p.\ objects with lower-case letters, $for:contra-\mathcal{S} \xhookrightarrow{} mod-\mathcal{S}$ is actually an equivalence of categories. Moreover, the above theorem induces the following similar characterization.

\begin{proposition}{(\cite{Positselski&Trlifaj},9.1)}
In the previous notation, $\lim add(M)$ coincides with the class of all $R$-modules of the form $F \odot_\mathcal{S} M$, where $F\in \lim (contra-\mathcal{S})_{proj}$.
\end{proposition}

Moreover, one can prove directly the following simpler characterization via the ordinary tensor product.

\begin{theorem}{(\cite{Positselski&Trlifaj},3.3)}
In the previous notation, $\lim add(M)$ coincides with the class of all right $R$-modules of the form $F \otimes_\mathcal{S} M$ where $F$ is a flat right $\mathcal{S}$-module.
\end{theorem}

We can finally present the main applications of the tools introduced so far, that turn out to provide valuable insights on the following issue.

Given a class  of modules $\mathcal{C}$, we aim at investigating whether the inclusion $\lim add(\mathcal{C}) \subseteq \lim Add(\mathcal{C})$ is in general a strict one or not.

There are several interesting examples of equality, although a general answer is still unknown even for projective modules.

It is therefore remarkable the recently discovered link between our question and the properties of the class of projectives in $Contra-\mathcal{S}$, as highlighted by Positselski and Trlifaj. More explicitly, the previous discussion yields the following result.

\begin{corollary}{(\cite{Positselski&Trlifaj},9.4)}
Provided that $\lim (Contra-\mathcal{S})_{proj} = \lim (contra-\mathcal{S})_{proj}$, it holds $\lim Add(M) = \lim add(M)$.
\end{corollary}

Thus, the notion of contramodules seems to be a promising direction for future research in investigating the properties of $Add$ classes, in the aim of gaining also a better understanding of the issues raised by Enochs' Conjecture.

\end{document}